\documentclass[11pt]{article}%
\usepackage{graphics}
\usepackage{indentfirst}
\usepackage{color}
\usepackage{multirow}
\usepackage[latin1]{inputenc}
\usepackage[T1]{fontenc}
\usepackage[english]{babel}
\usepackage{amsfonts}
\usepackage{amssymb}
\usepackage{graphicx}
\usepackage{dsfont}
\usepackage{epsfig}
\usepackage{ae}
\usepackage{latexsym}
\usepackage{amsmath}
\usepackage{setspace}
\usepackage{amsthm}
\usepackage{appendix}
\usepackage{fontenc}
\usepackage[mathscr]{euscript}
\usepackage{url}
\usepackage[breaklinks]{hyperref}
\usepackage{harvard}%
\providecommand{\U}[1]{\protect\rule{.1in}{.1in}}

\hypersetup{
colorlinks=true,           linkcolor=black,    citecolor=black,    filecolor=black,    urlcolor=black           }
\renewcommand{\tilde}{\widetilde}
\newtheorem{theo}{Theorem}[section]
\newtheorem{pr}{Proposition}[section]
\newtheorem{lem}{Lemma}[section]
\newtheorem{co}{Corollary}[section]

\theoremstyle{remark}
\newtheorem{re}{Remark}[section]
\renewcommand{\geq}{\geqslant}
\renewcommand{\leq}{\leqslant}
\renewcommand{\ge}{\geqslant}
\renewcommand{\le}{\leqslant}
\def\be{\begin{equation}}
\def\ee{\end{equation}}
\citationstyle{agsm}
\begin{document}

\title{On the martingale property in stochastic volatility models based on
time-homogeneous diffusions\thanks{C. Bernard acknowledges support from the
Natural Sciences and Engineering Research Council of Canada. Z. Cui
acknowledges support from the German Academic Exchange Service(DAAD)
scholarship at Summer Academy 2012 on ``Advanced Stochastic Methods to Model
Risk" at Ulm University, Germany, where the paper was presented. D.L. McLeish
acknowledges support from the Natural Sciences and Engineering Research
Council of Canada. The authors are grateful to an anonymous referee and the
associate editor for their careful reading and very helpful suggestions, which
improved the paper. The authors thank Antoine Jacquier and seminar participants Christian Benes,
Peter Carr, Travis Fisher, Olympia Hadjiliadis, Adam Kolkiewicz, Elena
Kosygina, Jay Rosen, David Saunders, Mikhail Urusov, and Jiming Yu for helpful
discussions. The usual disclaimer applies.}}
\author{\textsc{Carole Bernard} \thanks{ C. Bernard is with the department of
Statistics and Actuarial Science at the University of Waterloo, Email
\texttt{\ c3bernar@uwaterloo.ca}. }\ \textsc{Zhenyu Cui} \thanks{Corresponding
author. Z. Cui is with the department of Mathematics at the Brooklyn College
of the City University of New York, Email \texttt{zhenyucui@brooklyn.cuny.edu}%
. }\ \ and\ \textsc{Don McLeish} \thanks{D.L. McLeish is with the department
of Statistics and Actuarial Science at the University of Waterloo, Email
\texttt{dlmcleis@uwaterloo.ca}. }\ }
\date{\today}
\maketitle

\begin{abstract}
Lions and Musiela \citeyear{LM07} give sufficient conditions to verify when a
stochastic exponential of a continuous local martingale is a martingale or a
uniformly integrable martingale. Blei and Engelbert \citeyear{BE09} and
Mijatovi\'c and Urusov \citeyear{MU12PTRF} give necessary and sufficient
conditions in the case of perfect correlation ($\rho=1$). For financial
applications, such as checking the martingale property of the stock price
process in correlated stochastic volatility models, we extend their work to
the arbitrary correlation case ($-1\leq\rho\leq1$). We give a complete
classification of the convergence properties of both perpetual and capped
integral functionals of time-homogeneous diffusions and generalize results in
Mijatovi\'c and Urusov \citeyear{MU12FS} \citeyear{MU12PTRF} with direct
proofs avoiding the use of \textit{separating times} (concept introduced by
Cherny and Urusov \citeyear{CU04} and extensively used in the proofs of
Mijatovi\'c and Urusov \citeyear{MU12PTRF}).

\end{abstract}


\vspace{3mm}

\newpage

\section{Introduction}

There are several recent papers proposing sufficient conditions (Lions and
Musiela \citeyear{LM07}) or necessary and sufficient conditions (Blei and
Engelbert \citeyear{BE09}, Delbaen and Shirakawa \citeyear{DS}, Mijatovi\'c
and Urusov \citeyear{MU12PTRF}, Mijatovi\'c, Novak and Urusov
\citeyear{MNU12}) to verify when the stochastic exponential of a continuous
local martingale is a true martingale or a uniformly integrable(UI)
martingale. A relevant application in finance is to check if the discounted
stock price is a true martingale in a general stochastic volatility model with
arbitrary correlation.

This problem has been extensively studied and dates back from Girsanov
\citeyear{G60} who poses the problem of deciding whether a stochastic
exponential is a true martingale or not. Gikhman and Skorohod \citeyear{GS72},
Liptser and Shiryaev \citeyear{LS72}, Novikov \citeyear{N72} and Kazamaki
\citeyear{K77} provide sufficient conditions for the martingale property of a
stochastic exponential. Novikov's criterion is easy to apply in practical
situations, but it may not always be verified in models in mathematical
finance. In the setting of Brownian motions, refer to Kramkov and Shiryaev
\citeyear{KS98}, Cherny and Shiryaev \citeyear{CS01} and Ruf \citeyear{R12}
for improvements of the criteria of Novikov \citeyear{N72} and Kazamaki
\citeyear{K77}. For affine processes, similar questions are considered by
Kallsen and Shiryaev \citeyear{KS02}, Kallsen and Muhle-Karbe
\citeyear{KMK10}, and Mayerhofer, Muhle-Karbe, and Smirnov \citeyear{MMS11}.
Kotani \citeyear{K06} and Hulley and Platen \citeyear{HP11} obtain necessary
and sufficient conditions for a one-dimensional regular strong Markov
continuous local martingale to be a true martingale. In the strand of
stochastic exponentials based on time-homogeneous diffusions, Engelbert and
Schmidt \citeyear{ES84} provide analytic conditions for the martingale
property, and Stummer \citeyear{S93} gives further analytic conditions when
the diffusion coefficient is the identity. Delbaen and Shirakawa \citeyear{DS}
first provide deterministic criteria to check if a stochastic exponential is a
true martingale under a slightly restrictive assumption requiring certain
functions to be locally bounded on $(0,\infty)$. Mijatovi\'c and Urusov
\citeyear{MU12PTRF} removed the restriction of locally boundedness and extend
their results utilizing a new tool called \textit{separating times} introduced
in Cherny and Urusov \citeyear{CU04}. In the context of stochastic volatility
models, Sin \citeyear{S98}, Andersen and Piterbarg \citeyear{AP07}, and Lions
and Musiela \citeyear{LM07} provide easily verifiable sufficient conditions.
Blanchet and Ruf \citeyear{BR13} describe a method to decide on the martingale
property of a non-negative local martingale based on weak convergence
arguments. Through the study of the classical solutions to the valuation
partial differential equation associated with the stochastic volatility model,
Bayraktar, Kardaras and Xing \citeyear{BKX12} establish a necessary and
sufficient condition when the asset price is a martingale. In the context of
stochastic differential equations(SDE), Doss and Lenglart \citeyear{DL78}
provide a detailed study of their asymptotics and other properties. Ruf
\citeyear{R13b} studies the martingale property of a non-negative local
martingale that is given as a nonanticipative functional of a solution to a
SDE. A recent paper by Karatzas and Ruf \citeyear{KR13} provides the precise
relationship between explosions of one-dimensional stochastic differential
equations and the martingale properties of related stochastic exponentials.
For an overview of stochastic exponentials and the related problem of
martingale properties, refer to Rheinl{\"a}nder \citeyear{R10} and the
references therein.

This paper makes two contributions to the current literature. First, we
provide a complete classification of the convergence or divergence properties
of perpetual and capped integral functionals of time-homogeneous diffusions
based on the local integrability of certain deterministic test functions.
Theorem \ref{int} provides similar necessary and sufficient conditions weaker
than those in Salminen and Yor \citeyear{SY06}, Khoshnevisan, Salminen, and
Yor \citeyear{KSY06}.  Mijatovi\'{c} and Urusov \citeyear{MU12WP} provide a similar result. Theorem \ref{int} permits two absorbing boundaries,
while Engelbert and Tittel \citeyear{ET02} assume that there is exactly one
absorbing boundary. Theorem \ref{xit} concerns the \textit{capped} integral
functional and, to the best of authors' knowledge, is new. We also extend some
results in Mijatovi\'{c} and Urusov (\citeyear*{MU12FS}, \citeyear*{MU12PTRF})
from the case $\rho=1$ to the case $-1\leq\rho\leq1$ (see Proposition \ref{mm}
and Proposition \ref{ui}). Our proofs  do not require the concept of
\textit{separating times} introduced by Cherny and Urusov \citeyear{CU04}. As
examples, we give necessary and sufficient conditions for the (uniformly
integrable) martingale property of the stock price in popular stochastic
volatility models (Hull-White \citeyear{HW87}, (stopped) Heston \citeyear{H},
Sch{ö}bel and Zhu \citeyear{SZ99}, and 3/2 models).

Section \ref{ms2} uses the probabilistic setting and technical tools of Ruf
\citeyear{R12} and Carr, Fisher and Ruf \citeyear{CFR13}. Section \ref{ms3}
provides a complete classification of the convergence or divergence properties
of perpetual and capped integral functionals of time-homogeneous diffusions.
The main result of the paper is given in Section \ref{ms4}: we generalize some
results in Mijatovi\'{c} and Urusov (\citeyear*{MU12FS}, \citeyear*{MU12PTRF})
to the arbitrary correlation case with new direct proofs. Section \ref{ms6}
studies in detail the martingale properties in four popular stochastic
volatility models. Section \ref{ms10} concludes.

\section[Necessary and sufficient conditions]{Necessary and sufficient
conditions for the martingale property \label{ms2}}

\subsection{Probabilistic setup}

\label{pset} Throughout the paper, we fix a time horizon $T\in(0,\infty]$. As
in Carr, Fisher and Ruf \citeyear{CFR13}, we define a stochastic basis by
$(\Omega,\mathcal{F}_{T},\{\mathcal{F}_{t}\}_{t\in\lbrack0,T]},P)$ with a
right-continuous filtration $\{\mathcal{F}_{t}\}_{t\in\lbrack0,T]}.$ \ This
basis is assumed rich enough to support the processes described below and
satisfies the regularity conditions outlined in Appendix A. For any stopping
time $\tau$, we define $\mathcal{F}_{\tau}:=\{A\in\mathcal{F}_{T}\mid
A\cap\{\tau\leq t\}\in\mathcal{F}_{t}\text{ for all }t\in\lbrack0,T]\}$ and
$\mathcal{F}_{\tau-}:=\sigma(\{A\cap\{\tau>t\}\in\mathcal{F}_{T}\mid
A\in\mathcal{F}_{t}\text{\ for some }t\in\lbrack0,T]\cup\mathcal{F}_{0}\})$.
In general, non-negative random variables are permitted to take values in the
set $[0,\infty]$ and stopping times $\tau$ are permitted to take values in the
set $[0,\infty]\cup\mathscr{T}$ for some transfinite time $\mathscr{T}>T$ as
in Appendix A of Carr, Fisher and Ruf \citeyear{CFR13}. In special cases, we
will restrict the range.

For an $\mathcal{F}_{t}$-adapted Brownian motion process $W_{t}$, assume that
$Y$ satisfies the SDE
\begin{equation}
dY_{t}=\mu(Y_{t})dt+\sigma(Y_{t})dW_{t},\quad Y_{0}=x_{0}, \label{y}%
\end{equation}
where $\mu,\sigma:J\rightarrow\mathbb{R}$ are Borel functions, $x_{0}\in J$,
and that $\mu,\sigma$ satisfy the Engelbert-Schmidt condition
\begin{equation}
\forall x\in J,\quad\sigma(x)\neq0,\quad\text{ and }\quad\frac{1}{\sigma
^{2}(\cdot)},\quad\frac{\mu(\cdot)}{\sigma^{2}(\cdot)}\in L_{loc}^{1}(J).
\label{cond1}%
\end{equation}

Here $\ L_{loc}^{1}(J)$ denotes the class of locally integrable functions,
i.e. the functions $J\rightarrow\mathbb{R}$ that are integrable on compact
subsets of the state space, $J=(\ell,r),-\infty\leq\ell<r\leq\infty,$ of the
process $Y=(Y_{t})_{t\in\lbrack0,T]}$. We set $\bar{J}=[\ell,r]$.

The Engelbert-Schmidt condition \eqref{cond1} guarantees that the SDE
\eqref{y} has a unique in law weak solution that possibly exits its state
space $J$ (see Theorem $5.15$, page $341$, Karatzas and Shreve \citeyear{KS91}).
Denote the possible exit time\footnote{Refer to Karatzas and Ruf
\citeyear{KR13} for a detailed study of the distribution of this exit time in
a one-dimensional time-homogeneous diffusion setting.} of $Y$ from its state
space by $\zeta$, i.e. $\zeta=\inf\{u>0,Y_{u}\not \in J\}$, $P$-a.s. which
means that on $\{\zeta=\infty\}$ the trajectories of $Y$ do not exit $J$,
$P$-a.s., and on $\{\zeta<\infty\}$, $\lim_{t\rightarrow\zeta}Y_{t}=r$ or
$\lim_{t\rightarrow\zeta}Y_{t}=\ell$, $P$-a.s.. Observe that $Y$ is defined
such that it stays at its exit point, which means that $\ell$ and $r$ are
absorbing boundaries. The following terminology will be used:
\textquotedblleft$Y$ may exit the state space $J$ at $r$" means $P\left(
\zeta<\infty,\lim\limits_{t\rightarrow\zeta}Y_{t}=r\right)  >0$.

Then we introduce a standard Brownian motion $W^{(2)}$ independent of $(Y,W)$.
Let $Z=(Z_{t})_{t\in\lbrack0,T]}$ denote the (discounted) stock price with
$Z_{0}=1$, and define
\begin{equation}
Z_{t}=\exp\left\{  \rho\int_{0}^{t\wedge\zeta}b(Y_{u})dW_{u}+\sqrt{1-\rho^{2}%
}\int_{0}^{t\wedge\zeta}b(Y_{u})dW_{u}^{(2)}-\frac{1}{2}\int_{0}^{t\wedge
\zeta}b^{2}(Y_{u})du\right\}  ,\quad t\in\lbrack0,\infty), \label{eqz1ex}%
\end{equation}
where $b:J\rightarrow\mathbb{R}$ is a Borel function, and the constant
correlation satisfies $-1\leq\rho\leq1$.

\bigskip

Denote $W_{\cdot}^{(1)}=\rho W_{\cdot}+\sqrt{1-\rho^{2}}W_{\cdot}^{(2)}$, we
have
\begin{equation}
Z_{t}=\exp\left\{  \int_{0}^{t\wedge\zeta}b(Y_{u})dW_{u}^{(1)}-\frac{1}{2}%
\int_{0}^{t\wedge\zeta}b^{2}(Y_{u})du\right\}  ,\quad t\in\lbrack0,\infty),
\label{eqz1}%
\end{equation}
and it is easy to verify that $Z$ and $Y$ satisfy the following system of
SDEs
\begin{align}
dZ_{t}  &  =Z_{t}b(Y_{t})dW_{t}^{(1)},\quad Z_{0}=1,\nonumber\\
dY_{t}  &  =\mu(Y_{t})dt+\sigma(Y_{t})dW_{t},\quad Y_{0}=x_{0}. \label{eq1tc}%
\end{align}

The Borel sigma algebra $\mathcal{B}(\mathbb{R})$ in $\mathbb{R}$ is the
smallest $\sigma$-algebra that contains the open intervals of $\mathbb{R}$. In
what follows, $\lambda(\cdot)$ denotes the Lebesgue measure on $\mathcal{B}%
(\mathbb{R})$. We require that\footnote{Note that this is the same condition
as in Mijatovi\'{c} and Urusov (\citeyear*{MU12FS}, \citeyear*{MU12PTRF}), and
Cherny and Urusov \citeyear{CU06}.} $\lambda(x\in(\ell,r):b^{2}(x)>0)>0$, and
assume the following local integrability condition
\begin{equation}
\forall x\in J,\quad\sigma(x)\neq0,\quad\text{ and }\quad\frac{b^{2}(\cdot
)}{\sigma^{2}(\cdot)}\in L_{loc}^{1}(J). \label{cond2}%
\end{equation}

\begin{re}
In the literature (e.g. Andersen and Piterbarg \citeyear{AP07}), there is a
more general class of stochastic volatility models where the (discounted)
stock price has a non-linear diffusion coefficient in $Z$. For example, a
general model is as follows
\begin{align}
dZ_{t}  &  =Z_{t}^{\alpha}b(Y_{t})\mathds{1}_{t\in\lbrack0,\zeta)}dW_{t}%
^{(1)},\quad Z_{0}=1,\nonumber\\
dY_{t}  &  =\mu(Y_{t})\mathds{1}_{t\in\lbrack0,\zeta)}dt+\sigma(Y_{t}%
)\mathds{1}_{t\in\lbrack0,\zeta)}dW_{t},\quad Y_{0}=x_{0},\nonumber
\end{align}
where $W_{t}^{(1)}$ and $W_{t}$ are standard $\mathcal{F}_{t}$-Brownian
motions, with $\mathbb{E}[dW_{t}^{(1)}dW_{t}]=\rho dt$. $\rho$ is the constant
correlation coefficient and $-1\leq\rho\leq1$. Here $1\leq\alpha\leq2$. The
difficulty of dealing with this model lies mainly in obtaining an explicit
representation of $Z$ in terms of functionals of only $Y$. Thus, in this
paper, we only focus on model \eqref{eq1tc}.

\end{re}

\begin{lem}
(Mijatovi\'{c} and Urusov \citeyear{MU12PTRF})\label{ll1}. Assume conditions
\eqref{cond1} and \eqref{cond2}, and $0<t<\infty.$ Then
\[
\int_{0}^{t}b^{2}(Y_{u})du<\infty\text{ P-a.s. on }\left\{  t<\zeta\right\}
,\quad\label{nine}%
\]

\end{lem}

Fix an arbitrary constant $c\in J$ and introduce the scale function $s(\cdot)$
of the SDE \eqref{y} under $P$
\begin{equation}
s(x):=\int_{c}^{x}\exp\left\{  -\int_{c}^{y}\frac{2\mu}{\sigma^{2}%
}(u)du\right\}  dy,\quad x\in\bar{J}. \label{scale}%
\end{equation}

The following result and its proof can be found in Cherny and Urusov
\citeyear{CU06}, here translated into our notation.

\begin{lem}
\label{cu}(Lemma $5.7$, page 149 of Cherny and Urusov \citeyear{CU06}). Assume
conditions \eqref{cond1} and \eqref{cond2} for the SDE \eqref{y}, and
$s(\ell)=-\infty$, $s(r)=\infty$. Then $\int_{0}^{\infty}b^{2}(Y_{u}%
)du=\infty$, $P$-a.s.
\end{lem}

\subsection{Properties of non-negative continuous local martingales}

In this section, we fix a time horizon $T\in(0,\infty]$, and work under the
canonical probability space $(\Omega,\mathcal{F}_{T},(\mathcal{F}_{t}%
)_{t\in\lbrack0,T]},P)$. This space must be rich enough to support processes
with distributions described below, and the filtration $(\mathcal{F}%
_{t})_{t\in\lbrack0,T]}$ must satisfy additional conditions outlined in
Appendix A. We begin by applying some results from Ruf \citeyear{R12} and
Carr, Fisher and Ruf \citeyear{CFR13} concerning non-negative continuous local
martingales to time-homogeneous diffusions as in \eqref{eqz1}. Ruf
\citeyear{R12} does not specify the form of the continuous local martingale
$(L_{t})_{t\in\lbrack0,T)}$, which is, in our setting
\begin{equation}
L_{t}=\int_{0}^{t\wedge\zeta}b(Y_{u})dW_{u}^{(1)}. \label{Lt}%
\end{equation}

To cast the setting of Ruf \citeyear{R12} into the current notation, the
process in \eqref{eqz1} under $P$ can be rewritten as $Z_{t}=\mathcal{E}%
(L_{t})=\exp\left(  L_{t}-\langle L\rangle_{t}/2\right)  $ where $L_{t}$ in
\eqref{Lt} is a continuous local martingale under $P$.

\begin{lem}
\label{r1} (Lemma 1, Ruf \citeyear{R12}) Assume conditions \eqref{cond1} and
\eqref{cond2} for the SDE \eqref{y}. Under $P$, consider the continuous local
martingale $(L_{t})_{t\in\lbrack0,T]}$ given in \eqref{Lt}, and its quadratic
variation $\langle L\rangle_{t}=\int_{0}^{t\wedge\zeta}b^{2}(Y_{u})du$. For a
predictable positive stopping time $0<\tau\leq\infty$, define $Z_{t}%
=\mathcal{E}(L_{t}),t\in\lbrack0,\tau)$. Then the random variable $Z_{\tau
}:=\lim\limits_{t\uparrow\tau}Z_{t}$ exists, is non-negative and satisfies
\[
\left\{  \int_{0}^{\tau\wedge\zeta}b^{2}(Y_{u})du<\infty\right\}  =\left\{
Z_{\tau}>0\right\}  ,\quad\text{$P$-a.s}.
\]

\end{lem}

As an application of Lemma \ref{r1}, we have the following result.

\begin{co}
\label{z0} Assume\footnote{This is stated without proof after equation $(7)$
on page 4, Mijatovi\'{c} and Urusov \citeyear{MU12PTRF}, and after equation
$(2.4)$ on page 228, Mijatovi\'{c} and Urusov \citeyear{MU12FS}. Here we provide a
proof.} conditions \eqref{cond1} and \eqref{cond2} for the SDE \eqref{y}.
Under $P$, with the process $Z$ defined in \eqref{eqz1}, for $t\in
\lbrack0,T]$
\[
\left\{  Z_{t}=0\right\}  =\left\{  \zeta\leq t,\int_{0}^{\zeta}b^{2}%
(Y_{u})du=\infty\right\}  ,\quad\text{$P$-a.s.}%
\]

\end{co}

\proof From Lemma \ref{r1},
\[
\left\{  Z_{t}=0\right\}  =\left\{  \int_{0}^{t\wedge\zeta}b^{2}%
(Y_{u})du=\infty\right\}  ,\quad\text{$P$-a.s.}%
\]
From Lemma \ref{ll1}, $P\left(  \int_{0}^{t\wedge\zeta}b^{2}(Y_{u}%
)du<\infty\right)  =P\left(  \int_{0}^{t}b^{2}(Y_{u})du<\infty\right)  =1$ on
the set $\left\{  t<\zeta,t\in\lbrack0,T]\right\}  $. Therefore
\[
\left\{  Z_{t}=0\right\}  =\left\{  \zeta\leq t,\int_{0}^{t\wedge\zeta}%
b^{2}(Y_{u})du=\infty\right\}  ,\quad\text{$P$-a.s.}%
\]
\qed

In the following, for notation convenience, denote $T_{\infty}:=R$ and
$T_{0}:=S$ as the first hitting times to $\infty$ and $0$ respectively by $Z$,
where $R$ and $S$ are defined in Section \ref{pset}. Both may take values in
$[0,\infty]\cup\mathscr{T}$. The next result is Theorem 2.1 of Carr, Fisher
and Ruf \citeyear{CFR13} and given in our notation$.$

\begin{pr}
\label{r2} (Theorem 2.1, of Carr, Fisher and Ruf\footnote{Theorem $2.1$, page 6 of
Carr, Fisher and Ruf \citeyear{CFR13} is a general result for non-negative
local martingales. See also Ruf \citeyear{R12} for a similar result for
\textit{continuous} non-negative local martingales.} \citeyear{CFR13}).
Consider the canonical probability space $(\Omega,\mathcal{F}_{T}%
,(\mathcal{F}_{t})_{t\in\lbrack0,T]},P)$, with the process $Z$ defined in
\eqref{eqz1} (so that $Z_{0}=1)$ and assume conditions \eqref{cond1} and
\eqref{cond2}. Then there exists a unique probability measure $\widetilde{P}$
on $(\Omega,\mathcal{F}_{T_{\infty}-})$ such that, for any stopping time
$0<\nu<\infty,$ \newline(1)
\begin{equation}
\tilde{P}(A\cap\lbrace T_{\infty}>\nu\wedge T\rbrace)=\mathbb{E}%
^{P}[\mathds{1}_{A}Z_{\nu\wedge T}] \label{2.1_cond0}%
\end{equation}
for all $A\in\mathcal{F}_{\nu\wedge T}$. \newline(2) for all non-negative
$\mathcal{F}_{\nu\wedge T}$-measurable random variables $U$ taking values in
$[0,\infty],$
\begin{equation}
\mathbb{E}^{\widetilde{P}}\left[  U\mathds{1}_{\lbrace T_{\infty}>\nu\wedge
T\rbrace}\right]  =\mathbb{E}^{P}\left[  U Z_{\nu\wedge T}\mathds{1}_{\lbrace
T_{0}>\nu\wedge T\rbrace}\right]  , \label{2.1_cond1}%
\end{equation}
and, with $\tilde{Z}_{t}=\frac{1}{Z_{t}}\mathds{1}_{\lbrace T_{\infty
}>t\rbrace}$\footnote{By definition this is 0 whenever $t\geq T_{\infty}$ even
if $Z_{t}=0.$},
\begin{equation}
\mathbb{E}^{P}\left[  U\mathds{1}_{\lbrace T_{0}>\nu\wedge T\rbrace}\right]
=\mathbb{E}^{\tilde{P}}\left[  U\tilde{Z}_{\nu\wedge T}\right]  ,
\label{2.1_cond2}%
\end{equation}
(3) $Z$ is a uniformly integrable $P$ martingale on $[0,T]$ if and only if
\begin{equation}
\tilde{P}(T_{\infty}>T)=1. \label{2.1_cond3}%
\end{equation}

\end{pr}

Notice that from (\ref{2.1_cond0}), for any stopping time $\nu<T,$ $\tilde
{P}(Z_{\nu}=0)=0$ so that the measure $\tilde{P}$ assigns zero mass to paths
$Z_{t}$ that hit 0. The condition (\ref{2.1_cond3}) is equivalent to
\begin{align*}
\tilde{P}\left(  \sup_{t\in(0,T]}Z_{t}<\infty\right)   &  =\tilde{P}\left(
\inf_{t\in(0,T]}\tilde{Z}_{t}>0\right)  =1\text{ or } \tilde{P}\left(
\sup_{t\in[0,T]}Z_{t}=\infty\right)  =0.
\end{align*}

\begin{pr}
\label{expZ_}(1) Under $P$, for $t\in\lbrack0,T_{0})$, define the continuous
$P$-local martingale $L_{t}$ as in \eqref{Lt}. Then under $\widetilde{P}$, for
$t\in\lbrack0,T_{\infty})$, $\widetilde{L}_{t}:=L_{t}-\langle L\rangle
_{t}=\int_{0}^{t\wedge\zeta}b(Y_{u})dW_{u}^{(1)}-\int_{0}^{t\wedge\zeta}%
b^{2}(Y_{u})du$ is a continuous $\widetilde{P}$-local martingale.

(2) Under $\widetilde{P}$, for $t\in\lbrack0,T_{\infty})$
\[
\tilde{Z}_{t}=\mathcal{E}(-\widetilde{L}_{t})=\exp\left\{  -\int_{0}%
^{t\wedge\zeta}b(Y_{u})dW_{u}^{(1)}+\frac{1}{2}\int_{0}^{t\wedge\zeta}%
b^{2}(Y_{u})du\right\}  .
\]

\end{pr}

\proof

For statement (1), we need to show that $\widetilde{L}_{t}=\int_{0}%
^{t\wedge\zeta}b(Y_{u})dW_{u}^{(1)}-\int_{0}^{t\wedge\zeta}b^{2}(Y_{u})du$ is
a $\widetilde{P}$-local-martingale on $[0,T_{\infty})$. Recall $R_{n}$ is the
first hitting time of $Z_{t}$ to the level $n$, and put $\tau_{n}=R_{n}\wedge
n$ for all $n\in\mathbb{N}$. We will show that $\widetilde{L}_{t\wedge\tau
_{n}}=\int_{0}^{t\wedge\zeta\wedge\tau_{n}}b(Y_{u})dW_{u}^{(1)}-\int
_{0}^{t\wedge\zeta\wedge\tau_{n}}b^{2}(Y_{u})du$ is a $\widetilde{P}$-local
martingale. This follows from the Girsanov theorem (Ch.VIII, Theorem $1.4$ in
Revuz and Yor \citeyear{RY99}), the facts that $\widetilde{P}<<P$ on
$\mathcal{F}_{\tau_{n}}$ and $\widetilde{P}(\lim_{n\rightarrow\infty}\tau
_{n}=T_{\infty})=1.$

For statement (2), under $\widetilde{P}$, for $t<T_{\infty}$
\begin{align}
\tilde{Z}_{t}  &  =\exp\left\{  -\int_{0}^{t\wedge\zeta}b(Y_{u})dW_{u}%
^{(1)}+\frac{1}{2}\int_{0}^{t\wedge\zeta}b^{2}(Y_{u})du\right\} \nonumber\\
&  =\exp\left\{  -\int_{0}^{t\wedge\zeta}b(Y_{u})dW_{u}^{(1)}+\int
_{0}^{t\wedge\zeta}b^{2}(Y_{u})du-\frac{1}{2}\int_{0}^{t\wedge\zeta}%
b^{2}(Y_{u})du\right\} \nonumber\\
&  =\mathcal{E}(-\widetilde{L}_{t}^{\ast}).\nonumber
\end{align}
\qed

Now we seek to determine the SDE satisfied by $Y$ under $\widetilde{P}$.

\begin{pr}
\label{newsde} Assume conditions \eqref{cond1} and \eqref{cond2} for the SDE
\eqref{y}. Under $\widetilde{P}$, for $-1\leq\rho\leq1$, the diffusion $Y$
satisfies the following SDE up to $\zeta$
\begin{equation}
dY_{t}=(\mu(Y_{t})+\rho b(Y_{t})\sigma(Y_{t}))\mathds{1}_{t\in\lbrack0,\zeta
)}dt+\sigma(Y_{t})\mathds{1}_{t\in\lbrack0,\zeta)}d\widetilde{W}_{t},\quad
Y_{0}=x_{0}. \label{yo}%
\end{equation}

\end{pr}

\proof Consider the system of SDEs in \eqref{eq1tc}, from the Cholesky
decomposition, $dW_{t}^{(1)}=\rho dW_{t}+\sqrt{1-\rho^{2}}dW_{t}^{(2)}$, where
$W$ and $W^{(2)}$ are standard independent Brownian motions under $P$. Define
for $t\in\lbrack0,T]$
\begin{equation}
\widetilde{W}_{t}:=%
\begin{cases}
W_{t}-\rho\int_{0}^{t}b(Y_{u})du,\quad & \text{if $t<\zeta$},\label{keyw}\\
W_{\zeta}-\rho\int_{0}^{\zeta}b(Y_{u})du+\widetilde{\beta}_{t-\zeta},\quad &
\text{if $t\geq\zeta$},
\end{cases}
\end{equation}
where $\widetilde{\beta}$ is a standard $\widetilde{P}$-Brownian motion
independent of $W$ with $\widetilde{\beta}_{0}=0$.

Define $\xi_{n}=\zeta\wedge\tau_{n}$, where $\tau_{n}=R_{n}\wedge n$ and
consider the process $\widetilde{W}$ up to $\xi_{n}$. Since $\mathcal{F}%
_{\xi_{n}}\subset\mathcal{F}_{\tau_{n}}$, it follows from Proposition \ref{r2}
that $\widetilde{P}$ restricted to $\mathcal{F}_{\xi_{n}}$ is absolutely
continuous with respect to $P$ restricted to $\mathcal{F}_{\xi_{n}}$ for
$n\in\mathbb{N}$. Then from Girsanov Theorem (Ch.VIII, Theorem 1.12, page $331$ of
Revuz and Yor \citeyear{RY99})
\begin{align}
\widetilde{W}_{t}  &  :=W_{t}-\langle W_{t},\int_{0}^{t}b(Y_{u})dW_{u}%
^{(1)}\rangle\nonumber\\
&  =W_{t}-\langle W_{t},\rho\int_{0}^{t}b(Y_{u})dW_{u}\rangle-\langle
W_{t},\sqrt{1-\rho^{2}}\int_{0}^{t}b(Y_{u})dW_{u}^{(2)}\rangle\nonumber\\
&  =W_{t}-\rho\int_{0}^{t}b(Y_{u})du,\nonumber
\end{align}
is a $\widetilde{P}$-Brownian motion for $t\in\lbrack0,\xi_{n})$ and
$n\in\mathbb{N}$. It is easy to see from the construction \eqref{keyw}, the
finite dimensional distributions of $\widetilde{W}$ are those of a Brownian
motion under $\widetilde{P}$ on $[0,\xi_{n})$. Thus $Y$ is governed by the
following SDE under $\widetilde{P}$ for $t\in\lbrack0,\xi_{n})$
\begin{align}
dY_{t}  &  =\mu(Y_{t})dt+\sigma(Y_{t})\left(  d\widetilde{W}_{t}+\rho
b(Y_{t})dt\right) \nonumber\\
&  =(\mu(Y_{t})+\rho b(Y_{t})\sigma(Y_{t}))dt+\sigma(Y_{t})d\widetilde{W}%
_{t},\quad Y_{0}=x_{0}.
\end{align}
The result will follow from the following lemma which shows that $\xi
_{n}=\zeta\wedge R_{n}\wedge n\rightarrow\zeta\wedge T_{\infty}=\zeta$,
$P$-a.s. \qed

\begin{lem}
\label{r11} Assume conditions \eqref{cond1} and \eqref{cond2}, then $\zeta\leq
T_{0}\wedge T_{\infty}$, $P$-a.s. and $\widetilde P$-a.s.
\end{lem}

\proof  We prove by contradiction that $P(T_{0}\wedge T_{\infty}<\zeta)=0$.
Suppose that $T_{\infty}<\zeta$ with positive probability so that for some
$t,$ $P(T_{\infty}<t<\zeta)>0.$ Since $T_{\infty}<t,$ $P(Z_{t}=\infty)>0.$ By
Lemma \ref{ll1},
\begin{equation}
P\left(  \int_{0}^{t}b^{2}(Y_{u})du<\infty,Z_{t}=\infty\right)>0.
\label{contr1}%
\end{equation}
Note that $Z_{t}=\exp(L_{t}-\frac{1}{2}\left\langle L\right\rangle
_{t})=\infty$ if and only if $L_{t}=\infty.$ By the Dambis-Dubins-Schwartz
theorem (Ch.V, Theorem $1.6$, Revuz and Yor \citeyear{RY99}), for some
Brownian motion $B$ on an extended probability space, we can write
$L_{t}-\frac{1}{2}\left\langle L\right\rangle _{t}=\left\langle L\right\rangle
_{t}\left(  \frac{B_{\left\langle L\right\rangle _{t}}}{\left\langle
L\right\rangle _{t}}-\frac{1}{2}\right)  $ and from the continuity of the
Brownian motion, $P(\left\langle L\right\rangle _{t}<\infty,L_{t}-\frac{1}%
{2}\left\langle L\right\rangle _{t}=\infty)=P(\left\langle L\right\rangle
_{t}<\infty,B_{\left\langle L\right\rangle _{t}}=\infty)=0$ so that%

\[
P\left(  \int_{0}^{t}b^{2}(Y_{u})du<\infty,Z_{t}=\infty\right)  =P\left(
\left\langle L\right\rangle _{t}<\infty,L_{t}-\frac{1}{2}\left\langle
L\right\rangle _{t}=\infty\right)  =0
\]
contradicting (\ref{contr1}). Similarly suppose that, for some $t,$
$P(T_{0}<t<\zeta)>0.$ Then $P(Z_{t}=0)>0$ and since $t<\zeta,$ from Lemma
\ref{ll1},
\[
P\left(  \int_{0}^{t}b^{2}(Y_{u})du<\infty,Z_{t}=0\right)  >0
\]
contradicting Lemma \ref{r1}. We have thus shown that $P(T_{\infty}%
<\zeta)=P(T_{0}<\zeta)=0.$ To demonstrate a similar statement under the
probability measure $\tilde{P},$ note that $\widetilde{P}$ is a probability
measure on $(\Omega,\mathcal{F}_{R-})$ such that, for a stopping time
$R_{n},$
\[
\tilde{P}( \zeta>T_{0}\wedge R_{n}) =\tilde{P}(\lbrace\zeta>T_{0}\wedge
R_{n}\rbrace\cap\lbrace T_{\infty}>R_{n}\rbrace)=\mathbb{E}^{P}%
[\mathds{1}_{\lbrace\zeta>T_{0}\wedge R_{n}\rbrace}Z_{R_{n}}]=0
\]
since $\lbrace\zeta>T_{0}\wedge R_{n}\rbrace$ is a $\mathcal{F}_{R_{n}}$
measurable event. The last equality holds since $P(\zeta>T_{0}\wedge
R_{n})=0.$ Then by monotone convergence
\[
\tilde{P}( \zeta>T_{0}\wedge T_{\infty}) =\lim_{n\rightarrow\infty}\tilde{P}(
\zeta>T_{0}\wedge R_{n}) =0.
\]
\qed

In view of Lemma \ref{r11} and the definition of $Z_{t}$ in (\ref{eqz1}),
there are only three possibilities almost surely under the measures $P$ and
$\widetilde{P}$:
\[
\zeta=T_{0}<T_{\infty}=\mathscr{T}\text{ \ \ or \ \ }\zeta=T_{\infty}%
<T_{0}=\mathscr{T}\ \ \text{\ \ \ or \ \ \ }\zeta<T_{0}=T_{\infty
}=\mathscr{T}.
\]

In order to verify $\mathbb{E}^{P}[Z_{T}]=1$ for $T\in\lbrack0,\infty]$, the
equivalent condition in Proposition \ref{r2}, (3) can be transformed into a
condition related to integral functionals of $Y$ under $\widetilde{P}$ as
shown in the following proposition.

\begin{pr}
\label{r4} Assume\footnote{A similar result for the general setting of
multi-dimensional diffusions appears in Theorem $1$ of Ruf \citeyear{R13b}.}
conditions \eqref{cond1} and \eqref{cond2}, and $T\in\lbrack0,\infty]$. Then
$Z_{t}$ is a (uniformly integrable) $P$-martingale for $t\in\lbrack0,T]$, i.e.
$\mathbb{E}^{P}[Z_{T}]=1$, if and only if $\widetilde{P}\left(  \int
_{0}^{T\wedge\zeta}b^{2}(Y_{u})du<\infty\right)  =1$.
\end{pr}

\proof\ By Proposition \ref{r2} (3), we have a uniformly integrable martingale
satisfying $\mathbb{E}^{P}[Z_{T}]=1$ if and only if
\[
\tilde{P}\left(  T_{\infty}>T\right)  =\tilde{P}\left(  0<\inf_{t\in
\lbrack0,T]}\tilde{Z}_{t}\right)  =1\text{.}%
\]
But by Proposition \ref{expZ_} (2), under the measure $\tilde{P}$ \
\[
\tilde{Z}_{t}=\mathcal{E}(-\widetilde{L}_{t})=\exp\left\{  -\int_{0}%
^{t\wedge\zeta}b(Y_{u})dW_{u}^{(1)}+\frac{1}{2}\int_{0}^{t\wedge\zeta}%
b^{2}(Y_{u})du\right\}
\]
is a continuous local martingale and for a stopping time $\tau=T\wedge
T_{\infty}\wedge\zeta=T\wedge\zeta,$ by Lemma \ref{r1}, $\lbrace\tilde
{Z}_{\tau}>0\rbrace=\lbrace\int_{0}^{\tau\wedge\zeta}b^{2}(Y_{u}%
)du<\infty\rbrace.$ The result follows. \qed
\begin{re}
Since \ $\widetilde{P}\left(  \int_{0}^{T\wedge\zeta}b^{2}(Y_{u}%
)du<\infty\right)  =\lim_{q\rightarrow0}E^{\tilde{P}}\left[  e^{-q\int
_{0}^{T\wedge\zeta}b^{2}(Y_{u})du}\right]  $  is the right limit of the
Laplace transform $\mathcal{L}(q)$ of $\int_{0}^{T\wedge\zeta}b^{2}(Y_{u})du$
at $0$ \ under the measure $\widetilde{P}$ (and defining $\mathcal{L}(0)=1$), we have the alternative formulation of Proposition \ref{r4} \ that $Z_{t}$
is a (uniformly integrable) $P$-martingale on $[0,T]$, i.e. $\mathbb{E}%
^{P}[Z_{T}]=1$, if and only if $\mathcal{L}(q)$ is right continuous at $0.$
\end{re}

\section[Classification of convergence properties]{Classification of
convergence properties of integral functionals of time-homogeneous diffusions
\label{ms3}}

The Engelbert-Schmidt zero-one law was initially proved in the Brownian motion
case (see Engelbert and Schmidt \citeyear{ES81} or Proposition $3.6.27 $,
page 216 of Karatzas and Shreve \citeyear{KS91}). Engelbert and Tittel
\citeyear{ET02} obtain a generalized Engelbert-Schmidt type zero-one law for
the integral functional $\int_{0}^{t} f(X_{s})ds$, where $f$ is a non-negative
Borel function and $X$ is a strong Markov continuous local martingale. In an
expository paper, Mijatovi\'c and Urusov \citeyear{MU12WP} consider the case
of a one-dimensional time-homogeneous diffusion and the zero-one law is given
in their Theorem $2.11$. They provide two proofs that circumvent the use of
Jeulin's lemma\footnote{The first proof is based on William's theorem
(Ch.VII, Corollary $4.6$, page $317 $, Revuz and Yor \citeyear{RY99}). The second
proof is based on the first Ray-Knight theorem (Ch.XI, Theorem $2.2$, page $455$,
Revuz and Yor \citeyear{RY99}).}. Through stochastic time-change, Cui
\citeyear{C14} proposes a new proof under a slightly stronger assumption.

Recall the scale function $s(\cdot)$ defined in \eqref{scale}, and introduce
the following test functions for $x\in\bar{J}$, with a constant $c\in J$.
\begin{align}
v(x)  &  := \int_{c}^{x} (s(x)-s(y)) \frac{2}{s^{\prime}(y) \sigma^{2}%
(y)}dy,\nonumber\\
v_{b} (x)  &  := \int_{c}^{x} (s(x)-s(y)) \frac{2 b^{2}(y)}{s^{\prime}(y)
\sigma^{2}(y)}dy. \label{vb}%
\end{align}
Note that if $s(\infty)=\infty$, then $v(\infty)=\infty$ and $v_{b}%
(\infty)=\infty$ by the definition in \eqref{vb}. Define $\widetilde{s}%
(\cdot)$, $\widetilde{v}(\cdot)$ and $\widetilde{v}_{b}(\cdot)$ similarly
based on the SDE \eqref{yo} under $\widetilde P$. Throughout this section, we
assume that $\lambda(x\in(\ell,r): b^{2}(x) >0)>0$, which is assumed in
Mijatovi\'c and Urusov \citeyear{MU12WP}.

We have the following Engelbert-Schmidt type zero-one law for the SDE
\eqref{y} under $P$, which is Theorem $2.11$ of Mijatovi\'c and Urusov
\citeyear{MU12WP} with $f(\cdot)=b^{2}(\cdot)$ using our notation.

\begin{pr}
\label{zo2} (Engelbert-Schmidt type zero-one law for a time-homogeneous
diffusion, Theorem $2.11$ of Mijatovi\'c and Urusov \citeyear{MU12WP})

Assume conditions \eqref{cond1}, \eqref{cond2} and $s(r)<\infty$.

(i)If $v_{b}(r)<\infty$, then $\int_{0}^{\zeta} b^{2}(Y_{u})du<\infty$,
$P$-a.s. on $\left\{  \lim_{t\rightarrow\zeta} Y_{t} =r\right\}  $.

(ii)If $v_{b}(r)=\infty$, then $\int_{0}^{\zeta} b^{2}(Y_{u})du=\infty$,
$P$-a.s. on $\left\{  \lim_{t\rightarrow\zeta} Y_{t} =r\right\}  $.
\end{pr}

Analogous results on the set $\{\lim_{t\rightarrow\zeta} Y_{t} =\ell\}$ can be
similarly stated. Clearly the above proposition has a counterpart for the SDE
\eqref{yo} under $\widetilde{P}$ for the end points $r$ and $\ell$.

The following result is Proposition $5.5.22$ on page $345$ of Karatzas and Shreve
\citeyear{KS91} using our notation. It classifies possible exit behaviors of
the process $Y$ at the boundaries of its state space $J$ under $P$.

\begin{pr}
(Proposition $5.5.22$, Karatzas and Shreve \citeyear{KS91})\label{ksd} Assume
condition \eqref{cond1}. Let $Y$ be a weak solution of \eqref{y} in $J $ under
$P$, with nonrandom initial condition $Y_{0} =x_{0} \in J$. Distinguish four cases:

\noindent(a) If $s(\ell)=-\infty$ and $s(r)=\infty$, $P(\zeta=\infty
)=P(\sup_{0\leq t<\infty}Y_{t} =r)=P(\inf_{0\leq t<\infty}Y_{t} =\ell)=1. $

\noindent(b) If $s(\ell)>-\infty$ and $s(r)=\infty$, $P(\lim_{t\rightarrow
\zeta}Y_{t} =\ell)=P(\sup_{0\leq t<\zeta}Y_{t} <r)=1. $

\noindent(c) If $s(\ell)=-\infty$ and $s(r)<\infty$, $P(\lim_{t\rightarrow
\zeta}Y_{t} =r)=P(\inf_{0\leq t<\zeta}Y_{t} >\ell)=1. $

\noindent(d) If $s(\ell)>-\infty$ and $s(r)<\infty$, $P(\lim_{t\rightarrow
\zeta}Y_{t} =\ell)=1-P(\lim_{t\rightarrow\zeta}Y_{t} =r)=\frac{s(r)-s(x_{0}%
)}{s(r)-s(\ell)}. $ Note that $0<\frac{s(r)-s(x_{0})}{s(r)-s(\ell)}<1$.
\end{pr}

Analogous results also hold for the SDE \eqref{yo} under $\widetilde{P}$.

\begin{re}
In the conditions $(b)$, $(c)$ and $(d)$ above, we make no claim concerning
the finiteness of $\zeta$. See Remark $5.5.23$ on page $345$ of Karatzas and
Shreve \citeyear{KS91}. Note that conditions $(b)$ and $(c)$ are consequences
of the expression in condition $(d)$ by letting either $s(r)=\infty$ or
$s(\ell)=-\infty$.
\end{re}

Similar to the statements in Proposition \ref{ksd}, for the study of the
convergence or divergence properties of integral functionals of
time-homogeneous diffusions, we distinguish the following four exhaustive and
disjoint cases under $P$:

\noindent$\bullet$ Case (1): $s(\ell)=-\infty$, $s(r)=\infty$.

\noindent$\bullet$ Case (2): $s(\ell)=-\infty$, $s(r)<\infty$.

\noindent$\bullet$ Case (3): $s(\ell)>-\infty$, $s(r)=\infty$.

\noindent$\bullet$ Case (4): $s(\ell)>-\infty$, $s(r)<\infty$.

Further divide each case above into the following subcases based on the
finiteness of $v_{b}(r)$ and $v_{b}(\ell)$ as defined in \eqref{vb}:

Case (2)(i): $s(\ell)=-\infty$, $s(r)<\infty$, $v_{b}(r)=\infty$.

Case (2)(ii): $s(\ell)=-\infty$, $s(r)<\infty$, $v_{b}(r)<\infty$.

Case (3)(i): $s(\ell)>-\infty$, $s(r)=\infty$, $v_{b}(\ell)=\infty$.

Case (3)(ii): $s(\ell)>-\infty$, $s(r)=\infty$, $v_{b}(\ell)<\infty$.

Case (4)(i): $s(\ell)>-\infty$, $s(r)<\infty$, $v_{b}(r)=\infty$, $v_{b}%
(\ell)=\infty$.

Case (4)(ii): $s(\ell)>-\infty$, $s(r)<\infty$, $v_{b}(r)<\infty$, $v_{b}%
(\ell)=\infty$.

Case (4)(iii): $s(\ell)>-\infty$, $s(r)<\infty$, $v_{b}(r)=\infty$,
$v_{b}(\ell)<\infty$.

Case (4)(iv): $s(\ell)>-\infty$, $s(r)<\infty$, $v_{b}(r)<\infty$, $v_{b}%
(\ell)<\infty$.

\noindent Define
\begin{equation}
\label{starphi}\varphi_{t}:=\int_{0}^{t} b^{2}(Y_{u})du,
\end{equation}
for $t\in[0,\zeta]$. Recall that $b^{2}(\cdot)$ is a non-negative Borel
function, thus $\varphi_{t}$ is a non-decreasing function for $t\in[0,\zeta]$.
Because $\varphi_{t}$ is an integral, it is continuous for $t\in[0,\zeta)$,
and is left continuous at $t=\zeta$.
We now apply the Engelbert-Schmidt type zero-one law under $P$ as in
Proposition \ref{zo2} to determine whether $P(\varphi_{\zeta} <\infty)=1$ or
$P(\varphi_{\zeta} =\infty)=1$ in each of the cases above. We first prove two lemmas.

\begin{lem}
\label{sl1} Assume\footnote{Lemma $5.7$, page 149 of Cherny and Urusov
\citeyear{CU06} is a one-sided version of the current result, namely, ``if
$s(\ell)=\infty$ and $s(r)=\infty$(which implies $v_{b}(\ell)=\infty$ and
$v_{b}(r)=\infty$), then $P(\varphi_{\zeta} =\infty)=1$."} conditions
\eqref{cond1} and \eqref{cond2}, then ``$v_{b}(\ell)=\infty$ and
$v_{b}(r)=\infty$" are necessary and sufficient for $P(\varphi_{\zeta}
=\infty)=1$. \textit{\ }
\end{lem}

\proof
For the sufficiency, assume $v_{b}(r)=\infty$ and $v_{b}(\ell)=\infty$ and
consider the following four distinct cases:

\noindent$\bullet$ Case (1): $s(\ell)=-\infty$, $s(r)=\infty$. From
Proposition \ref{ksd} (a), we have $P(\zeta=\infty)=1$. This, combined with
Lemma \ref{cu} implies $P(\varphi_{\zeta}=\infty)=1$.

\noindent$\bullet$ Case (2): $s(\ell)=-\infty$, $s(r)<\infty$. From
Proposition \ref{ksd} (c), $P(\lim_{t\rightarrow\zeta}Y_{t}=r)=1$. Since
$v_{b}(r)=\infty$, then from Proposition \ref{zo2} $P(\varphi_{\zeta}%
=\infty)=P(\varphi_{\zeta}=\infty,\lim_{t\rightarrow\zeta}Y_{t}=r)$ and from
Proposition \ref{ksd}, $P(\varphi_{\zeta}=\infty,\lim_{t\rightarrow\zeta}%
Y_{t}=r)=P(\lim_{t\rightarrow\zeta}Y_{t}=r)=1$.

\noindent$\bullet$ Case (3): $s(\ell)>-\infty$, $s(r)=\infty$. The proof is
similar to Case (2) above by switching the roles of $\ell$ and $r$, and
applying Proposition \ref{ksd} (b) and Proposition \ref{zo2}.

\noindent$\bullet$ Case (4): $s(\ell)>-\infty$, $s(r)<\infty$. From
Proposition \ref{ksd} (d), $0<p=P(\lim_{t\rightarrow\zeta}Y_{t}=r)<1$. Since
$v_{b}(r)=\infty$ and $v_{b}(\ell)=\infty$, from Proposition \ref{zo2}
\begin{align*}
P(\varphi_{\zeta}  &  =\infty)=P(\varphi_{\zeta}=\infty,\lim_{t\rightarrow
\zeta}Y_{t}=r)+P(\varphi_{\zeta}=\infty,\lim_{t\rightarrow\zeta}Y_{t}=\ell)\\
&  =P(\lim_{t\rightarrow\zeta}Y_{t}=r)+P(\lim_{t\rightarrow\zeta}Y_{t}%
=\ell)=1.
\end{align*}
For the necessity, we only need to prove the contrapositive statement:
\textquotedblleft If at least one of $v_{b}(\ell)$ or $v_{b}(r)$ is finite,
then $P(\varphi_{\zeta}=\infty)<1$." Note that case (a) of Proposition
\ref{ksd} is ruled out here so that we are assured that $P(\lim_{t\rightarrow
\zeta}Y_{t}=r)+P(\lim_{t\rightarrow\zeta}Y_{t}=\ell)=1.$ Without loss of
generality, assume that $v_{b}(\ell)<\infty$, because the case $v_{b}%
(r)<\infty$ can be similarly proved. Then
\begin{align*}
P(\varphi_{\zeta}  &  =\infty)=P(\varphi_{\zeta}=\infty,\lim_{t\rightarrow
\zeta}Y_{t}=\ell)+P(\varphi_{\zeta}=\infty,\lim_{t\rightarrow\zeta}Y_{t}=r)\\
&  =P(\varphi_{\zeta}=\infty,\lim_{t\rightarrow\zeta}Y_{t}=r)\\
&  \leq P(\lim_{t\rightarrow\zeta}Y_{t}=r)
\end{align*}
where the second line follows since from Proposition \ref{zo2}, $P(\varphi
_{\zeta}=\infty,\lim_{t\rightarrow\zeta}Y_{t}=\ell)=0$. There are now two
possibilities for $s(r).$ If $s(r)=\infty,$ since $s(\ell)>-\infty$, we have
$P(\lim_{t\rightarrow\zeta}Y_{t}=r)=0$ from Proposition \ref{ksd} (b).
Alternatively, if $s(r)<\infty,$ since also $s(\ell)>-\infty$ we have from
Proposition \ref{ksd} (d), $0<p=P(\lim_{t\rightarrow\zeta}Y_{t}=r)<1$. In both
cases $P(\lim_{t\rightarrow\zeta}Y_{t}=r)<1$, thus $P(\varphi_{\zeta}%
=\infty)<1$, and the necessity follows. \qed

\begin{lem}
\label{sl2} Assume\footnote{Theorem $2.11$ on page 61 of Mijatovi\'c and Urusov
\citeyear{MU12WP} give a similar result.} conditions
\eqref{cond1} and \eqref{cond2}, and $s(\ell)>-\infty$, $s(r)<\infty$, then
``$v_{b}(\ell)<\infty$ and $v_{b}(r)<\infty$" are necessary and sufficient for
$P(\varphi_{\zeta} <\infty)=1$.
\end{lem}

\proof
With $s(\ell)>-\infty$ and $s(r)<\infty$, denote $p=P(\lim_{t\rightarrow\zeta
}Y_{t} =r)=1-P(\lim_{t\rightarrow\zeta}Y_{t} =\ell)$. From Proposition
\ref{ksd} (d), $0<p<1$.

For the sufficiency, assume that $v_{b}(\ell)<\infty$ and $v_{b}(r)<\infty$
hold. We aim to prove that $P(\varphi_{\zeta} <\infty)=1$ where $\phi_{\zeta
}=\int_{0}^{\zeta}b^{2}(Y_{u})du$ according to its definition \eqref{starphi}.

From Proposition \ref{zo2}, $P(\varphi_{\zeta} <\infty, \lim_{t\rightarrow
\zeta} Y_{t} =r)=P(\lim_{t\rightarrow\zeta} Y_{t} =r)$ and $P(\varphi_{\zeta}
<\infty, \lim_{t\rightarrow\zeta} Y_{t} =\ell)=P(\lim_{t\rightarrow\zeta}
Y_{t} =\ell)$. Then
\begin{align}
P(\varphi_{\zeta} <\infty)  &  =P(\varphi_{\zeta} <\infty, \lim_{t\rightarrow
\zeta} Y_{t} =r) + P(\varphi_{\zeta} <\infty, \lim_{t\rightarrow\zeta} Y_{t}
=\ell)\nonumber\\
&  =P(\lim_{t\rightarrow\zeta} Y_{t} =r)+P(\lim_{t\rightarrow\zeta} Y_{t}
=\ell)=1.\nonumber
\end{align}
For the necessity, we only need to prove the contrapositive argument: ``If at
least one of $v_{b}(\ell)$ and $v_{b}(r)$ is infinite, then $P(\varphi_{\zeta}
<\infty)<1$." Without loss of generality, assume that $v_{b}(r)=\infty$,
because the case $v_{b}(\ell)=\infty$ can be similarly proved. From
Proposition \ref{zo2}, $P(\varphi_{\zeta} <\infty, \lim_{t\rightarrow\zeta}
Y_{t} =r)=0$, and
\begin{align}
P(\varphi_{\zeta} <\infty)  &  =P(\varphi_{\zeta} <\infty, \lim_{t\rightarrow
\zeta} Y_{t} =r) + P(\varphi_{\zeta} <\infty, \lim_{t\rightarrow\zeta} Y_{t}
=\ell)\nonumber\\
&  =P(\varphi_{\zeta} <\infty, \lim_{t\rightarrow\zeta} Y_{t} =\ell
)\nonumber\\
&  \leq P(\lim_{t\rightarrow\zeta} Y_{t} =\ell)<1,\quad
\hbox{from Proposition
\ref{ksd}.}\nonumber
\end{align}
Thus the necessity follows. \qed

We now give a detailed study of the function $\varphi_{t},t\in\lbrack0,\zeta]$
under $P$ using the Engelbert-Schmidt type zero-one law. Theorem \ref{int}
completely characterizes the convergence or divergence property of
$\varphi_{t},t\in\lbrack0,\zeta]$, and several results from the literature are
one-sided versions of it: Theorem \ref{int} (i) is Lemma \ref{ll1}, which is
stated and proved after equation (9) on page 5 of Mijatovi\'{c} and Urusov
\citeyear{MU12PTRF}. Theorem 2 on page 3 of Khoshnevisan, Salminen, and Yor
\citeyear{KSY06} provides\footnote{Salminen and Yor \citeyear{SY06} give
similar conditions for a Brownian motion with drift, and Khoshnevisan,
Salminen, and Yor \citeyear{KSY06} extend it to time-homogeneous diffusions.}
the necessary and sufficient conditions for $P(\varphi_{\zeta}<\infty)=1$,
which corresponds to Theorem \ref{int} (ii). However, they make use of the
stochastic time change and It$\bar{\hbox{o}}$'s lemma in their proof, and thus
need to assume the twice differentiability of a function $g(\cdot)$ defined in
their paper. Our proof is based on Engelbert-Schmidt type zero-one laws of
Mijatovi\'{c} and Urusov \citeyear{MU12WP}, and our weaker assumptions concern
the local integrability of certain deterministic functions. Under these assumptions, Mijatovi\'{c} and Urusov \citeyear{MU12WP} give a result similar to Theorem \ref{int} (ii)
 (in their Theorem
$2.11$). In a parallel paper, Engelbert
and Tittel \citeyear{ET02} consider a strong Markov continuous local
martingale and is broader in scope. As a comparison, their Proposition $3.7$
gives necessary and sufficient conditions for the integral functional to be
convergent or divergent, but assume in Proposition $3.7$  that the process
$X$ has exactly one absorbing point  whereas in our setting and that of
Mijatovi\'{c} and Urusov \citeyear{MU12WP}, it is assumed that the process $Y$
can be absorbed at either boundary $\ell$ or $r$.


\begin{theo}
\label{int}
Under conditions \eqref{cond1} and \eqref{cond2}, the following
properties
for  $\varphi_{t}, t\in[0,\zeta]$ hold:

(i) $\varphi_{t} <\infty$ $P$-a.s. on $\left\{  0\leq t<\zeta\right\}  $.

(ii) $P(\varphi_{\zeta} <\infty)=1$ if and only if at least one of the
following conditions is satisfied:

$\quad$ (a) $v_{b}(r)<\infty$ and $s(\ell)=-\infty$,

$\quad$ (b) $v_{b}(\ell)<\infty$ and $s(r)=\infty$,

$\quad$ (c) $v_{b}(r)<\infty$ and $v_{b}(\ell)<\infty$.

(iii) $P(\varphi_{\zeta} =\infty)=1$ if and only if $v_{b}(r)=\infty$ and
$v_{b}(\ell)=\infty$.
\end{theo}

We summarize the results of Theorem \ref{int} in Table \ref{table1fs}
hereafter. Note that $P(\varphi_{\zeta}<\infty)=P(Z_{\infty}>0)$ always holds
by taking $\tau=\infty$ in Lemma \ref{r1}, and the last two columns in Table
\ref{table1fs} agree.

\proof
Statement (i) follows from Lemma \ref{ll1}. For statement (ii), the detailed
proof for each of the cases in Table \ref{table1fs} is as follows:

\noindent$\bullet$ In Case (1), $s(\ell)=-\infty$ and $s(r)=\infty$ and so
from Lemma \ref{cu}, $P(\varphi_{\zeta}=\infty)=1$.

\noindent$\bullet$ In Case (2), $s(\ell)=-\infty$ and $s(r)<\infty$ and so
from Proposition \ref{ksd}, $P(\lim_{t\rightarrow\zeta}Y_{t}=r)=1.$ There are
two possible subcases. First, in Case (2)(i), $v_{b}(r)<\infty$ and it follows
from Lemma \ref{ll1} that $P(\varphi_{\zeta}=\infty)=1.$ In Case (2)(ii),
since $v_{b}(r)<\infty,$ we have from Lemma \ref{sl1} that $\varphi_{\zeta
}<\infty$ a.s. on the set $\lbrace\lim_{t\rightarrow\zeta}Y_{t}=r\rbrace.$
Moreover, from Proposition \ref{ksd}, $P(\lim_{t\rightarrow\zeta}Y_{t}=r)=1.$
It follows that $P(\varphi_{\zeta}<\infty)=1.$

\noindent$\bullet$ In Case (3), $s(\ell)>-\infty$ and $s(r)=\infty$ and so
from Proposition \ref{ksd}, $P(\lim_{t\rightarrow\zeta}Y_{t}=\ell)=1.$ Again
there are two possible subcases, but they are the reverse of cases in (2);
Case (3)(i) is exactly the reverse of (2)(i) with $\ell$ and $r$ interchanged
and similarly, Case (3)(ii) is exactly the reverse of (2)(ii) so the proofs in
Case (2) suffice.

\noindent$\bullet$ In Case (4): $s(\ell)>-\infty$ and $s(r)<\infty$. Then,
from Proposition \ref{ksd}, $1>p=P(\lim_{t\rightarrow\zeta}Y_{t}%
=r)=1-P(\lim_{t\rightarrow\zeta}Y_{t}=\ell)>0$. For individual subcases, in
Case 4(i), Lemma \ref{sl1} implies $P(\varphi_{\zeta}=\infty)=1.$ In Case
(4)(ii), Proposition \ref{zo2} implies that $P(\varphi_{\zeta}=\infty)<1$ so
that $P(\varphi_{\zeta}<\infty)>0.$ By Lemma \ref{sl2}, we have $P(\varphi
_{\zeta}<\infty)<1.$ Case (4) (iii) is exactly the reverse of (4)(ii) with
$\ell$ and $r$ interchanged so the proof follows using this substitution. And
finally, for Case (4)(iv), $P(\varphi_{\zeta}<\infty)=1$ follows from Lemma
3.2. Therefore, we have the three distinct behaviors for $P(\varphi_{\zeta
}<\infty)$ as outlined in Table \ref{table1fs}. The necessity follows by
examination of Table \ref{table1fs}. \qed

\begin{table}[h]
\begin{center}%
\begin{tabular}
[c]{cc|c|c|c|c|c|c}%
{Case} &  & ${s(\ell)}$ & ${s(r)}$ & ${v}_{b}\mathbf{(\ell)}$ & ${v}_{b}{(r)}$
& $P(\varphi_{\zeta}<\infty)$ & $P(Z_{\infty}>0)$\\\hline
\multirow{1}{*}{(1)} &  & \multicolumn{1}{|c}{$-\infty$} &
\multicolumn{1}{|c}{$\infty$} & \multicolumn{1}{|c}{$\infty$} &
\multicolumn{1}{|c}{$\infty$} & \multicolumn{1}{|c}{$0$} &
\multicolumn{1}{|c}{$0$}\\\hline\hline
\multirow{2}{*}{(2)} & (i) & \multicolumn{1}{|c}{$-\infty$} &
\multicolumn{1}{|c}{$<\infty$} & \multicolumn{1}{|c}{$\infty$} &
\multicolumn{1}{|c}{$\infty$} & \multicolumn{1}{|c}{$0$} &
\multicolumn{1}{|c}{$0$}\\\cline{2-8}
& (ii) & \multicolumn{1}{|c}{$-\infty$} & \multicolumn{1}{|c}{$<\infty$} &
\multicolumn{1}{|c}{$\infty$} & \multicolumn{1}{|c}{$<\infty$} &
\multicolumn{1}{|c}{$1$} & \multicolumn{1}{|c}{$1$}\\\hline\hline
\multirow{2}{*}{(3)} & (i) & \multicolumn{1}{|c}{$>-\infty$} &
\multicolumn{1}{|c}{$\infty$} & \multicolumn{1}{|c}{$\infty$} &
\multicolumn{1}{|c}{$\infty$} & \multicolumn{1}{|c}{$0$} &
\multicolumn{1}{|c}{$0$}\\\cline{2-8}
& (ii) & \multicolumn{1}{|c}{$>-\infty$} & \multicolumn{1}{|c}{$\infty$} &
\multicolumn{1}{|c}{$<\infty$} & \multicolumn{1}{|c}{$\infty$} &
\multicolumn{1}{|c}{$1$} & \multicolumn{1}{|c}{$1$}\\\hline\hline
\multirow{4}{*}{(4)} & (i) & \multicolumn{1}{|c}{$>-\infty$} &
\multicolumn{1}{|c}{$<\infty$} & \multicolumn{1}{|c}{$\infty$} &
\multicolumn{1}{|c}{$\infty$} & \multicolumn{1}{|c}{$0$} &
\multicolumn{1}{|c}{$0$}\\\cline{2-8}
& (ii) & \multicolumn{1}{|c}{$>-\infty$} & \multicolumn{1}{|c}{$<\infty$} &
\multicolumn{1}{|c}{$\infty$} & \multicolumn{1}{|c}{$<\infty$} &
\multicolumn{1}{|c}{$(0,1)^{\ast}$} & \multicolumn{1}{|c}{$(0,1)^{\ast}$%
}\\\cline{2-8}
& (iii) & \multicolumn{1}{|c}{$>-\infty$} & \multicolumn{1}{|c}{$<\infty$} &
\multicolumn{1}{|c}{$<\infty$} & \multicolumn{1}{|c}{$\infty$} &
\multicolumn{1}{|c}{$(0,1)^{\ast}$} & \multicolumn{1}{|c}{$(0,1)^{\ast}$%
}\\\cline{2-8}
& (iv) & \multicolumn{1}{|c}{$>-\infty$} & \multicolumn{1}{|c}{$<\infty$} &
\multicolumn{1}{|c}{$<\infty$} & \multicolumn{1}{|c}{$<\infty$} &
\multicolumn{1}{|c}{$1$} & \multicolumn{1}{|c}{$1$}\\\hline
\end{tabular}
\end{center}
\caption{Table indicating the positivity of the stock price and the finiteness
of $\varphi_{\zeta}$. ($^{*}$ indicates that the probability lies in the open
interval (0,1)).}%
\label{table1fs}%
\end{table}

Similar results as Theorem \ref{int} hold under $\widetilde P$, and the
results are summarized in Table \ref{table2fs}. Note that $\mathbb{E}%
^{P}[Z_{\infty}]=\widetilde{P}(\varphi_{\zeta}<\infty)$ from Proposition
\ref{r4}, and the second-to-last and third-to-last columns in Table
\ref{table2fs} are equal.

\newpage\begin{table}[h]
\begin{center}%
\begin{tabular}
[c]{cc|c|c|c|c|c|c|c}%
{Case} &  & $\widetilde{s}(\ell)$ & $\widetilde{s}(r)$ & $\widetilde{v}%
_{b}(\ell)$ & $\widetilde{v}_{b}(r)$ & $\widetilde{P}(\varphi_{\zeta}<\infty)$
& $\mathbb{E}^{P}(Z_{\infty})$ & UI Mart.\\\hline
(1) &  & \multicolumn{1}{|c}{$-\infty$} & \multicolumn{1}{|c}{$\infty$} &
\multicolumn{1}{|c}{$\infty$} & \multicolumn{1}{|c}{$\infty$} &
\multicolumn{1}{|c}{$0$} & \multicolumn{1}{|c}{$<1$} & \multicolumn{1}{|c}{No}%
\\\hline\hline
\multirow{2}{*}{(2)} & (i) & \multicolumn{1}{|c}{$-\infty$} &
\multicolumn{1}{|c}{$<\infty$} & \multicolumn{1}{|c}{$\infty$} &
\multicolumn{1}{|c}{$\infty$} & \multicolumn{1}{|c}{$0$} &
\multicolumn{1}{|c}{$<1$} & \multicolumn{1}{|c}{No}\\\cline{2-9}
& (ii) & \multicolumn{1}{|c}{$-\infty$} & \multicolumn{1}{|c}{$<\infty$} &
\multicolumn{1}{|c}{$\infty$} & \multicolumn{1}{|c}{$<\infty$} &
\multicolumn{1}{|c}{$1$} & \multicolumn{1}{|c}{$1$} & \multicolumn{1}{|c}{Yes}%
\\\hline\hline
\multirow{2}{*}{(3)} & (i) & \multicolumn{1}{|c}{$>-\infty$} &
\multicolumn{1}{|c}{$\infty$} & \multicolumn{1}{|c}{$\infty$} &
\multicolumn{1}{|c}{$\infty$} & \multicolumn{1}{|c}{$0$} &
\multicolumn{1}{|c}{$<1$} & \multicolumn{1}{|c}{No}\\\cline{2-9}
& (ii) & \multicolumn{1}{|c}{$>-\infty$} & \multicolumn{1}{|c}{$\infty$} &
\multicolumn{1}{|c}{$<\infty$} & \multicolumn{1}{|c}{$\infty$} &
\multicolumn{1}{|c}{$1$} & \multicolumn{1}{|c}{$1$} & \multicolumn{1}{|c}{Yes}%
\\\hline\hline
\multirow{4}{*}{(4)} & (i) & \multicolumn{1}{|c}{$>-\infty$} &
\multicolumn{1}{|c}{$<\infty$} & \multicolumn{1}{|c}{$\infty$} &
\multicolumn{1}{|c}{$\infty$} & \multicolumn{1}{|c}{$0$} &
\multicolumn{1}{|c}{$<1$} & \multicolumn{1}{|c}{No}\\\cline{2-9}
& (ii) & \multicolumn{1}{|c}{$>-\infty$} & \multicolumn{1}{|c}{$<\infty$} &
\multicolumn{1}{|c}{$\infty$} & \multicolumn{1}{|c}{$<\infty$} &
\multicolumn{1}{|c}{$(0,1)^{\ast}$} & \multicolumn{1}{|c}{$<1$} &
\multicolumn{1}{|c}{No}\\\cline{2-9}
& (iii) & \multicolumn{1}{|c}{$>-\infty$} & \multicolumn{1}{|c}{$<\infty$} &
\multicolumn{1}{|c}{$<\infty$} & \multicolumn{1}{|c}{$\infty$} &
\multicolumn{1}{|c}{$(0,1)^{\ast}$} & \multicolumn{1}{|c}{$<1$} &
\multicolumn{1}{|c}{No}\\\cline{2-9}
& (iv) & \multicolumn{1}{|c}{$>-\infty$} & \multicolumn{1}{|c}{$<\infty$} &
\multicolumn{1}{|c}{$<\infty$} & \multicolumn{1}{|c}{$<\infty$} &
\multicolumn{1}{|c}{$1$} & \multicolumn{1}{|c}{$1$} & \multicolumn{1}{|c}{Yes}%
\\\hline
\end{tabular}
\end{center}
\caption{Table indicating $\mathbb{E}^{P}\mathbf{(Z}_{\infty}\mathbf{)}$ and
the uniform integrability of $Z$. (*indicates that the probability lies in the
open interval $(0,1)$)}%
\label{table2fs}%
\end{table}

The following result provides necessary and sufficient conditions for
$P(\varphi_{\zeta\wedge T}<\infty)=1$, for $T\in(0,\infty)$.

\begin{theo}
\label{xit} Assume conditions \eqref{cond1} and \eqref{cond2}. 
 $$P(\varphi_{\zeta\wedge T}<\infty)=P\left(\int_{0}^{\zeta\wedge
T} b^{2}(Y_{u})du<\infty\right)=1$$ for all $T\in(0,\infty)$ if and only if at least one of the following
conditions is satisfied:

(a) $v(\ell)=v(r)=\infty$,

(b) $v_{b}(r)<\infty$ and $v(\ell)=\infty$,

(c) $v_{b}(\ell)<\infty$ and $v(r)=\infty$,

(d) $v_{b}(r)<\infty$ and $v_{b}(\ell)<\infty$.
\end{theo}

\proof
The conditions state that $\left(  \left\{  v(\ell)=\infty\right\}  \text{ or
}\left\{  v_{b}(\ell)<\infty\right\}  \right)  $ and $\left(  \left\{
v(r)=\infty\right\}  \text{ or } \left\{  v_{b}(r)<\infty\right\}  \right)  .$
For a given $T<\infty$, define the events $A_{T}=\left\{  \varphi_{\zeta\wedge
T}<\infty\right\}  ,A=\left\{  \varphi_{\zeta}<\infty\right\}  $ and
$B=\left\{  \zeta<\infty\right\}  .$ Notice that the sets $A_{T}\cap B $ form
a decreasing sequence of sets (as $T\rightarrow\infty$ through a countable
set) so that $\bigcap_{T}(A_{T}\cap B)=A\cap B.$ Therefore,
\begin{equation}
P(A_{T}\cap B)\downarrow P(A\cap B)\text{ as }T\rightarrow\infty. \label{lim1}%
\end{equation}

Moreover, from Theorem \ref{int} (i), for each $T<\infty,$
\begin{equation}
P(A_{T}\cap\overline{B})=P(\overline{B}). \label{PAT}%
\end{equation}

We wish to find necessary and sufficient conditions for $P(A_{T})=1$ for all
$T<\infty.$ In view of (\ref{lim1}) and (\ref{PAT}), this is equivalent to the
condition
\begin{align}
P(A_{T}\cap B)+P(\overline{B})  &  =1\text{ for all }T\text{ or}\nonumber\\
P(A\cap B)+P(\overline{B})  &  =1\text{ or } P(B\cap\overline{A}) =0.
\end{align}

In other words, we seek necessary and sufficient conditions to ensure that
\begin{equation}
P(\zeta<\infty,\varphi_{\zeta}=\infty)=0. \label{show_this}%
\end{equation}

We first show the\textbf{\ sufficiency }of the above conditions. Condition (a)
and Feller's test for explosions implies $P(\zeta<\infty)=0$ and so
(\ref{show_this}) follows. $P(\varphi_{\zeta}=\infty)=0$ is implied in the
cases 2(ii), 3(ii) or 4(iv) of Table \ref{table1fs}. These conditions are
special cases of conditions (b), (c) and (d) as indicated in Table \ref{suff}
below.\newline

\begin{table}[h]
\begin{center}%
\begin{tabular}
[c]{c|c|c}%
{Case} & {implies} & {Cases in Table \ref{table1fs}}\\\hline
(b)\ $v_{b}(r)<\infty\text{ and } v(\ell)=\infty$ & $s(r)<\infty$ &
\text{2(ii),4(ii),4(iv)}\\\hline
(c)\ $v_{b}(\ell)<\infty\text{ and } v(r)=\infty$ & $s(\ell)>-\infty$ &
\text{3(ii),4(iii),4(iv)}\\\hline
(d)\ $v_{b}(r)<\infty\text{ and } v_{b}(\ell)<\infty$ & $s(r)<\infty,
s(\ell)>-\infty$ & \text{4(iv)}\\\hline
\end{tabular}
\end{center}
\caption{Correspondence between conditions (b), (c) and (d) for the
sufficiency case and cases in Table \ref{table1fs}.}%
\label{suff}%
\end{table}

It remains to show (\ref{show_this}) in case 4(ii), i.e. $v(\ell
)=\infty,s(\ell)>-\infty,v_{b}(r)<\infty,s(r)<\infty$ and in case 4(iii), i.e.
$v_{b}(\ell)<\infty$, $s(\ell)>-\infty,v(r)=\infty,s(r)<\infty.$ By
interchanging the role of $\ell$ and $r,$ it suffices to show the first of
these. By Proposition \ref{zo2}, $\varphi_{\zeta}<\infty$ $P-$a.s. on the set
$\left\{  \lim\limits_{t\rightarrow\zeta}Y_{t}=r\right\}  $ or
\[
P(\varphi_{\zeta}=\infty,\lim_{t\rightarrow\zeta}Y_{t}=r)=0.
\]
From Feller's test of explosions, $v(\ell)<\infty$ if and only if
$P(\zeta<\infty,\lim_{t\rightarrow\zeta}Y_{t}=\ell)>0$, and so in this case
$v(\ell)=\infty$ implies
\[
P(\zeta<\infty,\lim_{t\rightarrow\zeta}Y_{t}=\ell)=0.
\]
It follows that%
\begin{align}
P(\zeta<\infty,\varphi_{\zeta}=\infty)  &  =P(\zeta<\infty,\varphi_{\zeta
}=\infty,\lim_{t\rightarrow\zeta}Y_{t}=\ell)+P(\zeta<\infty,\varphi_{\zeta
}=\infty,\lim_{t\rightarrow\zeta}Y_{t}=r)\nonumber\\
&  \leq P(\zeta<\infty,\lim_{t\rightarrow\zeta}Y_{t}=\ell)+P(\varphi_{\zeta
}=\infty,\lim_{t\rightarrow\zeta}Y_{t}=r)=0.
\end{align}

For the \textbf{necessity}, we wish to show the contrapositive: if $\left\{
v(\ell)<\infty\text{ and } v_{b}(\ell)=\infty\right\}  $ OR $\left\{
v(r)<\infty\text{ and } v_{b}(r)=\infty\right\}  $ (i.e. at least one of the
two boundaries, $v$ is finite and $v_{b}$ infinite), then (\ref{show_this})
fails, that is
\[
P(\zeta<\infty,\varphi_{\zeta}=\infty)>0.
\]

The contrapositive is consistent with Table \ref{table1fs}, cases 2(i), 3(i),
4(i), 4(ii), 4(iii) as indicated in Table \ref{contra} below.\newline

\begin{table}[h]
\begin{center}%
\begin{tabular}
[c]{c|c|c|c}%
{Contrapositive Case} & {implies} &  & {Cases in Table \ref{table1fs}}\\\hline
$\left\{  v(\ell)<\infty\right\}  \text{ and }\left\{  v_{b}(\ell
)=\infty\right\}  $ & $s(\ell)>-\infty$ & \text{Consistent with} &
\text{3(i),4(i),4(ii)}\\\hline
$\left\{  v(r)<\infty\right\}  \text{ and }\left\{  v_{b}(r)=\infty\right\}  $
& $s(r)<\infty$ & \text{Consistent with} & \text{2(i),4(i),4(iii)}\\\hline
\end{tabular}
\end{center}
\caption{Correspondence between the two conditions from the contrapositive
case and cases in Table \ref{table1fs}.}%
\label{contra}%
\end{table}

Consider the first row above when $v(\ell)<\infty$, $v_{b}(\ell)=\infty$,
$s(\ell)>-\infty.$ By Feller's test, $v(\ell)<\infty$ implies $P(\zeta
<\infty,\lim_{t\rightarrow\zeta}Y_{t}=\ell)>0$ and by Proposition \ref{zo2},
since $v_{b}(\ell)=\infty,$ $\varphi_{\zeta}=\infty$ $P-$a.s. on the set
$\left\{  \lim_{t\rightarrow\zeta}Y_{t}=\ell\right\}  $ and
\[
P(\zeta<\infty,\varphi_{\zeta}=\infty)\geq P(\zeta<\infty,\varphi_{\zeta
}=\infty,\lim_{t\rightarrow\zeta}Y_{t}=\ell)=P(\zeta<\infty,\lim
_{t\rightarrow\zeta}Y_{t}=\ell)>0.
\]
The proof in the second case $v(r)<\infty,$ $v_{b}(r)=\infty$ follows once
again by interchanging the roles of $\ell$ and $r.$ \qed

Similarly statements as Theorem \ref{xit} hold under $\widetilde{P}$ with SDE \eqref{yo}.

\section{Generalization of some results in Mijatovi\'c and Urusov \label{ms4}}

In this section, we generalize the main results in Mijatovi\'c and Urusov
(\citeyear*{MU12FS}, \citeyear*{MU12PTRF}) and provide new unified proofs
without the concepts of ``separating times". Note that Mijatovi\'c and Urusov
(\citeyear*{MU12FS}, \citeyear*{MU12PTRF}) work in the $\rho= 1$ case, and we
generalize it to the arbitrary correlation case.

Consider the stochastic exponential $Z$ defined in \eqref{eqz1}. The following
proposition provides the necessary and sufficient condition for $Z_{T}$ to be
a $P$-martingale for all $T\in(0,\infty)$, when $-1\leq\rho\leq1$. Note that
Theorem $2.1$ in Mijatovi\'c and Urusov \citeyear{MU12PTRF} is the case
$\rho=1$ of the following proposition.

\begin{pr}
\label{mm} Assume conditions \eqref{cond1} and \eqref{cond2}, then for all
$T\in(0,\infty)$, $\mathbb{E}^{P}[Z_{T}]=1$ if and only if at least one of the
conditions (1)-(4) below is satisfied:

(1) $\widetilde{v}(\ell)=\widetilde{v}(r)=\infty$,

(2) $\widetilde{v}_{b}(r)<\infty$ and $\widetilde{v}(\ell)=\infty$,

(3) $\widetilde{v}_{b} (\ell)<\infty$ and $\widetilde{v}(r)=\infty$,

(4) $\widetilde{v}_{b}(r)<\infty$ and $\widetilde{v}_{b} (\ell)<\infty$.
\end{pr}

\proof
From Proposition \ref{r4}, for all $T\in(0,\infty)$, $\mathbb{E}^{P}[Z_{T}]=1$
if and only if $\widetilde{P}(\int_{0}^{\zeta\wedge T} b^{2}(Y_{u})du
<\infty)=1$. Then the statement follows from Theorem \ref{xit} applied to
$\widetilde{P} $. \qed

We have the following necessary and sufficient condition for $Z$ to be a
uniformly integrable $P$-martingale on $[0,\infty]$, when $-1\leq\rho\leq1 $.
Note that Theorem $2.3$ of Mijatovi\'c and Urusov \citeyear{MU12PTRF} proves
the case $\rho=1$ of the following proposition.

\begin{pr}
\label{ui} Assume conditions \eqref{cond1} and \eqref{cond2}, then
$\mathbb{E}^{P}[Z_{\infty}]=1$ if and only if at least one of the conditions
$(A^{\prime})-(D^{\prime})$ below is satisfied:

$(A^{\prime})$ $b=0$ a.e. on $J$ with respect to the Lebesgue measure,

$(B^{\prime})$ $\widetilde{v}_{b}(r)<\infty$ and $\widetilde{s}(\ell)=-\infty$,

$(C^{\prime})$ $\widetilde{v}_{b} (\ell)<\infty$ and $\widetilde{s}(r)=\infty$,

$(D^{\prime})$ $\widetilde{v}_{b}(r)<\infty$ and $\widetilde{v}_{b}
(\ell)<\infty$.
\end{pr}

\proof
From Proposition \ref{r4}, $\mathbb{E}^{P}[Z_{\infty}]=1$ if and only if
$\widetilde{P}(\int_{0}^{\zeta} b^{2}(Y_{u})du <\infty)=1$. Condition
$(A^{\prime})$ is a trivial case and it is easy to verify. From Theorem
\ref{int} applied to $\widetilde{P}$ and the classification in Table
\ref{table2fs}, $\mathbb{E}^{P}[Z_{\infty}]=1$ if and only if at least one of
the conditions $(B^{\prime})$, $(C^{\prime})$ or $(D^{\prime})$ holds. \qed

Here we generalize some results in Mijatovi\'c and Urusov \citeyear{MU12FS} to
the arbitrary correlation case and provide new proofs without the concept of
\textit{separating times}. Precisely, Theorem $2.1$ of Mijatovi\'c and Urusov
\citeyear{MU12FS} is the case $\rho=1$ of the following proposition.

\begin{pr}
\label{mfs} Assume conditions \eqref{cond1} and \eqref{cond2}, then for all
$T\in(0,\infty)$, $Z_{T}>0$ $P$-a.s. if and only if at least one of the
conditions\footnote{Note that conditions (1)-(4) in Proposition \ref{mfs} do
not depend on the correlation $\rho$, which means that the positivity of the
(discounted) stock price does not depend on the correlation. Similar remarks
hold for Proposition \ref{t4fs} and Proposition \ref{t5fs}.} (1)-(4) below is satisfied:

(1) $v(\ell)=v(r)=\infty$,

(2) $v_{b}(r)<\infty$ and $v(\ell)=\infty$,

(3) $v_{b} (\ell)<\infty$ and $v(r)=\infty$,

(4) $v_{b}(r)<\infty$ and $v_{b} (\ell)<\infty$.
\end{pr}

\proof
From Lemma \ref{r1}, for all $T\in(0,\infty)$, $Z_{T}>0$, $P$-a.s. if and only
if $P\left(  \int_{0}^{\zeta\wedge T}b^{2}(Y_{u})du<\infty\right)  =1$. Then
the statement follows from Theorem \ref{xit}. \qed

Note that Theorem $2.3$ of Mijatovi\'c and Urusov \citeyear{MU12FS} proves the
case $\rho=1$ of the following proposition.

\begin{pr}
\label{t4fs} Let the functions $\mu$, $\sigma$ and $b$ satisfy conditions
$(2.1)$, $(2.3)$ and $(2.5)$ of Mijatovi\'c and Urusov \citeyear{MU12FS}
(equivalently conditions \eqref{cond1} and \eqref{cond2} in this paper), and
let $Y$ be a (possibly explosive) solution of the SDE \eqref{y} under $P$,
with $Z$ defined in \eqref{eqz1}, then $Z_{\infty} >0$, $P$-a.s. if and only
if at least one of the conditions (I)-(IV) below is satisfied:

(I) $b=0$ a.e. on $J$ with respect to the Lebesgue measure,

(II) $v_{b}(r)<\infty$ and $s(\ell)=-\infty$,

(III) $v_{b}(\ell)<\infty$ and $s(r)=\infty$,

(IV) $v_{b}(r)<\infty$ and $v_{b}(\ell)<\infty$.
\end{pr}

\proof
Condition (I) is a trivial case and it is easy to verify. From Lemma \ref{r1},
$Z_{\infty} >0$, $P$-a.s. if and only if $P\left(  \int_{0}^{\zeta}b^{2}%
(Y_{s})ds<\infty\right)  =1$. Then the proof follows from Theorem \ref{int}
and the classification in Table \ref{table1fs}. \qed

Note that Theorem $2.5$ of Mijatovi\'c and Urusov \citeyear{MU12FS} is a
special case of the following proposition when $\rho=1$.

\begin{pr}
\label{t5fs} Let the functions $\mu$, $\sigma$ and $b$ satisfy conditions
$(2.1)$, $(2.3)$ and $(2.5)$ of Mijatovi\'c and Urusov \citeyear{MU12FS}
(equivalently conditions \eqref{cond1} and \eqref{cond2} in this paper), and
let $Y$ be a (possibly explosive) solution of the SDE \eqref{y} under $P$,
with $Z$ defined in \eqref{eqz1}. Then $Z_{\infty} =0$, $P$-a.s. if and only
if both conditions (i) and (ii) below are satisfied:

(i) $b$ is not identically zero with respect to Lebesgue measure on $(\ell,r)
$,

(ii) $v_{b}(\ell)=v_{b}(r)=\infty$.
\end{pr}

\proof
Condition (i) is a trivial case and it is easy to verify. From Lemma \ref{r1},
$Z_{\infty} =0$, $P$-a.s. if and only if $P\left(  \int_{0}^{\zeta}b^{2}%
(Y_{u})du =\infty\right)  =P(\varphi_{\zeta}=\infty)=1$. From Theorem
\ref{int} (iii), this is equivalent to checking the condition (ii) here.
\qed

\section{Examples of correlated stochastic volatility models \label{ms6}}

In this section, we apply the results in Section \ref{ms4} to the study of
martingale properties of (discounted) stock prices\footnote{Equivalently, we
may assume that the risk-free interest rate is zero.} in four popular
correlated stochastic volatility models: the (stopped) Heston\footnote{The
volatility is stopped whenever it hits the boundary $0$. When $2\kappa\theta>
\xi^{2}$(zero is unattainable), our model coincides with the usual Heston
model.}, the $3/2$, the Sch{\"o}bel-Zhu and the Hull-White models. The results
are summarized at the end of the section in Table \ref{SUMM} and Table
\ref{SUMM2}.

\subsection{Stopped Heston stochastic volatility model}

Suppose that under a probability measure $P$, the (correlated) stopped Heston
stochastic volatility model has the following diffusive dynamics
\begin{align}
dS_{t}  &  =S_{t} \sqrt{Y_{t}}\mathds{1}_{t\in[0,\zeta)} dW_{t}^{(1)}, \quad
S_{0}=1.\nonumber\\
dY_{t}  &  = \kappa(\theta-Y_{t})\mathds{1}_{t\in[0,\zeta)} dt +\xi\sqrt
{Y_{t}}\mathds{1}_{t\in[0,\zeta)} dW_{t}, \quad Y_{0}=x_{0}>0, \label{hes2b}%
\end{align}
with $\mathbb{E}^{P}[dW_{t}^{(1)} dW_{t}]=\rho dt$, $-1\leq\rho\leq1$,
$\kappa>0, \theta>0, \xi>0$. The natural state space for $Y$ is $J=(\ell
,r)=(0,\infty)$. $\zeta$ is the possible exit time of the process $Y$ from its
state space $J$. The model \eqref{hes2b} belongs to the general stochastic
volatility model considered in \eqref{eq1tc} with $\mu(x)=\kappa(\theta-x)$,
$\sigma(x)=\xi\sqrt{x}$, and $b(x)= \sqrt{x}$. Clearly $\sigma(x)=\xi\sqrt{x}
\neq0, x\in J$, $\frac{1}{\sigma^{2}(x)}=\frac{1}{\xi^{2} x}\in L^{1}%
_{loc}(J)$, $\frac{\mu(x)}{\sigma^{2}(x)}=\frac{\kappa(\theta-x)}{\xi^{2}
x}\in L^{1}_{loc}(J)$, and $\frac{b^{2}(x)}{\sigma^{2}(x)}=\frac{1}{\xi^{2}%
}\in L^{1}_{loc}(J)$ are satisfied. Thus, the conditions \eqref{cond1} and
\eqref{cond2} are satisfied. From Proposition \ref{newsde}, under
$\widetilde{P}$, the diffusion $Y$ satisfies the following SDE
\begin{align}
dY_{t}  &  =\widetilde{\kappa}(\widetilde{\theta}-Y_{t})\mathds{1}_{t\in
[0,\zeta)} dt +\xi\sqrt{Y_{t}}\mathds{1}_{t\in[0,\zeta)} d\widetilde{W}_{t},
\quad Y_{0}=x_{0}>0,\nonumber
\end{align}
where $\widetilde{\kappa}=\kappa-\rho\xi$ and $\widetilde{\theta}=\frac
{\kappa\theta}{\kappa-\rho\xi}$.

For a constant $c\in J$, the scale functions of the SDE \eqref{y} and the SDE
\eqref{yo} are respectively
\begin{align}
s(x)  &  =e^{\frac{2\kappa c}{\xi^{2}}}c^{\frac{2\kappa\theta}{\xi^{2}}}
\int_{c}^{x} y^{-\frac{2\kappa\theta}{\xi^{2}}} e^{\frac{2\kappa y}{\xi^{2}}%
}dy=C_{1} \int_{c}^{x} y^{-\alpha} e^{\beta y}dy,\nonumber\\
\widetilde{s}(x)  &  =e^{\frac{2\widetilde{\kappa} c}{\xi^{2}}}c^{\frac
{2\widetilde{\kappa}\widetilde{\theta}}{\xi^{2}}}\int_{c}^{x} y^{-\frac
{2\widetilde{\kappa}\widetilde{\theta}}{\xi^{2}}} e^{\frac{2\widetilde{\kappa}
y}{\xi^{2}}}dy =C_{2} \int_{c}^{x} y^{-\alpha} e^{\gamma y}dy, \label{sches}%
\end{align}
with $\alpha=\frac{2\kappa\theta}{\xi^{2}}$, $\beta=\frac{2\kappa}{\xi^{2}}%
>0$, $\gamma=\frac{2\kappa}{\xi^{2}}-\frac{2\rho}{\xi}$, and the constant
terms are $C_{1} =e^{\frac{2\kappa c}{\xi^{2}}}c^{\frac{2\kappa\theta}{\xi
^{2}}}>0$ and $C_{2} =e^{\frac{2\kappa c}{\xi^{2}}-\frac{2\rho c}{\xi}%
}c^{\frac{2\kappa\theta}{\xi^{2}}}>0$. Under $\widetilde{P}$, we have the
following test functions for $x\in\bar{J}$,
\[
\widetilde{v}(x) =\frac{2}{\xi^{2}}\int_{c}^{x} \frac{ \int_{y}^{x}
z^{-\alpha}e^{\gamma z}dz }{y^{1-\alpha}e^{\gamma y}}dy, \quad\widetilde
{v}_{b}(x) =\frac{2}{\xi^{2}}\int_{c}^{x} \frac{ \int_{y}^{x} z^{-\alpha
}e^{\gamma z}dz }{y^{-\alpha}e^{\gamma y}}dy.
\]

\begin{pr}
\label{hm1} For\footnote{Proposition \ref{hm1} is consistent with Proposition
$2.5$, page $34$ of Andersen and Piterbarg \citeyear{AP07}, also see Remark $4.2$,
page $2052$ of Del Ba{ñ}o Rollin et al. \citeyear{DBR10}.} the stopped Heston
model \eqref{hes2b}, the underlying stock price $(S_{t})_{0\leq t<\infty}$ is
a true $P$-martingale.
\end{pr}

\proof The proof of Proposition \ref{hm1} is elementary and details are given in Appendix \ref{hm1app}.
 To prove it, we check the conditions of Proposition
\ref{mm}: the results are summarized in Table \ref{tablehes}.\newline%
\begin{table}[h]
\begin{center}%
\begin{tabular}
[c]{c|c|c|c|c}%
Case & $\widetilde v(\ell)$ & $\widetilde v(r)$ & $\widetilde v_{b}(\ell)$ &
$\widetilde v_{b}(r)$\\\hline
$\alpha\geq1$ & $\infty$ & $\infty$ & $\infty$ & $\infty$\\\hline
$\alpha< 1$ & $<\infty$ & $\infty$ & $<\infty$ & $\infty$\\\hline
\end{tabular}
\end{center}
\caption{First classification table for the Heston model}%
\label{tablehes}%
\end{table}

\noindent From Table \ref{tablehes} and Proposition \ref{mm}, $(S_{t})_{0\leq
t<\infty})$ is a true $P$-martingale.\hfill$\Box$

\begin{pr}
\label{uihes} For the stopped Heston model \eqref{hes2b}, the underlying stock
price $(S_{t})_{0\leq t\leq\infty}$ is a uniformly integrable $P$-martingale
if and only if $\rho\xi\leq\kappa<\frac{\xi^{2}}{2\theta}$.
\end{pr}

Note that the Feller condition has to be violated in order to have a UI martingale.

\proof
The proof of Proposition \ref{uihes}  is elementary and details are given  in Appendix \ref{hm2app}.
 To prove it, we check the conditions of Proposition \ref{ui}:
the results are summarized in Table \ref{tablehes2}.

\begin{table}[h]
\begin{center}%
\begin{tabular}
[c]{cc|c|c|c|c|c|c}%
Case &  & $\widetilde s(\ell)$ & $\widetilde s(r)$ & $\widetilde v(\ell)$ &
$\widetilde v(r)$ & $\widetilde v_{b}(\ell)$ & $\widetilde v_{b}(r)$\\\hline
\multirow{3}{*}{ $\alpha> 1$} & $\gamma<0$ & $-\infty$ & $<\infty$ & $\infty$
& $\infty$ & $\infty$ & $\infty$\\\cline{2-8}
& $\gamma=0$ & $-\infty$ & $<\infty$ & $\infty$ & $\infty$ & $\infty$ &
$\infty$\\\cline{2-8}
& $\gamma>0$ & $-\infty$ & $\infty$ & $\infty$ & $\infty$ & $\infty$ &
$\infty$\\\hline\hline
\multirow{3}{*}{ $\alpha= 1$} & $\gamma<0$ & $-\infty$ & $<\infty$ & $\infty$
& $\infty$ & $\infty$ & $\infty$\\\cline{2-8}
& $\gamma=0$ & $-\infty$ & $\infty$ & $\infty$ & $\infty$ & $\infty$ &
$\infty$\\\cline{2-8}
& $\gamma>0$ & $-\infty$ & $\infty$ & $\infty$ & $\infty$ & $\infty$ &
$\infty$\\\hline\hline
\multirow{3}{*}{ $\alpha<1$} & $\gamma<0$ & $>-\infty$ & $<\infty$ & $<\infty$
& $\infty$ & $<\infty$ & $\infty$\\\cline{2-8}
& $\gamma=0$ & $>-\infty$ & $\infty$ & $<\infty$ & $\infty$ & $<\infty$ &
$\infty$\\\cline{2-8}
& $\gamma>0$ & $>-\infty$ & $\infty$ & $<\infty$ & $\infty$ & $<\infty$ &
$\infty$\\\hline
\end{tabular}
\end{center}
\caption{Second classification table for the Heston model}%
\label{tablehes2}%
\end{table}

From Table \ref{tablehes2} and Proposition \ref{ui}, $(S_{t})_{0\leq
t\leq\infty}$ is a uniformly integrable $P$-martingale if and only if
$\alpha=\frac{2\kappa\theta}{\xi^{2}}<1,$ and $\gamma=\frac{2(\kappa-\rho\xi
)}{\xi^{2}}\geq0,$ which is equivalent to $\rho\xi\leq\kappa<\frac{\xi^{2}%
}{2\theta}.$ \. \qed

Under $P$, we have the following result on the positivity of the stock price
in the stopped Heston model.

\begin{pr}
\label{poshes} For the stopped Heston model \eqref{hes2b}, \newline(1)
$P(S_{T}>0)=1$ for all $T\in(0,\infty)$,\newline(2) $P(S_{\infty}>0)=1$ if and
only if $\kappa<\frac{\xi^{2}}{2\theta}$.
\end{pr}

\proof
Similar to the proofs of Proposition \ref{hm1} and Proposition \ref{uihes}
with $\gamma$ replaced by $\beta>0$ and $C_{2}$ by $C_{1}$, we obtain the
classification in Table \ref{tablehes3}. \newline\begin{table}[h]
\begin{center}%
\begin{tabular}
[c]{c|c|c|c|c|c|c}%
Case & $s(\ell)$ & $s(r)$ & $v(\ell)$ & $v(r)$ & $v_{b}(\ell)$ & $v_{b}%
(r)$\\\hline
$\alpha> 1$ & $-\infty$ & $\infty$ & $\infty$ & $\infty$ & $\infty$ & $\infty
$\\\hline
$\alpha= 1$ & $-\infty$ & $\infty$ & $\infty$ & $\infty$ & $\infty$ & $\infty
$\\\hline
$\alpha< 1$ & $>-\infty$ & $\infty$ & $<\infty$ & $\infty$ & $<\infty$ &
$\infty$\\\hline
\end{tabular}
\end{center}
\caption{Third classification table for the Heston model}%
\label{tablehes3}%
\end{table}

\noindent Based on Table \ref{tablehes3}, from Proposition \ref{mfs} and
Proposition \ref{t4fs}, we obtain the desired results. \qed

\subsection{$3/2$ stochastic volatility model}

Under $P$, the (correlated) 3/2 stochastic volatility model has the following
diffusive dynamics
\begin{align}
dS_{t}  &  =S_{t} \sqrt{Y_{t}}\mathds{1}_{t\in[0,\zeta)} dW_{t}^{(1)}, \quad
S_{0}=1,\nonumber\\
dY_{t}  &  =(\omega Y_{t}-\theta Y_{t}^{2})\mathds{1}_{t\in[0,\zeta)} dt +\xi
Y_{t}^{\frac{3}{2}}\mathds{1}_{t\in[0,\zeta)} dW_{t},\quad Y_{0}=x_{0}>0,
\label{threetc}%
\end{align}
where $\mathbb{E}^{P}[dW_{t}^{(1)}dW_{t}]=\rho dt$, $-1\leq\rho\leq1$,
$\omega> 0, \xi> 0, \theta\in\mathbb{R}$. The natural state space is given by
$J=(\ell,r)=(0,\infty)$. $\zeta$ is the possible exit time of the process $Y$
from its state space $J$. The model \eqref{threetc} belongs to the general
stochastic volatility model considered in \eqref{eq1tc} with $\mu(x)=\omega
x-\theta x^{2}, \sigma(x)=\xi x^{3/2}$, and $b(x)=\sqrt{x}$. Clearly
$\sigma(x)=\xi x^{3/2} \neq0, x\in J$, $\frac{1}{\sigma^{2}(x)}=\frac{1}%
{\xi^{2} x^{3}}\in L^{1}_{loc}(J)$, $\frac{\mu(x)}{\sigma^{2}(x)}=\frac
{\omega-\theta x}{\xi^{2} x^{2}}\in L^{1}_{loc}(J)$, and $\frac{b^{2}%
(x)}{\sigma^{2}(x)}=\frac{1}{\xi^{2} x^{2}}\in L^{1}_{loc}(J)$ are satisfied.
Thus, the conditions \eqref{cond1} and \eqref{cond2} are satisfied. From
Proposition \ref{newsde}, under $\widetilde{P}$, the diffusion $Y$ satisfies
the following SDE
\begin{align}
dY_{t}  &  =(\omega Y_{t} -\widetilde{\theta} Y_{t}^{2})\mathds{1}_{t\in
[0,\zeta)} dt +\xi Y_{t}^{\frac{3}{2}}\mathds{1}_{t\in[0,\zeta)}
d\widetilde{W}_{t}, \quad Y_{0}=x_{0}>0,\nonumber
\end{align}
where $\widetilde{\theta}=\theta-\rho\xi$. For a constant $c\in J$, the scale
functions of the SDE \eqref{y} and the SDE \eqref{yo} are respectively
\begin{align}
s(x)  &  =\frac{b}{c^{a}}\int_{c}^{x} y^{a}\exp\left(  \frac{d}{y}\right)  dy,
\quad\widetilde s(x)=\frac{b}{c^{\widetilde a}}\int_{c}^{x} y^{\widetilde
a}\exp\left(  \frac{d}{y}\right)  dy, \quad x\in\bar{J}, \label{sth}%
\end{align}
where $a=\frac{2\theta}{\xi^{2}}$, $b=\exp\left(  -\frac{2\omega}{c\xi^{2}%
}\right)  $, $d=\frac{2\omega}{\xi^{2}}$ and ${\widetilde a}=a-\frac{2\rho
}{\xi}$. Since the only difference between $s(\cdot)$ and $\widetilde{s}%
(\cdot)$ is in the parameters $a$ and $\widetilde a$, the analysis under
$\widetilde P$ is similar to the analysis under $P$, except with a change of
the parameter from $a$ to $\widetilde{a}$. Thus, we only need the results
under $P$. We have the following test functions
\begin{align}
v(x)  &  =\frac{2}{\xi^{2}}\int_{c}^{x}\frac{1}{y^{a+3}\exp\left(  \frac{d}%
{y}\right)  }\left(  \int_{y}^{x} z^{ a}\exp\left(  \frac{d}{z}\right)  dz
\right)  dy,\label{lop2s}\\
v_{b}(x)  &  =\frac{2}{\xi^{2}}\int_{c}^{x}\frac{1}{y^{a+2}\exp\left(
\frac{d}{y}\right)  }\left(  \int_{y}^{x} z^{ a}\exp\left(  \frac{d}%
{z}\right)  dz\right)  dy. \label{lop41}%
\end{align}

\begin{lem}
\label{p4} With $\omega>0$, the following properties are satisfied.
\[%
\begin{array}
[c]{cl|cl}%
(i) & a<-1 \Longleftrightarrow v(r)<\infty,\quad & (v) & \quad\widetilde a<-1
\Longleftrightarrow\widetilde v(r)<\infty.\\
(ii) & \forall a\in\mathbb{R}, \quad v_{b}(r)=\infty,\quad & (vi) &
\quad\forall\widetilde a\in\mathbb{R}, \quad\widetilde v_{b}(r)=\infty.\\
(iii) & \forall a\in\mathbb{R}, \quad v(\ell)=\infty,\quad & (vii) &
\quad\forall\widetilde a\in\mathbb{R}, \quad\widetilde v(\ell)=\infty.\\
(iv) & \forall a\in\mathbb{R}, \quad v_{b}(\ell)=\infty,\quad & (viii) &
\quad\forall\widetilde a\in\mathbb{R}, \quad\widetilde v_{b}(\ell)=\infty.
\end{array}
\]

\end{lem}

\proof
Details of the derivations can be found in  Appendix \ref{p4app}.
 \qed

\begin{pr}
\label{DR} For\footnote{Theorem $3$, page $110$ of Carr and Sun \citeyear{CS07}
proves sufficiency. See also Lewis \citeyear{L00}.} the $3/2$ model
\eqref{threetc}, the underlying stock price $(S_{t})_{0\leq t<\infty}$ is a
true $P$-martingale if and only if $\xi^{2} -2\rho\xi+2\theta\geq0$.
\end{pr}

\proof
From Lemma \ref{p4} and Proposition \ref{mm}, $(S_{t})_{0\leq t<\infty}$ is a
true $P$-martingale if and only if $\widetilde{a}\geq-1 $, which is equivalent
to $\xi^{2} -2\rho\xi+2\theta\geq0$ after some simplifications. \qed

\begin{pr}
\label{uithree} For the $3/2$ model \eqref{threetc}, the underlying stock
price $(S_{t})_{0\leq t\leq\infty}$ is not a uniformly integrable $P$-martingale.
\end{pr}

\proof
From Lemma \ref{p4}, for all $\widetilde a \in\mathbb{R}$, $\widetilde
v_{b}(r)=\infty$ and $\widetilde v_{b}(l)=\infty$ hold. From Proposition
\ref{ui}, $(S_{t})_{0\leq t\leq\infty}$ is not a uniformly integrable
$P$-martingale. \qed

Under $P$, we have the following result on the positivity of the stock price
in the $3/2$ model.

\begin{pr}
\label{posthree} For the $3/2$ model \eqref{threetc}, \newline(1)
$P(S_{T}>0)=1$ for all $T\in(0,\infty)$ if and only if $\xi^{2} +2\theta\geq
0$,\newline(2) $P(S_{\infty}>0)<1$.
\end{pr}

\proof
Similar to the proofs of Proposition \ref{DR} and Proposition \ref{uithree}
with $\widetilde a$ replaced by $a$, we obtain the classification in Table
\ref{tablethree2}. \begin{table}[h]
\begin{center}%
\begin{tabular}
[c]{c|c|c|c|c}%
Case & $v(\ell)$ & $v(r)$ & $v_{b}(\ell)$ & $v_{b}(r)$\\\hline
$a< -1$ & $\infty$ & $<\infty$ & $\infty$ & $\infty$\\\hline
$a \geq-1$ & $\infty$ & $\infty$ & $\infty$ & $\infty$\\\hline
\end{tabular}
\end{center}
\caption{Classification table for the $3/2$ model}%
\label{tablethree2}%
\end{table}Based on Table \ref{tablethree2}, Proposition \ref{mfs} and
Proposition \ref{t4fs}, we have the desired results. Note that $a\geq-1$ is
equivalent to $\xi^{2}+2\theta\geq0$. \qed

\subsection{Sch{\"o}bel-Zhu stochastic volatility model}

Under $P$, the correlated Sch{\"o}bel-Zhu stochastic volatility
model\footnote{It is the correlated version of the Stein-Stein \citeyear{SS91}
model. In Rheinl{\"a}nder \citeyear{R05}, the minimal entropy martingale
measure is studied in detail for this model, and its Proposition $3.1$ gives a
necessary and sufficient condition such that the associated stochastic
exponential is a true martingale. Here we provide deterministic criteria.}
(see Sch{\"o}bel and Zhu \citeyear{SZ99}) can be described by the following
diffusive dynamics
\begin{align}
dS_{t}  &  = S_{t} Y_{t} \mathds{1}_{t\in[0,\zeta)} dW_{t}^{(1)}, \quad
S_{0}=1,\nonumber\\
dY_{t}  &  =\kappa(\theta-Y_{t}) \mathds{1}_{t\in[0,\zeta)} dt +\gamma
\mathds{1}_{t\in[0,\zeta)} dW_{t}, \quad Y_{0}=x_{0}, \label{szmart}%
\end{align}
where $\mathbb{E}[dW^{(1)}_{t} dW_{t}]=\rho dt$, $-1\leq\rho\leq1$, $\kappa
>0$, $\theta>0$, $\gamma>0$. The process $Y$ is an Ornstein-Uhlenbeck process,
and this implies that its natural state space is $J=(\ell,r)=(-\infty,\infty)
$. $\zeta$ is the possible exit time of the process $Y$ from its state space
$J$. The model \eqref{szmart} belongs to the general stochastic volatility
model considered in \eqref{eq1tc} with $\mu(x)=\kappa(\theta-x)$,
$\sigma(x)=\gamma$, and $b(x)= x$. Clearly $\sigma(x)=\gamma\neq0, x\in J$,
then $\frac{1}{\sigma(x)^{2}}=\frac{1}{\gamma^{2}}\in L^{1}_{loc}(J)$,
$\frac{\mu(x)}{\sigma(x)^{2}}=\frac{\kappa(\theta-x)}{\gamma^{2} } \in
L^{1}_{loc}(J)$, and $\frac{b^{2}(x)}{\sigma^{2}(x)}=\frac{x^{2}}{\gamma^{2}%
}\in L^{1}_{loc}(J)$ are satisfied. Thus, the conditions \eqref{cond1} and
\eqref{cond2} are satisfied.

From Proposition \ref{newsde}, under $\widetilde{P}$, the diffusion $Y$
satisfies the following SDE
\begin{align}
dY_{t}  &  =(\kappa\theta-(\kappa-\rho\gamma)Y_{t})\mathds{1}_{t\in[0,\zeta)}
dt +\gamma\mathds{1}_{t\in[0,\zeta)} d\widetilde{W}_{t}, \quad Y_{0}%
=x_{0}.\nonumber
\end{align}

For a positive constant $c\in J$, denote $\alpha=\kappa-\rho\gamma$, and
compute the scale functions respectively of the SDE \eqref{y} and the SDE
\eqref{yo}
\begin{align}
s(x)  &  =\int_{c}^{x} e^{\frac{\kappa(y-\theta)^{2}-\kappa(c-\theta)^{2}
}{\gamma^{2}}} dy=C_{1} \int_{c}^{x} e^{\frac{\kappa(y-\theta)^{2}}{\gamma
^{2}}}dy,\nonumber\\
\widetilde{s}(x)  &  =\int_{c}^{x} e^{\frac{\alpha y^{2} -2\kappa\theta y
+2\kappa\theta c-\alpha c^{2} }{\gamma^{2}}}dy=
\begin{cases}
C_{2} \int_{c}^{x} e^{\frac{\alpha(y-\frac{\kappa\theta}{\alpha})^{2}}%
{\gamma^{2}}}dy, & \text{if $\alpha\neq0$},\notag\\
C_{3} \left(  e^{-\frac{2\kappa\theta c}{\gamma^{2}}}-e^{-\frac{2\kappa\theta
}{\gamma^{2}}x} \right)  , & \text{if $\alpha=0$},
\end{cases}
\end{align}
with constants $C_{1}=e^{-\kappa(c-\theta)^{2}/\gamma^{2}}>0$, $C_{2}%
=e^{(-\kappa^{2}\theta^{2}/\alpha+2\kappa\theta c -\alpha c^{2})/\gamma^{2}%
}>0$ for $\alpha\neq0$, and the constant $C_{3} =e^{2\kappa\theta c/\gamma
^{2}}\frac{\gamma^{2}}{2\kappa\theta}>0$ for $\alpha=0$. Since $\kappa>0$ by
assumption, $e^{\frac{\kappa(y-\theta)^{2}}{\gamma^{2}}}\geq1$ for any
$y\in[c,x]$, with $c\in J, x\in\bar{J}$, then $s(r)=s(\infty)=\infty$ always
holds, and consequently $v(r)=v(\infty)=\infty$.

\begin{pr}
\label{SZMMart} For the Sch{\"o}bel-Zhu model \eqref{szmart}, the underlying
stock price $(S_{t})_{0\leq t <\infty}$ is a true $P$-martingale.
\end{pr}

\proof
We now check the conditions in Proposition \ref{mm}. For the case of the right
endpoint $r$, depending on the sign of $\alpha=\kappa-\rho\gamma$, we obtain
the following classification
\begin{align}
\widetilde{s}(\infty)  &
\begin{cases}
<\infty, & \text{if $\alpha\leq0$},\notag\\
=\infty, & \text{if $\alpha>0$}.
\end{cases}
\end{align}

Details can be found in Appendix \ref{SZMMartapp}.
%
 \begin{table}[h]
\begin{center}%
\begin{tabular}
[c]{c|c|c|c|c|c|c}%
Case & $\widetilde s(\ell)$ & $\widetilde s(r)$ & $\widetilde v(\ell)$ &
$\widetilde v(r)$ & $\widetilde v_{b}(\ell)$ & $\widetilde v_{b}(r)$\\\hline
$\alpha\leq0$ & $>-\infty$ & $<\infty$ & $<\infty$ & $\infty$ & $<\infty$ &
$\infty$\\\hline
$\alpha>0$ & $>-\infty$ & $\infty$ & $<\infty$ & $\infty$ & $<\infty$ &
$\infty$\\\hline
\end{tabular}
\end{center}
\caption{First classification table for the Sch{\"o}bel-Zhu model}%
\label{tablesz1}%
\end{table}Above all, we can summarize the results in Table \ref{tablesz1}.
From Proposition \ref{mm} (3), $(S_{t})_{0\leq t<\infty}$ is a true
$P$-martingale. \qed

\begin{pr}
\label{uisz} For the Sch{\"o}bel-Zhu model \eqref{szmart}, the underlying
stock price $(S_{t})_{0\leq t\leq\infty}$ is a uniformly integrable
$P$-martingale if and only if $\kappa>\rho\gamma$.
\end{pr}

\proof
From Table \ref{tablesz1} and Proposition \ref{ui}, it follows that
$(S_{t})_{0\leq t\leq\infty}$ is a uniformly integrable $P$-martingale if and
only if $\alpha>0$, or equivalently $\kappa>\rho\gamma$. \qed


Under $P$, we obtain the following result on the positivity of the stock price
in the Sch{\"o}bel-Zhu model.

\begin{pr}
\label{posSZ} For the Sch{\"o}bel-Zhu model \eqref{szmart}, \newline(1)
$P(S_{T}>0)=1$ for all $T\in(0,\infty)$,\newline(2) $P(S_{\infty}>0)=1$.
\end{pr}

\proof
Similar to the proofs of Proposition \ref{SZMMart} and Proposition \ref{uisz}
with $\alpha$ replaced by $\kappa>0$, we obtain the classification given in
Table \ref{tablesz3}. \begin{table}[h]
\begin{center}%
\begin{tabular}
[c]{l|l|l|l|l|l|l}%
Case & $s(\ell)$ & $s(r)$ & $v(\ell)$ & $v(r)$ & $v_{b}(\ell)$ & $v_{b}%
(r)$\\\hline
$\alpha>0$ & $>-\infty$ & $\infty$ & $<\infty$ & $\infty$ & $<\infty$ &
$\infty$\\\hline
\end{tabular}
\end{center}
\caption{Second classification table for the Sch{\"o}bel-Zhu model}%
\label{tablesz3}%
\end{table}\qed

\subsection{Hull-White stochastic volatility model}

Under $P$, the correlated Hull-White stochastic volatility model (see Hull and
White \citeyear{HW87}) can be described by the following diffusive dynamics
\begin{align}
dS_{t}  &  = S_{t} \sqrt{Y_{t}} \mathds{1}_{t\in[0,\zeta)} dW_{t}^{(1)}, \quad
S_{0}=1,\nonumber\\
dY_{t}  &  =\mu Y_{t} \mathds{1}_{t\in[0,\zeta)} dt +\sigma Y_{t}
\mathds{1}_{t\in[0,\zeta)} dW_{t}, \quad Y_{0}=x_{0}>0, \label{hwmart}%
\end{align}
where $\mathbb{E}[dW^{(1)}_{t} dW_{t}]=\rho dt$, $-1\leq\rho\leq1$, $\mu>0$,
and $\sigma>0$. The process $Y$ is a geometric Brownian motion process, and
this implies that its natural state space is $J=(\ell,r)=(0,\infty)$. $\zeta$
is the possible exit time of the process $Y$ from its state space $J$. The
model \eqref{hwmart} belongs to the general stochastic volatility model
considered in \eqref{eq1tc} with $\mu(x)=\mu x$, $\sigma(x)=\sigma x $, and
$b(x)= \sqrt{x}$. Clearly $\sigma(x)=\sigma x \neq0, x\in J$, $\frac{1}%
{\sigma(x)^{2}}=\frac{1}{\sigma^{2} x^{2}}\in L^{1}_{loc}(J)$, $\frac{\mu
(x)}{\sigma(x)^{2}}=\frac{\mu}{\sigma^{2} x } \in L^{1}_{loc}(J)$, and
$\frac{b^{2}(x)}{\sigma^{2}(x)}=\frac{1}{\sigma^{2} x}\in L^{1}_{loc}(J)$ are
satisfied. Thus, the conditions \eqref{cond1} and \eqref{cond2} are satisfied.

From Proposition \ref{newsde}, under $\widetilde{P}$, the diffusion $Y$
satisfies the following SDE
\begin{align}
dY_{t}  &  =(\mu Y_{t}+\rho\sigma Y_{t}^{\frac{3}{2}})\mathds{1}_{t\in
[0,\zeta)} dt +\sigma Y_{t} \mathds{1}_{t\in[0,\zeta)} d\widetilde{W}_{t},
\quad Y_{0}=x_{0}>0, \label{fo}%
\end{align}

Denote $\alpha=\frac{4\mu}{\sigma^{2}}-1$ and $\gamma=\frac{4\rho}{\sigma}$.
For a constant $c\in J$, compute the scale functions of the SDE \eqref{yo}
\begin{align}
\widetilde s (x)  &  =\int_{c}^{x} e^{-\int_{c}^{y} \frac{2\mu u +2\rho\sigma
u^{3/2}}{\sigma^{2} u^{2}}du }dy =C_{1} \int_{c}^{x} y^{-\frac{\alpha+1}{2}}
e^{-\gamma\sqrt{y}}dy,\quad x\in\bar{J}, \label{schw}%
\end{align}
where $C_{1} =c^{\frac{2\mu}{\sigma^{2}}} e^{\frac{4\rho}{\sigma}\sqrt{c}}$ is
a positive constant. From the definition in \eqref{vb} and the scale function
in \eqref{schw}
\begin{align}
\widetilde v(x)  &  =\int_{c}^{x} \frac{2(\widetilde s(x)-\widetilde
s(y))}{\widetilde{s}^{\prime}(y) \widetilde\sigma^{2}(y)}dy=\frac{2}%
{\sigma^{2}} \int_{c}^{x} y^{\frac{\alpha-3}{2}} e^{\gamma\sqrt{y}} \left(
\int_{y}^{x} z^{-\frac{\alpha+1}{2}} e^{-\gamma\sqrt{z}}dz\right)  dy,
\label{vhw}%
\end{align}
and
\begin{align}
\widetilde v_{b}(x)  &  =\frac{2}{\sigma^{2}} \int_{c}^{x} y^{\frac{\alpha
-1}{2}} e^{\gamma\sqrt{y}} \left(  \int_{y}^{x} z^{-\frac{\alpha+1}{2}}
e^{-\gamma\sqrt{z}}dz\right)  dy. \label{vbhw}%
\end{align}

\begin{pr}
\label{hwmarthm} For\footnote{Proposition \ref{hwmarthm} is consistent with
Theorem $1$ of Jourdain \citeyear{J04}, and Proposition $2.5.$, page 34 of
Andersen and Piterbarg \citeyear{AP07}.} the Hull-White model \eqref{hwmart},
the underlying stock price $(S_{t})_{0\leq t<\infty}$ is a true $P$-martingale
if and only if $\rho\leq0$.
\end{pr}

\proof
We distinguish 3 situations: (I): $\mu>\frac{1}{2}\sigma^{2}$, (II):
$\mu=\frac{1}{2}\sigma^{2}$ and (III): $\mu<\frac{1}{2}\sigma^{2}$. Results
are summarized in Table \ref{tablehw1}. Details can be found in  Appendix \ref{hwapp}.
 The results in Table \ref{tablehw1}, combined with Proposition
\ref{mm} allow us to conclude if $(S_{t})_{0\leq t\leq T}, T\in(0,\infty)$ is
a true $P$-martingale. For $2\mu/\sigma^{2} >1$ ($\alpha>1$), $(S_{t})_{0\leq
t\leq T}, T\in(0,\infty)$ is a true $P$-martingale if and only if $\widetilde
v(r)=\infty$. This is equivalent to $\gamma\leq0$, and further equivalent to
$\rho\leq0$ from the definition of $\gamma$. When $2\mu/\sigma^{2}=1$
($\alpha=1$), $(S_{t})_{0\leq t\leq T}, T\in(0,\infty)$ is a true
$P$-martingale if and only if $\widetilde v(r)=\infty$, equivalently
$\gamma\leq0$, that is $\rho\leq0$. When $2\mu/\sigma^{2}<1$ ($\alpha<1$),
$(S_{t})_{0\leq t\leq T}, T\in(0,\infty)$ is a true $P$-martingale if and only
if $\widetilde v(r)=\infty$, equivalently $\gamma\leq0$, that is $\rho\leq0$.
\qed
\begin{table}[h]
\begin{center}%
\begin{tabular}
[c]{cc|c|c|c|c}%
Case &  & $\widetilde v(\ell)$ & $\widetilde v(r)$ & $\widetilde v_{b}(\ell)$
& $\widetilde v_{b}(r)$\\\hline
\multirow{2}{*}{ (I) $\mu>\frac{\sigma^2}{2}$} & $\gamma\leq0$ & $\infty$ &
$\infty$ & $\infty$ & $\infty$\\\cline{2-6}
& $\gamma> 0$ & $\infty$ & $<\infty$ & $\infty$ & $\infty$\\\hline\hline
\multirow{2}{*}{ (II) $\mu=\frac{\sigma^2}{2}$} & $\gamma\leq0$ & $\infty$ &
$\infty$ & $\infty$ & $\infty$\\\cline{2-6}
& $\gamma> 0$ & $\infty$ & $<\infty$ & $\infty$ & $\infty$\\\hline\hline
\multirow{2}{*}{ (III) $\mu<\frac{\sigma^2}{2}$} & $\gamma\leq0$ & $\infty$ &
$\infty$ & $<\infty$ & $\infty$\\\cline{2-6}
& $\gamma> 0$ & $\infty$ & $<\infty$ & $<\infty$ & $\infty$\\\hline
\end{tabular}
\end{center}
\caption{First classification table for the Hull-White model}%
\label{tablehw1}%
\end{table}

\begin{pr}
\label{uihw} For the Hull-White model \eqref{hwmart}, the underlying stock
price $(S_{t})_{0\leq t\leq\infty}$ is a uniformly integrable $P$-martingale
if and only if $\mu<\frac{1}{2}\sigma^{2}$ and $\rho\leq0$.
\end{pr}

\proof
The proof of Proposition \ref{uihw} requires the same 3 cases as Proposition
\ref{hwmarthm}. Results are summarized in Table \ref{tablehw4}. Details can be found in Appendix \ref{uihwapp}. \qed
\begin{table}[h]
\begin{center}%
\begin{tabular}
[c]{cc|c|c|c|c|c|c}%
Case &  & $\widetilde s(\ell)$ & $\widetilde s(r)$ & $\widetilde v(\ell)$ &
$\widetilde v(r)$ & $\widetilde v_{b}(\ell)$ & $\widetilde v_{b}(r)$\\\hline
\multirow{3}{*}{ (I) $\mu>\frac{\sigma^2}{2}$} & $\gamma>0$ & $-\infty$ &
$<\infty$ & $\infty$ & $<\infty$ & $\infty$ & $\infty$\\\cline{2-8}
& $\gamma=0$ & $-\infty$ & $<\infty$ & $\infty$ & $\infty$ & $\infty$ &
$\infty$\\\cline{2-8}
& $\gamma<0$ & $-\infty$ & $\infty$ & $\infty$ & $\infty$ & $\infty$ &
$\infty$\\\hline\hline
\multirow{3}{*}{ (II) $\mu=\frac{\sigma^2}{2}$} & $\gamma>0$ & $-\infty$ &
$<\infty$ & $\infty$ & $<\infty$ & $\infty$ & $\infty$\\\cline{2-8}
& $\gamma=0$ & $-\infty$ & $\infty$ & $\infty$ & $\infty$ & $\infty$ &
$\infty$\\\cline{2-8}
& $\gamma<0$ & $-\infty$ & $\infty$ & $\infty$ & $\infty$ & $\infty$ &
$\infty$\\\hline\hline
\multirow{3}{*}{ (III) $\mu<\frac{\sigma^2}{2}$} & $\gamma>0$ & $>-\infty$ &
$<\infty$ & $\infty$ & $<\infty$ & $<\infty$ & $\infty$\\\cline{2-8}
& $\gamma=0$ & $>-\infty$ & $\infty$ & $\infty$ & $\infty$ & $<\infty$ &
$\infty$\\\cline{2-8}
& $\gamma<0$ & $>-\infty$ & $\infty$ & $\infty$ & $\infty$ & $<\infty$ &
$\infty$\\\hline
\end{tabular}
\end{center}
\caption{Second classification table for the Hull-White model}%
\label{tablehw4}%
\end{table}


Under $P$, we have the following result on the positivity of the stock price
in the Hull-White model.

\begin{pr}
\label{poshw} For the Hull-White model \eqref{hwmart}, \newline(1)
$P(S_{T}>0)=1$ for all $T\in(0,\infty)$,\newline(2) $P(S_{\infty}>0)=1$ if and
only if $\frac{2\mu}{\sigma^{2}}<1$.
\end{pr}

\proof
Similar to the proofs of Proposition \ref{hwmarthm} and Proposition \ref{uihw}
with $\gamma=0$, we have the classification in Table \ref{tablehw7}%
.\newline\begin{table}[h]
\begin{center}%
\begin{tabular}
[c]{c|c|c|c|c|c|c}%
Case & $s(\ell)$ & $s(r)$ & $v(\ell)$ & $v(r)$ & $v_{b}(\ell)$ & $v_{b}%
(r)$\\\hline
(I) $\mu>\frac{\sigma^{2}}{2}$ & $-\infty$ & $<\infty$ & $\infty$ & $\infty$ &
$\infty$ & $\infty$\\\hline
(II) $\mu=\frac{\sigma^{2}}{2}$ & $-\infty$ & $\infty$ & $\infty$ & $\infty$ &
$\infty$ & $\infty$\\\hline
(III) $\mu<\frac{\sigma^{2}}{2}$ & $>-\infty$ & $\infty$ & $\infty$ & $\infty$
& $<\infty$ & $\infty$\\\hline
\end{tabular}
\end{center}
\caption{Third classification table for the Hull-White model}%
\label{tablehw7}%
\end{table}

From Table \ref{tablehw7}, Proposition \ref{mfs} and Proposition \ref{t4fs},
we obtain the desired results. \qed

\subsection{Summary of the Examples}

Table \ref{SUMM} and Table \ref{SUMM2} summarize the results obtained
throughout Section \ref{ms6}. In all cases, we study the ``stopped'' price
process as we assume that there are two absorbing barriers at $\ell$ and $r$.
Conditions for uniformly integrable martingales are stronger than those for a
true martingale on $(0,\infty)$. Similar remarks hold for the positivity of
$S_{T}$ and $S_{\infty}$, where $0<T<\infty$. \begin{table}[tbh]
\begin{center}%
\begin{tabular}
[c]{c|c|c}%
{Model} & True martingale on $(0,\infty)$ & UI martingale on $[0,\infty
$]\\\hline
\multirow{2}{*}{Heston Model \eqref{hes2b}} & under Feller condition & never
under Feller condition\\
& $\kappa> \frac{\xi^{2}}{2\theta}$ & $\rho\xi\leq\kappa<\frac{\xi^{2}%
}{2\theta}$\\\hline
3/2 Model \eqref{threetc} & $\xi^{2}+2\theta\ge\max(0,2\rho\xi)$ &
Never\\\hline
Sch\"obel-Zhu Model \eqref{szmart} & Always & when $\kappa>\rho\gamma$\\\hline
Hull White Model \eqref{hwmart} & $\rho\le0$ & $\mu<\frac{\sigma^{2}}{2}$ and
$\rho\le0$\\\hline
\end{tabular}
\end{center}
\caption{Summary of conditions for (uniformly integrable) martingales.}%
\label{SUMM}%
\end{table}

\begin{table}[tbh]
\begin{center}%
\begin{tabular}
[c]{c|c|c}%
{Model} & $S_{T}$ positive $P$-a.s. & $S_{\infty}$ positive $P$-a.s.\\\hline
\multirow{2}{*}{Heston Model \eqref{hes2b}} & Always & never under Feller
condition\\
&  & $\kappa<\frac{\xi^{2}}{2\theta}$\\\hline
3/2 Model \eqref{threetc} & $\xi^{2}+2\theta\ge0$ & Never\\\hline
Sch\"obel-Zhu Model \eqref{szmart} & Always & Always\\\hline
Hull White Model \eqref{hwmart} & Always & $\mu<\frac{\sigma^{2}}{2}$\\\hline
\end{tabular}
\end{center}
\caption{Summary of conditions for positivity of stock prices.}%
\label{SUMM2}%
\end{table}

\section{Concluding Remarks \label{ms10}}

This paper generalizes some results of Mijatovi\'c and Urusov
(\citeyear*{MU12FS}, \citeyear*{MU12PTRF}) concerning the (uniformly
integrable) martingale property of the asset price from the case $\rho=1$ to the case $-1\leq\rho\leq1$, and provides new direct proofs without
using the concept of ``separating times". We also obtain deterministic
criteria for the convergence or divergence of both perpetual and capped
integral functionals of time-homogeneous diffusions. Explicit deterministic
criteria for checking the (uniformly integrable) martingale properties for
four stochastic volatility models are provided. Future research directions
include finding necessary and sufficient deterministic conditions for the
martingale property of time-changed L\'evy processes with non-zero correlation
(Carr and Wu \citeyear{CW04}), of which the time-homogeneous stochastic
volatility models considered in this paper are special cases.

\bibliographystyle{plain}
\bibliography{muc}

\newpage\appendix

\section{Technical Conditions on the Probability space and Filtration}

Throughout the paper we assume a space accommodating all four processes
$(Y,Z,W,W^{(1)})$ in \eqref{eq1tc}. This is described below, following closely
the presentation in Appendix B of Carr, Fisher and Ruf \citeyear{CFR13}. For a
fixed time horizon $T\in(0,\infty]$, we require a stochastic basis
$(\Omega,\mathcal{F}_{T},\{\mathcal{F}_{t}\}_{t\in\lbrack0,T]},P)$ with a
right-continuous filtration $\{\mathcal{F}_{t}\}_{t\in\lbrack0,T]}$. As in
Carr, Fisher and Ruf \citeyear{CFR13} and Föllmer \citeyear{F72}, page 156,
\ for any stopping time $\tau$, we define $\mathcal{F}_{\tau}:=\{A\in
\mathcal{F}_{T}\mid A\cap\{\tau\leq t\}\in\mathcal{F}_{t}\text{ for all }%
t\in\lbrack0,T]\}$ and $\mathcal{F}_{\tau-}:=\sigma(\{A\cap\{\tau
>t\}\in\mathcal{F}_{T}\mid A\in\mathcal{F}_{t}\text{\ for some }t\in
\lbrack0,T]\cup\mathcal{F}_{0}\})$. In general, non-negative random variables
are permitted to take values in the set $[0,\infty]$ and stopping times $\tau$
are permitted to take values in the set $[0,\infty]\cup\mathscr{T}$ for some
transfinite time $\mathscr{T}>T$. In special cases we may restrict the range
of stopping times.

Let $\Omega_{1}$ denote the space of continuous paths $\omega_{1}%
:[0,\infty)\rightarrow\bar{J}$ with $\omega_{1}(0)\in J$. \ Define
$\zeta(\omega_{1}):=\inf\{t\in\lbrack0,T]:\omega_{1}(t)\not \in J\}$ with the
convention $\inf\emptyset=\mathscr{T}$ . Assume that $\omega_{1}$ stays at
either $\ell$ or $r$ once it hits it, i.e. that $\omega_{1}(\zeta
+s)=\omega_{1}(\zeta)$ for all $s>0$ on the set $\{\zeta<\infty\}$. Let
$\Omega_{2}$ denote the space of continuous\footnote{Continuity where the
function takes the value $\infty$ is defined as usual through a
compactification: if $\lim_{t\rightarrow t_{0}}\omega_{2}(t)=\infty$, then
$\omega_{2}$ is continuous at $t_{0}.$} paths $\omega_{2}:[0,\infty
)\rightarrow\lbrack0,\infty]$ with $\omega_{2}(0)=1$. As in Carr, Fisher and
Ruf \citeyear{CFR13}, define for all $i\in\mathbb{N}$, $R_{i}:=\inf
\{t\in\lbrack0,T]:\omega_{2}(t)>i\}$, and $S_{i}:=\inf\{t\in\lbrack
0,T]:\omega_{2}(t)<\frac{1}{i}\}$. Then $R:=\lim\limits_{i\uparrow\infty}%
R_{i}$, $S:=\lim\limits_{i\uparrow\infty}S_{i}$ are respectively the first
hitting time of infinity and zero by $\omega_{2}$, with the convention
$\inf\emptyset=\mathscr{T}$. Assume that $\omega_{2}(R+s)=\omega_{2}(R)$ for
all $s>0$ on $\{R<\infty\}$ and similarly $\omega_{2}(S+s)=\omega_{2}(S)$ for
all $s>0$ on $\{S<\infty\}$, so that $\omega_{2}$ stays at zero or infinity
once it hits it. Let $\Omega_{3}$ denote the space of continuous paths
$\omega_{3}:[0,\infty)\rightarrow\mathbb{R}$ with $\omega_{3}(0)=0$. Similarly
let $\Omega_{4}$ denote the space of continuous paths $\omega_{4}%
:[0,\infty)\rightarrow\mathbb{R}$ with $\omega_{4}(0)=0$. \

Denote $\Omega=\prod_{i=1}^{4}\Omega_{i}$ and $\omega=(\omega_{1},\omega
_{2},\omega_{3},\omega_{4})$. As in Appendix B of Carr, Fisher and Ruf
\citeyear{CFR13},  we require that $\{\mathcal{F}_{R_{i}-}\}_{i\in\mathbb{N}}$
is a standard system, see Remark $6.1.1$ of Föllmer \citeyear{F72}, so that in
the proof of Proposition \ref{r2}, the Extension Theorem V$4.1$ of
Parthasarathy \citeyear{P67} can be applied, and any probability measure on
$\mathcal{F}_{R-}$ has a (possibly non-unique) extension to a probability
measure on $\mathcal{F}_{T}.$ Such a canonical filtration can be constructed
as in Appendix B of Carr, Fisher and Ruf \citeyear{CFR13}. 
 
Given the canonical space $(\Omega,\mathcal{F}_{T},(\mathcal{F}_{t}%
)_{t\in\lbrack0,T]})$, the processes $(Y,Z,W,W^{(1)})$ in \eqref{eq1tc}
correspond respectively to the four components of $\omega$ and are formally
functions of $\omega$. We assume that processes $Y,Z$ are adapted to the
filtration $\{\mathcal{F}_{t}\}_{t\in\lbrack0,T]}$, as are $W,W^{(1)}$,  which
are assumed to be Brownian motions with respect to the same filtration.

\section{Online Supplement: Proofs of the Examples in Section \ref{ms6}}

\subsection{Proof of Proposition \ref{hm1} in the Heston model\label{hm1app}}

\proof
We first establish the following lemma.

\begin{lem}
\label{lem2}
\begin{gather}
\int_{c}^{\infty}y^{-\alpha}e^{\gamma y}dy%
\begin{cases}
<\infty, & \text{if }\gamma\text{$<0$ or }\gamma=0,\alpha>1\\
=\infty, & \text{if }\gamma>0\text{ or }\gamma=0,\alpha\leq1
\end{cases}
\text{ for }c>0\label{L2}\\
\int_{x}^{\infty}y^{-\alpha}e^{\gamma y}dy\sim\frac{1}{-\gamma}x^{-\alpha
}e^{\gamma x}, \quad\text{if }\gamma<0,\text{ as }x\rightarrow\infty,
\quad\label{L3}\\
\int_{0}^{x}y^{-\alpha}e^{\gamma y}dy%
\begin{cases}
<\infty, & \text{if $\alpha<1$},\\
=\infty, & \text{if $\alpha\geq1$}.
\end{cases}
\label{L4}\\
\int_{0}^{x}y^{-\alpha}e^{\gamma y}dy\sim%
\begin{cases}
\frac{x^{1-\alpha}}{\alpha-1}, & \text{if $\alpha>1$},\\
\ln(x), & \text{if }\alpha=1.
\end{cases}
\text{ as }x\rightarrow0. \label{L5}%
\end{gather}

\end{lem}

\proof Let us verify (\ref{L2}) and (\ref{L3}) for example. From L'Hôpital's
rule, since the numerator and denominator both approach $0$ if $\gamma<0$, \
\[
\lim_{x\rightarrow\infty}\frac{\int_{x}^{\infty}y^{-\alpha}e^{\gamma y}%
dy}{x^{-\alpha}e^{\gamma x}} =\lim_{x\rightarrow\infty}\frac{-x^{-\alpha
}e^{\gamma x}}{(\gamma x^{-\alpha}-\alpha x^{-\alpha-1})e^{\gamma x}}
=\lim_{x\rightarrow\infty}\frac{-1}{\gamma-\alpha/x} =-\frac{1}{\gamma}.
\]
The other asymptotics are similarly obtained. \qed

Then we will check the conditions of Proposition \ref{mm}. Here
\[
\widetilde{v}(x) \equiv\frac{2}{\xi^{2}}\int_{c}^{x}y^{\alpha-1}e^{-\gamma
y}\left(  \int_{y}^{x}z^{-\alpha}e^{\gamma z}dz\right)  dy,\quad\widetilde
{v}_{b}(x) \equiv\frac{2}{\xi^{2}}\int_{c}^{x}y^{\alpha}e^{-\gamma y}\left(
\int_{y}^{x}z^{-\alpha}e^{\gamma z}dz\right)  dy.
\]

It follows from the Lemma \ref{lem2} that $\widetilde{s}(0)$ is finite if and
only if $\alpha<1$ and $\widetilde{s}(\infty)$ is finite if and only if
$\gamma<0$ or $\gamma=0$ and $\alpha>1.$ We consider several cases.

\begin{itemize}
\item $\alpha=\frac{2\kappa\theta}{\xi^{2}}>1,$ $\gamma>0.$ In this case,
$\widetilde{s}(0)=-\infty,\widetilde{v}(0)=\infty$ and $\widetilde{v}%
_{b}(0)=\infty,$ $\widetilde{s}(\infty)=\infty$ so that $\widetilde{v}%
(\infty)=\infty.$ Therefore using Proposition \ref{mm}, (1) $E(S_{T})=1$ since
$\widetilde{v}(\infty)=\infty$ holds.

\item $\alpha>1,$ $\gamma=0.$ In this case $\widetilde{s}(0)=-\infty
,\widetilde{v}(0)=\infty$ and $\widetilde{v}_{b}(0)=\infty.$ Moreover%
\begin{equation}
\widetilde{v}(\infty)\equiv\frac{2}{\xi^{2}}\int_{c}^{\infty}y^{\alpha
-1}\left(  \int_{y}^{\infty}z^{-\alpha}dz\right)  dy=C\int_{c}^{\infty
}y^{\alpha-1}y^{1-\alpha}dy=\infty
\end{equation}
Again from Proposition \ref{mm}, (i) $E(S_{T})=1$ if and only if either
$\widetilde{v}(\infty)=\infty$ holds or $\widetilde{v}_{b}(\infty)<\infty$.

\item $\alpha>1,$ $\gamma<0.$ In this case $\widetilde{s}(0)=-\infty
,\widetilde{v}(0)=\infty$ and $\widetilde{v}_{b}(0)=\infty.$ Again we will
show that $\widetilde{v}(\infty)=\infty$. By the Lemma \ref{lem2}, for some
positive constant $C,$
\[
\widetilde{v}(\infty)\equiv\frac{2}{\xi^{2}}\int_{c}^{\infty}y^{\alpha
-1}e^{-\gamma y}\left(  \int_{y}^{\infty}z^{-\alpha}e^{\gamma z}dz\right)
dy\geq C\int_{c}^{\infty}y^{\alpha-1}e^{-\gamma y}y^{-\alpha}e^{\gamma
y}dy=\infty
\]

\item $\alpha=\frac{2\kappa\theta}{\xi^{2}}=1,$ $\gamma\geq0.$ In this case
$\widetilde{s}(0)=-\infty$ since the integral $\int_{0}^{c}y^{-\alpha
}e^{\gamma y}dy$ is divergent at 0, so that $\widetilde{v}(0)=\infty.$ Also
since $\gamma\geq0,$ $\widetilde{s}(\infty)=C_{2}\int_{c}^{\infty}y^{-\alpha
}e^{\gamma y}dy\geq C_{2}\int_{c}^{\infty}y^{-\alpha}dy=\infty$ so that
$\widetilde{v}(\infty)=\infty$ and Proposition \ref{mm} (1) applies.

\item $\alpha=\frac{2\kappa\theta}{\xi^{2}}=1,$ $\gamma<0.$ Here
$\widetilde{s}(0)=-\infty$ since the integral $\int_{0}^{c}y^{-1}e^{\gamma
y}dy$ diverges around 0 so that $\widetilde{v}(0)=\infty.$ However
$\widetilde{s}(\infty)=C\int_{c}^{\infty}y^{-1}e^{\gamma y}dy<\infty.$ Then
\[
\widetilde{v}(\infty)\equiv\frac{2}{\xi^{2}}\int_{c}^{\infty}e^{-\gamma
y}\left(  \int_{y}^{\infty}z^{-1}e^{\gamma z}dz\right)  dy\geq C\int
_{c}^{\infty}e^{-\gamma y}y^{-1}e^{\gamma y}dy=\infty
\]
so again in this case $\widetilde{v}(0)=\widetilde{v}(\infty)=\infty$ and
Proposition \ref{mm} (1) applies.

\item $\alpha=\frac{2\kappa\theta}{\xi^{2}}<1,$ $\gamma\geq0.$ In this case
$\widetilde{s}(0)$ is finite since the integral $\int_{0}^{c}y^{-\alpha
}e^{\gamma y}dy$ is convergent at 0. Also since $\gamma\geq0,$ $\widetilde
{s}(\infty)=C_{2}\int_{c}^{\infty}y^{-\alpha}e^{\gamma y}dy\geq C_{2}\int
_{c}^{\infty}y^{-\alpha}dy=\infty$ so that $\widetilde{v}(\infty)=\infty.$ In
this case, from the Lemma \ref{lem2}, $\int_{0}^{y}z^{-\alpha}e^{\gamma
z}dz\sim\frac{y^{1-\alpha}}{1-\alpha}$ and
\begin{align*}
\widetilde{v}(0)  &  \equiv\frac{2}{\xi^{2}}\int_{0}^{c}y^{\alpha-1}e^{-\gamma
y}\left(  \int_{0}^{y}z^{-\alpha}e^{\gamma z}dz\right)  dy\leq C\int_{0}%
^{c}y^{\alpha-1}e^{-\gamma y}\frac{y^{1-\alpha}}{1-\alpha}dy<\infty
\end{align*}
Similarly,
\begin{align*}
\widetilde{v}_{b}(0)  &  \equiv\frac{2}{\xi^{2}}\int_{0}^{c}y^{\alpha
}e^{-\gamma y}\left(  \int_{0}^{y}z^{-\alpha}e^{\gamma z}dz\right)  dy \leq
C\int_{0}^{c}y^{\alpha}e^{-\gamma y}\frac{y^{1-\alpha}}{1-\alpha}dy<\infty
\end{align*}
So in this case $\widetilde{v}(\infty)=\infty$ and $\widetilde{v}%
_{b}(0)<\infty$ so that Proposition \ref{mm} (3) applies.

\item $\alpha=\frac{2\kappa\theta}{\xi^{2}}<1,$ $\gamma<0.$ In this case
$\widetilde{s}(0)$ is finite since the integral $\int_{0}^{c}y^{-\alpha
}e^{\gamma y}dy$ is convergent at 0. Also
\begin{align*}
\widetilde{v}(0)  &  \equiv\frac{2}{\xi^{2}}\int_{0}^{c}y^{\alpha-1}e^{-\gamma
y}\left(  \int_{0}^{y}z^{-\alpha}e^{\gamma z}dz\right)  dy \leq C\int_{0}%
^{c}y^{\alpha-1}e^{-\gamma y}\frac{y^{1-\alpha}}{1-\alpha}dy<\infty\\
\widetilde{v}_{b}(0)  &  \equiv\frac{2}{\xi^{2}}\int_{0}^{c}y^{\alpha
}e^{-\gamma y}\left(  \int_{0}^{y}z^{-\alpha}e^{\gamma z}dz\right)  dy \leq
C\int_{0}^{c}y^{\alpha}e^{-\gamma y}\frac{y^{1-\alpha}}{1-\alpha}dy<\infty
\end{align*}
Also, since $\gamma<0,$ $\widetilde{s}(\infty)=C_{2}\int_{c}^{\infty
}y^{-\alpha}e^{\gamma y}dy<\infty.$ Then
\[
\widetilde{v}(\infty)\equiv\frac{2}{\xi^{2}}\int_{c}^{\infty}y^{\alpha
-1}e^{-\gamma y}\left(  \int_{y}^{\infty}z^{-\alpha}e^{\gamma z}dz\right)
dy\geq C\int_{c}^{\infty}y^{\alpha-1}e^{-\gamma y}y^{-\alpha}e^{\gamma
y}dy=\infty
\]
so that $\widetilde{v}(\infty)=\infty.$ In this case $\widetilde{v}%
_{b}(0)<\infty$ and $\widetilde{v}(\infty)=\infty$ so that Proposition
\ref{mm} (3) applies.
\end{itemize}

In summary, for the Heston model, $\{S_{t}\}_{t\leq T}$ is always a
martingale. \qed

\subsection{Proof of Proposition \ref{uihes} in the Heston model
\label{hm2app}}

\proof
We will check the conditions of Proposition \ref{ui}.

It follows from the Lemma \ref{lem2} that $\widetilde{s}(0)$ is finite if and
only if $\alpha<1$ and $\widetilde{s}(\infty)$ is finite if and only if
$\gamma<0$ or $\gamma=0$ and $\alpha>1.$ Note that
\[
\widetilde{v}_{b}(\infty)\equiv\frac{2}{\xi^{2}}\int_{c}^{\infty}y^{\alpha
}e^{-\gamma y}\left(  \int_{y}^{\infty}z^{-\alpha}e^{\gamma z}dz\right)  dy.
\]
Also,%
\[
\widetilde{v}_{b}(0)\equiv\frac{2}{\xi^{2}}\int_{0}^{c}y^{\alpha}e^{-\gamma
y}\left(  \int_{0}^{y}z^{-\alpha}e^{\gamma z}dz\right)  dy<\infty
\]
if and only if $\alpha<1$.

Again consider several cases.

\begin{itemize}
\item $\alpha>1,$ $\gamma>0.$ In this case $\widetilde{s}(0)=-\infty$ and
$\widetilde{v}_{b}(0)=\infty,$ $\widetilde{s}(\infty)=\infty$ so that
$\widetilde{v}_{b}(\infty)=\infty.$ Therefore, since $\widetilde{v}%
_{b}(0)=\widetilde{v}_{b}(\infty)=\infty,$ none of the conditions of
Proposition \ref{ui} apply.

\item $\alpha>1,$ $\gamma=0.$ In this case, $\widetilde{s}(0)=-\infty
,\widetilde{v}_{b}(0)=\infty.$ Moreover,
\begin{equation}
\widetilde{v}_{b}(\infty)\equiv\frac{2}{\xi^{2}}\int_{c}^{\infty}y^{\alpha
}\left(  \int_{y}^{\infty}z^{-\alpha}dz\right)  dy=C\int_{c}^{\infty}%
y^{\alpha}y^{1-\alpha}dy=\infty.
\end{equation}

Again, since $\widetilde{v}_{b}(0)=\widetilde{v}_{b}(\infty)=\infty,$ none of
the conditions of Proposition \ref{ui} apply.

\item $\alpha>1,$ $\gamma<0.$ In this case $\widetilde{s}(0)=-\infty,$ and
$\widetilde{v}_{b}(0)=\infty.$ Again we will show that $\widetilde{v}%
_{b}(\infty)=\infty$. By the Lemma \ref{lem2}, for some positive constant
$C,$
\[
\widetilde{v}_{b}(\infty)\equiv\frac{2}{\xi^{2}}\int_{c}^{\infty}y^{\alpha
}e^{-\gamma y}\left(  \int_{y}^{\infty}z^{-\alpha}e^{\gamma z}dz\right)
dy\geq C\int_{c}^{\infty}y^{\alpha}e^{-\gamma y}y^{-\alpha}e^{\gamma
y}dy=C\int_{c}^{\infty}1dy=\infty
\]
Again, since $\widetilde{v}_{b}(0)=\widetilde{v}_{b}(\infty)=\infty,$ none of
the conditions of Proposition \ref{ui} apply.

\item $\alpha=1,$ $\gamma\geq0.$ In this case $\widetilde{s}(0)=-\infty$ since
the integral $\int_{0}^{c}y^{-\alpha}e^{\gamma y}dy$ is divergent at 0, so
that $\widetilde{v}_{b}(0)=\infty.$ Also since $\gamma\geq0,$ $\widetilde
{s}(\infty)=C_{2}\int_{c}^{\infty}y^{-\alpha}e^{\gamma y}dy\geq C_{2}\int
_{c}^{\infty}y^{-\alpha}dy=\infty$ so that $\widetilde{v}_{b}(\infty)=\infty$.
Since $\widetilde{v}_{b}(0)=\widetilde{v}_{b}(\infty)=\infty,$ none of the
conditions of Proposition \ref{ui} apply.

\item $\alpha=1,$ $\gamma<0.$ Here $\widetilde{s}(0)=-\infty$ since the
integral $\int_{0}^{c}y^{-1}e^{\gamma y}dy$ diverges around 0 so that
$\widetilde{v}_{b}(0)=\infty.$ However $\widetilde{s}(\infty)=C\int
_{c}^{\infty}y^{-1}e^{\gamma y}dy<\infty.$ Then
\[
\widetilde{v}_{b}(\infty)\equiv\frac{2}{\xi^{2}}\int_{c}^{\infty}%
y^{1}e^{-\gamma y}\left(  \int_{y}^{\infty}z^{-1}e^{\gamma z}dz\right)  dy\geq
C\int_{c}^{\infty}ye^{-\gamma y}y^{-1}e^{\gamma y}dy=\infty.
\]
Since $\widetilde{v}_{b}(0)=\widetilde{v}_{b}(\infty)=\infty,$ none of the
conditions of Proposition \ref{ui} apply.

\item $\alpha<1,$ $\gamma\geq0.$ In this case $\widetilde{s}(0)$ is finite
since the integral $\int_{0}^{c}y^{-\alpha}e^{\gamma y}dy$ is convergent at 0.
Also since $\gamma\geq0,$ $\widetilde{s}(\infty)=C_{2}\int_{c}^{\infty
}y^{-\alpha}e^{\gamma y}dy\geq C_{2}\int_{c}^{\infty}y^{-\alpha}dy=\infty$ so
that $\widetilde{v}_{b}(\infty)=\infty.$ In this case, from the Lemma
\ref{lem2}, $\int_{0}^{y}z^{-\alpha}e^{\gamma z}dz\sim\frac{y^{1-\alpha}%
}{1-\alpha}$ and
\begin{align*}
\widetilde{v}_{b}(0)  &  \equiv\frac{2}{\xi^{2}}\int_{0}^{c}y^{\alpha
}e^{-\gamma y}\left(  \int_{0}^{y}z^{-\alpha}e^{\gamma z}dz\right)  dy \leq
C\int_{0}^{c}y^{\alpha}e^{-\gamma y}\frac{y^{1-\alpha}}{1-\alpha}dy<\infty
\end{align*}
In this case $\widetilde{v}_{b}(0)<\infty$ and $\widetilde{s}(\infty)=\infty$
so that the condition $C^{\prime}$ of Proposition \ref{ui} is satisfied.

\item $\alpha<1,$ $\gamma<0.$ In this case $\widetilde{s}(0)$ is finite since
the integral $\int_{0}^{c}y^{-\alpha}e^{\gamma y}dy$ is convergent at 0. Also
\begin{align*}
\widetilde{v}_{b}(0)  &  \equiv\frac{2}{\xi^{2}}\int_{0}^{c}y^{\alpha
}e^{-\gamma y}\left(  \int_{0}^{y}z^{-\alpha}e^{\gamma z}dz\right)  dy \leq
C\int_{0}^{c}y^{\alpha}e^{-\gamma y}\frac{y^{1-\alpha}}{1-\alpha}dy<\infty
\end{align*}
and $\widetilde{s}(\infty)=C_{2}\int_{c}^{\infty}y^{-\alpha}e^{\gamma
y}dy<\infty.$ Therefore, condition $D^{\prime}$ of Proposition \ref{ui} holds
if and only if $\widetilde{v}_{b}(\infty)<\infty.$ But
\[
\widetilde{v}_{b}(\infty)\equiv\frac{2}{\xi^{2}}\int_{c}^{\infty}y^{\alpha
}e^{-\gamma y}\left(  \int_{y}^{\infty}z^{-\alpha}e^{\gamma z}dz\right)
dy\geq C\int_{c}^{\infty}y^{\alpha}e^{-\gamma y}y^{-\alpha}e^{\gamma
y}dy=C\int_{c}^{\infty}1dy=\infty,
\]
so the conditions of Proposition \ref{ui} fail.
\end{itemize}

In summary, for the Heston model, $\{S_{t};t\leq\infty\}$ is a uniformly
integrable martingale if and only if $\alpha=\frac{2\kappa\theta}{\xi^{2}}<1,$
and $\gamma=2\frac{\kappa-\rho\xi}{\xi^{2}}\geq0,$ i.e. if and only if
$\rho\xi\leq\kappa<\frac{\xi^{2}}{2\theta}.$ \qed

\subsection{Proof of Lemma \ref{p4}\label{p4app}}

\proof For the right boundary $r$, divide into two cases:

(i) When $a<-1$, $\lim\limits_{y\rightarrow\infty} y^{a+1}\exp\left(  \frac
{d}{y}\right)  =0$. From L'Hôpital's rule
\begin{align}
\lim_{y\rightarrow\infty} \frac{\int_{y}^{\infty} z^{a}\exp\left(  \frac{d}%
{z}\right)  dz}{y^{a+1}\exp\left(  \frac{d}{y}\right)  }  &  =-\frac{1}%
{a+1}>0. \label{may}%
\end{align}
Since $\int_{y}^{\infty} z^{a}\exp\left(  \frac{d}{z}\right)  dz$ is
decreasing in $y$, there exists $M>c>0$, such that for $y>M$
\begin{align}
\int_{y}^{\infty} z^{a}\exp\left(  \frac{d}{z}\right)  dz  &  < \frac{-2}{a+1}
y^{a+1}\exp\left(  \frac{d}{y}\right)  . \label{niub}%
\end{align}
Substitute \eqref{niub} into \eqref{lop2s} with $x=\infty$
\begin{align}
v(\infty)  &  =\frac{2}{\xi^{2}}\int_{c}^{\infty} \frac{1}{y^{a+3}\exp\left(
\frac{d}{y}\right)  }\left(  \int_{y}^{\infty} z^{ a}\exp\left(  \frac{d}%
{z}\right)  dz\right)  dy\nonumber\\
&  =\frac{2}{\xi^{2}}\int_{c}^{M} \frac{\int_{y}^{\infty} z^{ a}\exp\left(
\frac{d}{z}\right)  dz}{y^{a+3}\exp\left(  \frac{d}{y}\right)  } dy+\frac
{2}{\xi^{2}}\int_{M}^{\infty} \frac{ \int_{y}^{\infty} z^{ a}\exp\left(
\frac{d}{z}\right)  dz}{y^{a+3}\exp\left(  \frac{d}{y}\right)  } dy\nonumber\\
&  <\frac{2}{\xi^{2}}\int_{c}^{M} \frac{ \int_{y}^{\infty} z^{ a}\exp\left(
\frac{d}{z}\right)  dz}{y^{a+3}\exp\left(  \frac{d}{y}\right)  } dy+\frac
{2}{\xi^{2}}\int_{M}^{\infty} \frac{1}{y^{a+3}\exp\left(  \frac{d}{y}\right)
} \frac{-2}{a+1} y^{a+1}\exp\left(  \frac{d}{y}\right)  dy\nonumber\\
&  =\frac{2}{\xi^{2}}\int_{c}^{M} \frac{1}{y^{a+3}\exp\left(  \frac{d}%
{y}\right)  }\left(  \int_{y}^{\infty} z^{ a}\exp\left(  \frac{d}{z}\right)
dz\right)  dy+\frac{-4}{(a+1)\xi^{2}}\int_{M}^{\infty} \frac{1}{y^{2}%
}dy\nonumber\\
&  =\frac{2}{\xi^{2}}\int_{c}^{M} \frac{1}{y^{a+3}\exp\left(  \frac{d}%
{y}\right)  }\int_{y}^{\infty} z^{ a}\exp\left(  \frac{d}{z}\right)  dz
dy+\frac{-4}{(a+1)\xi^{2} M}\nonumber\\
&  <\infty.\nonumber
\end{align}

From \eqref{may}, there exists $M^{\prime}>c>0$, such that for $y>M^{\prime}$
\begin{align}
\int_{y}^{\infty} z^{a}\exp\left(  \frac{d}{z}\right)  dz  &  > \frac
{-1}{2(a+1)} y^{a+1}\exp\left(  \frac{d}{y}\right)  . \label{niub22}%
\end{align}
Similarly substitute \eqref{niub22} into \eqref{lop41} with $x=\infty$
\begin{align}
v_{b}(\infty)  &  =\frac{2}{\xi^{2}}\int_{c}^{\infty} \frac{1}{y^{a+2}%
\exp\left(  \frac{d}{y}\right)  }\left(  \int_{y}^{\infty} z^{ a}\exp\left(
\frac{d}{z}\right)  dz\right)  dy\nonumber\\
&  \geq\frac{2}{\xi^{2}}\int_{M^{\prime}}^{\infty} \frac{1}{y^{a+2}\exp\left(
\frac{d}{y}\right)  }\left(  \int_{y}^{\infty} z^{ a}\exp\left(  \frac{d}%
{z}\right)  dz\right)  dy\nonumber\\
&  >\frac{2}{\xi^{2}}\int_{M^{\prime}}^{\infty} \frac{1}{y^{a+2}\exp\left(
\frac{d}{y}\right)  }\left(  \frac{-1}{2(a+1)} y^{a+1}\exp\left(  \frac{d}%
{y}\right)  \right)  dy\nonumber\\
&  =\frac{-1}{\xi^{2}(a+1)}\int_{M^{\prime}}^{\infty} \frac{1}{y}
dy\nonumber\\
&  =\infty.\nonumber
\end{align}

(ii) When $a\geq-1$, since $d>0$, we have that $\exp\left(  \frac{d}%
{y}\right)  \geq1$, for $y>c>0$. Then
\begin{align}
s(\infty)  &  =\frac{b}{c^{a}}\int_{c}^{\infty} y^{a}\exp\left(  \frac{d}%
{y}\right)  dy \geq\frac{b}{c^{a}}\int_{c}^{\infty} y^{a} dy=\infty.\nonumber
\end{align}
Thus $v(\infty)=\infty$ and $v_{b}(\infty)=\infty$ in this case. To summarize,
$v(r)<\infty$ if and only if $a<-1$, and $v_{b}(r)=\infty$ for $a\in
\mathbb{R}$.

For the left endpoint $\ell$
\begin{align}
v(0)  &  =\frac{2}{\xi^{2}}\int_{0}^{c} \frac{1}{y^{a+3}\exp\left(  \frac
{d}{y}\right)  }\left(  \int_{0}^{y} z^{ a}\exp\left(  \frac{d}{z}\right)
dz\right)  dy, \label{lop4}%
\end{align}
and
\begin{align}
v_{b}(0)  &  =\frac{2}{\xi^{2}}\int_{0}^{c} \frac{1}{y^{a+2}\exp\left(
\frac{d}{y}\right)  }\left(  \int_{0}^{y} z^{ a}\exp\left(  \frac{d}%
{z}\right)  dz\right)  dy. \label{lop51}%
\end{align}

For $0\leq z\leq y$, we have $e^{\frac{d}{z}}\geq e^{\frac{d}{y}}$, and
substitute this inequality into \eqref{lop4}
\begin{align}
v(0)  &  \geq\frac{2}{\xi^{2}}\int_{0}^{c} \frac{1}{y^{a+3}\exp\left(
\frac{d}{y}\right)  }\left(  \int_{0}^{y} z^{ a} dz\right)  \exp\left(
\frac{d}{y}\right)  dy =\frac{2}{(a+1)\xi^{2}}\int_{0}^{c} \frac{1}{y^{2}}dy
=\infty.\nonumber
\end{align}

Similarly substitute this inequality into \eqref{lop51}
\begin{align}
v_{b}(0)  &  \geq\frac{2}{\xi^{2}}\int_{0}^{c} \frac{1}{y^{a+2}\exp\left(
\frac{d}{y}\right)  }\left(  \int_{0}^{y} z^{ a}dz\right)  \exp\left(
\frac{d}{y}\right)  dy\nonumber\\
&  =\frac{2}{(a+1)\xi^{2}}\int_{0}^{c} \frac{1}{y}dy=\infty.\nonumber
\end{align}
To summarize, $v(\ell)=\infty$ and $v_{b}(\ell)=\infty$ for $a\in\mathbb{R}$.
From \eqref{sth}, the above proofs also work for the case of $\widetilde{v}$
by substituting $a$ for $\widetilde{a}$. \qed

\subsection{Proof of Proposition \ref{SZMMart}\label{SZMMartapp}}

\proof
Divide into three cases:

(i) When $\alpha>0$, $\widetilde{s}(\infty)=\infty$, then $\widetilde
{v}(\infty)=\infty$ and $\widetilde{v}_{b}(\infty)=\infty$.

(ii) When $\alpha=0$
\begin{align}
\widetilde{v}(x)  &  =\frac{1}{\kappa\theta} \int_{c}^{x} \left(  1-
e^{-\frac{2\kappa\theta}{\gamma^{2}}(x-y)} \right)  dy =\frac{1}{\kappa\theta
}\left(  x+\frac{\gamma^{2}}{2\kappa\theta}e^{ \frac{2\kappa\theta}{\gamma
^{2}}(c-x)} -c-\frac{\gamma^{2}}{2\kappa\theta}\right)  .\nonumber
\end{align}
Then $\widetilde{v}(\infty)=\infty$. Similarly we can compute
\begin{align}
\widetilde{v}_{b}(x)  &  =\frac{1}{\kappa\theta} \int_{c}^{x} y^{2} \left(  1-
e^{-\frac{2\kappa\theta}{\gamma^{2}}(x-y)} \right)  dy =\frac{1}{3\kappa
\theta}x^{3} - e^{-\frac{2\kappa\theta}{\gamma^{2}}x}\int_{c}^{x} y^{2}
e^{\frac{2\kappa\theta}{\gamma^{2}} y} dy -\frac{c^{3}}{3\kappa\theta
}.\nonumber
\end{align}
Since $\int_{c}^{x} y^{2} e^{\frac{2\kappa\theta}{\gamma^{2}} y} dy\leq
\int_{c}^{x} x^{2} e^{\frac{2\kappa\theta}{\gamma^{2}} y} dy$, then
\begin{align}
\widetilde{v}_{b}(x)  &  \geq\frac{1}{3\kappa\theta}x^{3} - e^{-\frac
{2\kappa\theta}{\gamma^{2}}x}\int_{c}^{x} x^{2} e^{\frac{2\kappa\theta}%
{\gamma^{2}} y} dy -\frac{c^{3}}{3\kappa\theta} =\frac{1}{3\kappa\theta}x^{3}
- \frac{\gamma^{2}}{2\kappa\theta} x^{2} (1- e^{\frac{2\kappa\theta}%
{\gamma^{2}}(c-x)}) -\frac{c^{3}}{3\kappa\theta}. \label{rgh}%
\end{align}
Then $\widetilde{v}_{b}(\infty)=\infty$ can be verified, because the right
hand side of \eqref{rgh} tends to $\infty$ as $x\rightarrow\infty$.

(iii) When $\alpha<0$, the test function is
\begin{align}
\widetilde{v}(x)  &  =\frac{2}{\gamma^{2}} \int_{c}^{x} \frac{\int_{y}^{x}
\frac{\alpha}{\gamma^{2}}\left(  z-\frac{\kappa\theta}{\alpha}\right)  ^{2} dz
}{e^{\frac{\alpha}{\gamma^{2}}(y-\frac{\kappa\theta}{\alpha})^{2}}}dy=\frac
{2}{\gamma^{2}} \int_{c}^{x} e^{-\frac{\alpha}{\gamma^{2}}(y-\frac
{\kappa\theta}{\alpha})^{2}}\left(  \int_{y}^{x} \frac{\alpha}{\gamma^{2}%
}\left(  z-\frac{\kappa\theta}{\alpha}\right)  ^{2} dz\right)  dy.\nonumber
\end{align}
Then
\begin{align}
\widetilde{v}(\infty)  &  =\frac{2}{\gamma^{2}} \int_{c}^{\infty}
e^{-\frac{\alpha}{\gamma^{2}}(y-\frac{\kappa\theta}{\alpha})^{2}}\left(
\int_{y}^{\infty} \frac{\alpha}{\gamma^{2}}\left(  z-\frac{\kappa\theta
}{\alpha}\right)  ^{2} dz\right)  dy. \label{sze1}%
\end{align}

Since $\alpha<0$ is assumed here, then $\lim_{y\rightarrow\infty} y^{-1}
e^{\frac{\alpha}{\gamma^{2}}(y-\frac{\kappa\theta}{\alpha})^{2}} =0$, and we
can apply L'Hôpital's rule
\begin{align}
\lim_{y\rightarrow\infty} \frac{ \int_{y}^{\infty} \frac{\alpha}{\gamma^{2}%
}\left(  z-\frac{\kappa\theta}{\alpha}\right)  ^{2} dz }{y^{-1} e^{\frac
{\alpha}{\gamma^{2}}(y-\frac{\kappa\theta}{\alpha})^{2}}}  &  =\frac
{-\gamma^{2}}{2\alpha}>0.\nonumber
\end{align}
So as $y\rightarrow\infty$, there exists $M>c>0$, such that for $y>M$
\begin{align}
\int_{y}^{\infty} \frac{\alpha}{\gamma^{2}}\left(  z-\frac{\kappa\theta
}{\alpha}\right)  ^{2} dz> \frac{-\gamma^{2}}{4\alpha} y^{-1} e^{\frac{\alpha
}{\gamma^{2}}(y-\frac{\kappa\theta}{\alpha})^{2}}. \label{sze}%
\end{align}
Substitute \eqref{sze} into \eqref{sze1}
\begin{align}
\widetilde{v}(\infty)  &  \geq\frac{2}{\gamma^{2}} \int_{M}^{\infty}
e^{-\frac{\alpha}{\gamma^{2}}(y-\frac{\kappa\theta}{\alpha})^{2}}\left(
\int_{y}^{\infty} \frac{\alpha}{\gamma^{2}}\left(  z-\frac{\kappa\theta
}{\alpha}\right)  ^{2} dz\right)  dy\nonumber\\
&  >\frac{2}{\gamma^{2}} \int_{M}^{\infty} e^{-\frac{\alpha}{\gamma^{2}%
}(y-\frac{\kappa\theta}{\alpha})^{2}}\left(  \frac{-\gamma^{2}}{4\alpha}
y^{-1} e^{\frac{\alpha}{\gamma^{2}}(y-\frac{\kappa\theta}{\alpha})^{2}}
\right)  dy\nonumber\\
&  =\frac{-1}{2\alpha} \int_{M}^{\infty} y^{-1} dy =\infty.\nonumber
\end{align}
Thus $\widetilde{v}(\infty)=\infty$ in this case. Similarly, we can compute
\begin{align}
\widetilde{v}_{b}(\infty)  &  =\frac{2}{\gamma^{2}} \int_{c}^{\infty} y^{2}
e^{-\frac{\alpha}{\gamma^{2}}(y-\frac{\kappa\theta}{\alpha})^{2}}\left(
\int_{y}^{\infty} \frac{\alpha}{\gamma^{2}}\left(  z-\frac{\kappa\theta
}{\alpha}\right)  ^{2} dz\right)  dy. \label{sze2}%
\end{align}
With the same $M$ as above, substitute \eqref{sze} into \eqref{sze2}
\begin{align}
\widetilde{v}_{b}(\infty)  &  \geq\frac{2}{\gamma^{2}} \int_{M}^{\infty} y^{2}
e^{-\frac{\alpha}{\gamma^{2}}(y-\frac{\kappa\theta}{\alpha})^{2}}\left(
\int_{y}^{\infty} \frac{\alpha}{\gamma^{2}}\left(  z-\frac{\kappa\theta
}{\alpha}\right)  ^{2} dz\right)  dy\nonumber\\
&  >\frac{2}{\gamma^{2}} \int_{c}^{\infty} y^{2} e^{-\frac{\alpha}{\gamma^{2}%
}(y-\frac{\kappa\theta}{\alpha})^{2}}\left(  \frac{-\gamma^{2}}{4\alpha}
y^{-1} e^{\frac{\alpha}{\gamma^{2}}(y-\frac{\kappa\theta}{\alpha})^{2}}
\right)  dy\nonumber\\
&  =\frac{-1}{2\alpha} \int_{M}^{\infty} y dy =\infty.\nonumber
\end{align}
Thus $\widetilde{v}_{b}(\infty)=\infty$ in this case.

We then consider the case of the left endpoint $\ell$. From the definition of
$\widetilde s(\cdot)$, we have that $\widetilde s (0)>-\infty$ for $\alpha
\in\mathbb{R}$.

Similarly as above, we consider the following two cases:

(i) When $\alpha=0$, $\widetilde{v}(0)=\frac{1}{\kappa\theta}\left(
\frac{\gamma^{2}}{2\kappa\theta}e^{ \frac{2\kappa\theta}{\gamma^{2}}(c)}
-c-\frac{\gamma^{2}}{2\kappa\theta}\right)  <\infty.$

(ii) When $\alpha\neq0$,
\begin{align}
\widetilde{v}(0)  &  =\frac{2}{\gamma^{2}} \int_{0}^{c} e^{-\frac{\alpha
}{\gamma^{2}}(y-\frac{\kappa\theta}{\alpha})^{2}}\left(  \int_{0}^{y}
\frac{\alpha}{\gamma^{2}}\left(  z-\frac{\kappa\theta}{\alpha}\right)  ^{2}
dz\right)  dy. \label{sze3}%
\end{align}
Since $\lim_{y\rightarrow0} y e^{\frac{\alpha}{\gamma^{2}}(y-\frac
{\kappa\theta}{\alpha})^{2}}=0$, we can apply L'Hôpital's rule
\begin{align}
\lim_{y\rightarrow0} \frac{ \int_{0}^{y} \frac{\alpha}{\gamma^{2}}\left(
z-\frac{\kappa\theta}{\alpha}\right)  ^{2} dz }{y e^{\frac{\alpha}{\gamma^{2}%
}(y-\frac{\kappa\theta}{\alpha})^{2}}}  &  =1.\nonumber
\end{align}

So as $y\rightarrow0$, there exists $0<\varepsilon<c$, such that for $0\leq
y<\varepsilon$
\begin{align}
\int_{0}^{y} \frac{\alpha}{\gamma^{2}}\left(  z-\frac{\kappa\theta}{\alpha
}\right)  ^{2} dz  &  < 2 y e^{\frac{\alpha}{\gamma^{2}}(y-\frac{\kappa\theta
}{\alpha})^{2}}. \label{sze4}%
\end{align}
Substitute \eqref{sze4} into \eqref{sze3}
\begin{align}
\widetilde{v}(0)  &  =\frac{2}{\gamma^{2}} \int_{0}^{\varepsilon}
e^{-\frac{\alpha}{\gamma^{2}}(y-\frac{\kappa\theta}{\alpha})^{2}}\left(
\int_{0}^{y} \frac{\alpha}{\gamma^{2}}\left(  z-\frac{\kappa\theta}{\alpha
}\right)  ^{2} dz\right)  dy\nonumber\\
&  \quad\quad+\frac{2}{\gamma^{2}} \int_{\varepsilon}^{c} e^{-\frac{\alpha
}{\gamma^{2}}(y-\frac{\kappa\theta}{\alpha})^{2}}\left(  \int_{0}^{y}
\frac{\alpha}{\gamma^{2}}\left(  z-\frac{\kappa\theta}{\alpha}\right)  ^{2}
dz\right)  dy\nonumber\\
&  < \frac{2}{\gamma^{2}} \int_{0}^{\varepsilon} e^{-\frac{\alpha}{\gamma^{2}%
}(y-\frac{\kappa\theta}{\alpha})^{2}}\left(  2 y e^{\frac{\alpha}{\gamma^{2}%
}(y-\frac{\kappa\theta}{\alpha})^{2}}\right)  dy\nonumber\\
&  \quad\quad+\frac{2}{\gamma^{2}} \int_{\varepsilon}^{c} e^{-\frac{\alpha
}{\gamma^{2}}(y-\frac{\kappa\theta}{\alpha})^{2}}\left(  \int_{0}^{y}
\frac{\alpha}{\gamma^{2}}\left(  z-\frac{\kappa\theta}{\alpha}\right)  ^{2}
dz\right)  dy\nonumber
\end{align}
Then,
\begin{align}
\widetilde{v}(\infty)  &  <\frac{2}{\gamma^{2}} \int_{0}^{\varepsilon} 2 y
dy+\frac{2}{\gamma^{2}} \int_{\varepsilon}^{c} e^{-\frac{\alpha}{\gamma^{2}%
}(y-\frac{\kappa\theta}{\alpha})^{2}}\left(  \int_{0}^{y} \frac{\alpha}%
{\gamma^{2}}\left(  z-\frac{\kappa\theta}{\alpha}\right)  ^{2} dz\right)
dy\nonumber\\
&  =\frac{2\varepsilon^{2}}{\gamma^{2}}+\frac{2}{\gamma^{2}} \int
_{\varepsilon}^{c} e^{-\frac{\alpha}{\gamma^{2}}(y-\frac{\kappa\theta}{\alpha
})^{2}}\left(  \int_{0}^{y} \frac{\alpha}{\gamma^{2}}\left(  z-\frac
{\kappa\theta}{\alpha}\right)  ^{2} dz\right)  dy\nonumber\\
&  <\infty. \label{sze5}%
\end{align}
To summarize, $\widetilde{v}(\ell)<\infty$ for $\alpha\in\mathbb{R}$.
Similarly, when $\alpha=0$, $\widetilde{v}_{b}(0)=\int_{0}^{c} y^{2}
e^{\frac{2\kappa\theta}{\gamma^{2}}y}dy-\frac{c^{3}}{3\kappa\theta}<\infty.$

When $\alpha\neq0$,
\begin{align}
\widetilde{v}_{b}(0)  &  =\frac{2}{\gamma^{2}} \int_{0}^{c} y^{2}
e^{-\frac{\alpha}{\gamma^{2}}(y-\frac{\kappa\theta}{\alpha})^{2}}\left(
\int_{0}^{y} \frac{\alpha}{\gamma^{2}}\left(  z-\frac{\kappa\theta}{\alpha
}\right)  ^{2} dz\right)  dy. \label{sze6}%
\end{align}
Substitute \eqref{sze4} into \eqref{sze6}, and use the same $\varepsilon$ as
above. For $0\leq y<\varepsilon$
\begin{align}
\widetilde{v}_{b}(0)  &  =\frac{2}{\gamma^{2}} \int_{0}^{\varepsilon} y^{2}
e^{-\frac{\alpha}{\gamma^{2}}(y-\frac{\kappa\theta}{\alpha})^{2}}\left(
\int_{0}^{y} \frac{\alpha}{\gamma^{2}}\left(  z-\frac{\kappa\theta}{\alpha
}\right)  ^{2} dz\right)  dy\nonumber\\
&  \quad\quad+\frac{2}{\gamma^{2}} \int_{\varepsilon}^{c} y^{2} e^{-\frac
{\alpha}{\gamma^{2}}(y-\frac{\kappa\theta}{\alpha})^{2}}\left(  \int_{0}^{y}
\frac{\alpha}{\gamma^{2}}\left(  z-\frac{\kappa\theta}{\alpha}\right)  ^{2}
dz\right)  dy\nonumber\\
&  < \frac{2}{\gamma^{2}} \int_{0}^{\varepsilon} y^{2} e^{-\frac{\alpha
}{\gamma^{2}}(y-\frac{\kappa\theta}{\alpha})^{2}}\left(  2 y e^{\frac{\alpha
}{\gamma^{2}}(y-\frac{\kappa\theta}{\alpha})^{2}}\right)  dy\nonumber\\
&  \quad\quad+\frac{2}{\gamma^{2}} \int_{\varepsilon}^{c} y^{2} e^{-\frac
{\alpha}{\gamma^{2}}(y-\frac{\kappa\theta}{\alpha})^{2}}\left(  \int_{0}^{y}
\frac{\alpha}{\gamma^{2}}\left(  z-\frac{\kappa\theta}{\alpha}\right)  ^{2}
dz\right)  dy\nonumber\\
&  =\frac{2}{\gamma^{2}} \int_{0}^{\varepsilon} 2 y^{3} dy+\frac{2}{\gamma
^{2}} \int_{\varepsilon}^{c} y^{2} e^{-\frac{\alpha}{\gamma^{2}}%
(y-\frac{\kappa\theta}{\alpha})^{2}}\left(  \int_{0}^{y} \frac{\alpha}%
{\gamma^{2}}\left(  z-\frac{\kappa\theta}{\alpha}\right)  ^{2} dz\right)
dy\nonumber\\
&  =\frac{\varepsilon^{4}}{\gamma^{2}}+\frac{2}{\gamma^{2}} \int_{\varepsilon
}^{c} y^{2} e^{-\frac{\alpha}{\gamma^{2}}(y-\frac{\kappa\theta}{\alpha})^{2}%
}\left(  \int_{0}^{y} \frac{\alpha}{\gamma^{2}}\left(  z-\frac{\kappa\theta
}{\alpha}\right)  ^{2} dz\right)  dy\nonumber\\
&  <\infty. \label{sze7}%
\end{align}
To summarize, $\widetilde{v}_{b}(\ell)<\infty$, for $\alpha\in\mathbb{R}%
$.\qed

\subsection{Proof of Proposition \ref{hwmarthm} \label{hwapp}}

\proof
We distinguish three situations:

(I) $\mu>\frac{1}{2}\sigma^{2}$. Apply a change of variable $z=\sqrt{y}$. Then
$y=z^{2}$, $dy=2zdz$, and
\begin{align}
\widetilde s (x)  &  =2C_{1} \int_{\sqrt{c}}^{\sqrt{x}} z^{1-\frac{4\mu
}{\sigma^{2}}} e^{-\frac{2\rho}{\sigma}z}dz=2C_{1} \int_{\sqrt{c}}^{\sqrt{x}}
z^{-\alpha} e^{-\gamma z}dz,\quad x\in\bar{J}. \label{schw3}%
\end{align}

Note that the function in \eqref{schw3} is similar to the scale function in
\eqref{sches}, except that there is a $\sqrt{x}$ in place of $x$. From
\eqref{schw3}
\begin{align}
\widetilde s (\infty)  &  =2C_{1} \int_{\sqrt{c}}^{\infty} z^{-\alpha}
e^{-\gamma z}dz.\nonumber
\end{align}
From the property of the gamma function
\begin{align}
\widetilde{s}(\infty)  &
\begin{cases}
<\infty,\quad & \text{if $\gamma\geq0$,}\notag\\
=\infty,\quad & \text{if $\gamma<0$.}%
\end{cases}
\end{align}
Divide into three cases based on $\gamma$:

(i) When $\gamma<0$, $\widetilde{s}(\infty)=\infty$, then $\widetilde
v(\infty)=\infty$ and $\widetilde v_{b}(\infty)=\infty$.

(ii) When $\gamma=0$, $\widetilde v(\infty)$ and $\widetilde v_{b}(\infty)$
can be simplified and
\begin{align}
\widetilde v(\infty)  &  =\frac{2}{\sigma^{2}}\int_{c}^{\infty} y^{\frac
{\alpha-3}{2}} \left(  \int_{y}^{\infty} z^{-\frac{\alpha+1}{2}} dz\right)
dy=\frac{4}{\sigma^{2}(\alpha-1)}\int_{c}^{\infty} y^{-1}dy=\infty,\nonumber
\end{align}
and
\begin{align}
\widetilde v_{b}(\infty)=\frac{2}{\sigma^{2}}\int_{c}^{\infty} y^{\frac
{\alpha-1}{2}} \left(  \int_{y}^{\infty} z^{-\frac{\alpha+1}{2}} dz\right)
dy=\int_{c}^{\infty} \frac{4}{\sigma^{2}(\alpha-1)} dy=\infty.\nonumber
\end{align}

(iii) When $\gamma> 0$, from \eqref{vhw}
\begin{align}
\widetilde v(\infty)  &  =\frac{2}{\sigma^{2}}\int_{c}^{\infty} y^{\frac
{\alpha-3}{2}} e^{\gamma\sqrt{y}} \left(  \int_{y}^{\infty} z^{-\frac
{\alpha+1}{2}} e^{-\gamma\sqrt{z}}dz\right)  dy, \label{rh}%
\end{align}
and
\begin{align}
\widetilde v_{b}(\infty)  &  =\frac{2}{\sigma^{2}}\int_{c}^{\infty}
y^{\frac{\alpha-1}{2}} e^{\gamma\sqrt{y}} \left(  \int_{y}^{\infty}
z^{-\frac{\alpha+1}{2}} e^{-\gamma\sqrt{z}}dz\right)  dy. \label{rh2}%
\end{align}

Since $\alpha>1$, then $\lim_{y\rightarrow\infty} y^{-\frac{\alpha}{2}}
e^{-\gamma\sqrt{y}}=0$, and from L'Hôpital's rule
\begin{align}
\lim_{y\rightarrow\infty}\frac{ \int_{y}^{\infty} z^{-\frac{\alpha+1}{2}}
e^{-\gamma\sqrt{z}}dz }{y^{-\frac{\alpha}{2}} e^{-\gamma\sqrt{y}}}  &
=\lim_{y\rightarrow\infty} \frac{1}{\frac{\alpha}{2}y^{-1/2}+\frac{\gamma}{2}%
}=\frac{2}{\gamma}>0.\nonumber
\end{align}
As $y\rightarrow\infty$
\begin{align}
\int_{y}^{\infty} z^{-\frac{\alpha+1}{2}} e^{-\gamma\sqrt{z}}dz  &  \sim
\frac{2}{\gamma} y^{-\frac{\alpha}{2}} e^{-\gamma\sqrt{y}}. \label{ro1}%
\end{align}

From \eqref{ro1}, there exists $0<M<\infty$, such that for $y>M$
\begin{align}
\int_{y}^{\infty} z^{-\frac{\alpha+1}{2}} e^{-\gamma\sqrt{z}}dz  &  < \frac
{4}{\gamma} y^{-\frac{\alpha}{2}} e^{-\gamma\sqrt{y}}. \label{ro2}%
\end{align}
Substitute \eqref{ro2} into \eqref{rh}
\begin{align}
\widetilde v(\infty)  &  =\frac{2}{\sigma^{2}}\int_{c}^{\infty} y^{\frac
{\alpha-3}{2}} e^{\gamma\sqrt{y}} \left(  \int_{y}^{\infty} z^{-\frac
{\alpha+1}{2}} e^{-\gamma\sqrt{z}}dz\right)  dy\nonumber\\
&  =\frac{2}{\sigma^{2}}\int_{c}^{M} y^{\frac{\alpha-3}{2}} e^{\gamma\sqrt{y}}
\left(  \int_{y}^{\infty} z^{-\frac{\alpha+1}{2}} e^{-\gamma\sqrt{z}%
}dz\right)  dy+\frac{2}{\sigma^{2}}\int_{M}^{\infty} y^{\frac{\alpha-3}{2}}
e^{\gamma\sqrt{y}} \left(  \int_{y}^{\infty} z^{-\frac{\alpha+1}{2}}
e^{-\gamma\sqrt{z}}dz\right)  dy\nonumber\\
&  <\frac{2}{\sigma^{2}}\int_{c}^{M} y^{\frac{\alpha-3}{2}} e^{\gamma\sqrt{y}}
\left(  \int_{y}^{\infty} z^{-\frac{\alpha+1}{2}} e^{-\gamma\sqrt{z}%
}dz\right)  dy+\frac{2}{\sigma^{2}}\int_{M}^{\infty} y^{\frac{\alpha-3}{2}}
e^{\gamma\sqrt{y}} \left(  \frac{4}{\gamma} y^{-\frac{\alpha}{2}}
e^{-\gamma\sqrt{y}}\right)  dy\nonumber\\
&  =\frac{2}{\sigma^{2}}\int_{c}^{M} y^{\frac{\alpha-3}{2}} e^{\gamma\sqrt{y}}
\left(  \int_{y}^{\infty} z^{-\frac{\alpha+1}{2}} e^{-\gamma\sqrt{z}%
}dz\right)  dy +\frac{16}{\sqrt{M} \gamma\sigma^{2}}\nonumber\\
&  <\infty.\nonumber
\end{align}
Then $\widetilde v(\infty)<\infty$, for $\gamma> 0$.

From \eqref{ro1}, there exists $0<c<M^{\prime}<\infty$, such that for
$y>M^{\prime}$
\begin{align}
\int_{y}^{\infty} z^{-\frac{\alpha+1}{2}} e^{-\gamma\sqrt{z}}dz  &  > \frac
{1}{\gamma} y^{-\frac{\alpha}{2}} e^{-\gamma\sqrt{y}}. \label{ro3}%
\end{align}

Substitute \eqref{ro3} into \eqref{rh2}
\begin{align}
\widetilde v(\infty)  &  =\frac{2}{\sigma^{2}}\int_{c}^{\infty} y^{\frac
{\alpha-3}{2}} e^{\gamma\sqrt{y}} \left(  \int_{y}^{\infty} z^{-\frac
{\alpha+1}{2}} e^{-\gamma\sqrt{z}}dz\right)  dy\nonumber\\
&  \geq\frac{2}{\sigma^{2}}\int_{M^{\prime}}^{\infty} y^{\frac{\alpha-3}{2}}
e^{\gamma\sqrt{y}} \left(  \int_{y}^{\infty} z^{-\frac{\alpha+1}{2}}
e^{-\gamma\sqrt{z}}dz\right)  dy\nonumber\\
&  > \frac{2}{\sigma^{2}}\int_{M^{\prime}}^{\infty} y^{\frac{\alpha-3}{2}}
e^{\gamma\sqrt{y}} \left(  \frac{1}{\gamma} y^{-\frac{\alpha}{2}}
e^{-\gamma\sqrt{y}}\right)  dy\nonumber\\
&  =\frac{2}{\gamma\sigma^{2}}\int_{M^{\prime}}^{\infty} y^{-1} dy
=\infty.\nonumber
\end{align}
Then $\widetilde v_{b}(\infty)=\infty$, for $\gamma> 0$.

We now look at the case of the left boundary $\ell$. From \eqref{schw3}
\begin{align}
\widetilde s (0)  &  =-2C_{1} \int_{0}^{\sqrt{c}} z^{-\alpha} e^{-\gamma
z}dz.\nonumber
\end{align}
When $\gamma>0$, since $\alpha>1$, from the property of the gamma function, we
have $\widetilde s(0)=-\infty$. When $\gamma\leq0$, then $e^{-\gamma z}\geq1$,
and
\begin{align}
\widetilde s (0)  &  =-2C_{1} \int_{0}^{\sqrt{c}} z^{-\alpha} e^{-\gamma z}dz
\leq-2C_{1} \int_{0}^{\sqrt{c}} z^{-\alpha} dz=-\infty.\nonumber
\end{align}
To summarize, $\widetilde s(0)=-\infty$ for $\gamma\in\mathbb{R}$. Then
$\widetilde v(0)=\infty$ and $\widetilde v_{b}(0)=\infty$ hold.

(II) $\mu=\frac{1}{2}\sigma^{2}$. We consider the case when $\alpha=1$. Then
\begin{align}
\widetilde s (\infty)  &  =2C_{1} \int_{\sqrt{c}}^{\infty} z^{-1} e^{-\gamma
z}dz,\nonumber
\end{align}
Divide into two cases based on the value of $\gamma$. If $\gamma\leq0$, then
$e^{-\gamma z}\geq1$, and
\begin{align}
\widetilde s (\infty)  &  \geq2C_{1} \int_{\sqrt{c}}^{\infty} z^{-1}
dz=\infty.\nonumber
\end{align}
Then in this case, $\widetilde v(r)=\infty$ and $\widetilde v_{b}(r)=\infty$.

If $\gamma>0$, from properties of the gamma function, $\widetilde
s(\infty)<\infty$. To summarize, when $\alpha=1$
\begin{align}
\widetilde{s}(\infty)  &
\begin{cases}
=\infty, \quad & \text{if $\gamma\leq0$},\notag\\
<\infty, \quad & \text{if $\gamma>0$}.
\end{cases}
\end{align}
Similarly for the case of the left boundary $\ell$. If $\gamma>0$, from the
properties of the gamma function, $\widetilde{s}(0)=-\infty$. If $\gamma\leq
0$, then $e^{-\gamma z}\geq1$, and
\begin{align}
\widetilde s (0)  &  \leq-2C_{1} \int_{0}^{\sqrt{c}} z^{-1} dz=-\infty
.\nonumber
\end{align}
To summarize, when $\alpha=1$, we have $\widetilde s(\ell)=-\infty$, then
$\widetilde v(\ell)=\infty$ and $\widetilde v_{b}(\ell)=\infty$.

Consider the case when $\alpha=1$ and $\gamma>0$, from the above result, there
is $\widetilde s(\infty)<\infty$, and we study the properties of $\widetilde
v(\infty)$ and $\widetilde v_{b} (\infty)$. From the definition in
\eqref{vhw}
\begin{align}
\widetilde v(\infty)  &  =\frac{2}{\sigma^{2}}\int_{c}^{\infty} y^{-1}
e^{\gamma\sqrt{y}} \left(  \int_{y}^{\infty} z^{-1} e^{-\gamma\sqrt{z}%
}dz\right)  dy. \label{vhw1}%
\end{align}
Since $\gamma>0$, then $\lim_{y\rightarrow\infty} y^{-\frac{1}{2}}
e^{-\gamma\sqrt{y}}=0$, and from L'Hôpital's rule
\begin{align}
\lim_{y\rightarrow\infty}\frac{ \int_{y}^{\infty} z^{-1} e^{-\gamma\sqrt{z}}dz
}{y^{-\frac{1}{2}} e^{-\gamma\sqrt{y}}}  &  =\lim_{y\rightarrow\infty}
\frac{1}{\frac{1}{2}y^{-1/2}+\frac{\gamma}{2}}=\frac{2}{\gamma}>0.
\label{sihw}%
\end{align}
As $y\rightarrow\infty$, $\int_{y}^{\infty} z^{-1} e^{-\gamma\sqrt{z}}dz
\sim\frac{2}{\gamma} y^{-\frac{1}{2}} e^{-\gamma\sqrt{y}}$, and thus there
exists $M<\infty$, such that for $y>M$
\begin{align}
\int_{y}^{\infty} z^{-1} e^{-\gamma\sqrt{z}}dz  &  < \frac{4}{\gamma}
y^{-\frac{1}{2}} e^{-\gamma\sqrt{y}}. \label{mhw}%
\end{align}
Substitute \eqref{mhw} into \eqref{vhw1}
\begin{align}
\widetilde v(\infty)  &  =\frac{2}{\sigma^{2}}\int_{c}^{M} y^{-1}
e^{\gamma\sqrt{y}} \left(  \int_{y}^{\infty} z^{-1} e^{-\gamma\sqrt{z}%
}dz\right)  dy+\frac{2}{\sigma^{2}}\int_{M}^{\infty} y^{-1} e^{\gamma\sqrt{y}}
\left(  \int_{y}^{\infty} z^{-1} e^{-\gamma\sqrt{z}}dz\right)  dy\nonumber\\
&  <\frac{2}{\sigma^{2}}\int_{c}^{M} y^{-1} e^{\gamma\sqrt{y}} \left(
\int_{y}^{\infty} z^{-1} e^{-\gamma\sqrt{z}}dz\right)  dy+\frac{2}{\sigma^{2}%
}\int_{M}^{\infty} y^{-1} e^{\gamma\sqrt{y}} \left(  \frac{4}{\gamma}
y^{-\frac{1}{2}} e^{-\gamma\sqrt{y}} \right)  dy\nonumber\\
&  =\frac{2}{\sigma^{2}}\int_{c}^{M} y^{-1} e^{\gamma\sqrt{y}} \left(
\int_{y}^{\infty} z^{-1} e^{-\gamma\sqrt{z}}dz\right)  dy+\frac{8}%
{\gamma\sigma^{2}}\int_{M}^{\infty} y^{-\frac{3}{2}} dy\nonumber\\
&  <\infty.\nonumber
\end{align}
From the definition in \eqref{vhw}
\begin{align}
\widetilde v_{b}(\infty)  &  =\frac{2}{\sigma^{2}}\int_{c}^{\infty}
e^{\gamma\sqrt{y}} \left(  \int_{y}^{\infty} z^{-1} e^{-\gamma\sqrt{z}%
}dz\right)  dy, \label{vhw2}%
\end{align}
From \eqref{sihw}, there exits $M^{\prime}>c>0$, such that for $y>M^{\prime}$
\begin{align}
\int_{y}^{\infty} z^{-1} e^{-\gamma\sqrt{z}}dz  &  > \frac{1}{\gamma}
y^{-\frac{1}{2}} e^{-\gamma\sqrt{y}}. \label{mhw2}%
\end{align}
Substitute \eqref{mhw2} into \eqref{vhw2}
\begin{align}
\widetilde v_{b}(\infty)  &  \geq\frac{2}{\sigma^{2}}\int_{M^{\prime}}%
^{\infty} e^{\gamma\sqrt{y}} \left(  \int_{y}^{\infty} z^{-1} e^{-\gamma
\sqrt{z}}dz\right)  dy\nonumber\\
&  >\frac{2}{\sigma^{2}}\int_{M^{\prime}}^{\infty} e^{\gamma\sqrt{y}} \left(
\frac{1}{\gamma} y^{-\frac{1}{2}} e^{-\gamma\sqrt{y}}\right)  dy =\frac
{2}{\sigma^{2}}\int_{M^{\prime}}^{\infty} y^{-\frac{1}{2}} dy =\infty
.\nonumber
\end{align}

(III) $\mu<\frac{1}{2}\sigma^{2}$. We consider the case when $\alpha<1$. Since
$-\frac{\alpha+1}{2}>-1$, then from the property of the gamma function
\begin{align}
\widetilde s(0)  &  =-C_{1} \int_{0}^{c} y^{-\frac{\alpha+1}{2}}
e^{-\gamma\sqrt{y}}dy>-\infty.\nonumber
\end{align}
From \eqref{schw3}, we have $\widetilde{s}(\infty)=2 C_{1} \int_{\sqrt{c}%
}^{\infty} z^{-\alpha} e^{-\gamma z}dz$, and divide into three cases. If
$\gamma>0$, then from the property of gamma function, $\widetilde{s}%
(\infty)<\infty$. If $\gamma\leq0$, then $e^{-\gamma z}\geq1$, and
$\widetilde{s}(\infty)\geq2 C_{1} \int_{\sqrt{c}}^{\infty} z^{-\alpha}
dz=\infty$. To summarize, when $\alpha<1$
\begin{align}
\widetilde{s}(\infty)  &
\begin{cases}
=\infty, \quad & \text{if $\gamma\leq0$},\notag\\
<\infty, \quad & \text{if $\gamma>0$}.
\end{cases}
\end{align}
We first look at $\widetilde v(0)$ and $\widetilde v_{b} (0)$. From the
definition in \eqref{vhw}
\begin{align}
\widetilde v(0)  &  =\frac{2}{\sigma^{2}} \int_{0}^{c} y^{\frac{\alpha-3}{2}}
e^{\gamma\sqrt{y}} \left(  \int_{0}^{y} z^{-\frac{\alpha+1}{2}} e^{-\gamma
\sqrt{z}}dz\right)  dy, \label{6hw}%
\end{align}
and
\begin{align}
\widetilde v_{b}(0)  &  =\frac{2}{\sigma^{2}} \int_{0}^{c} y^{\frac{\alpha
-1}{2}} e^{\gamma\sqrt{y}} \left(  \int_{0}^{y} z^{-\frac{\alpha+1}{2}}
e^{-\gamma\sqrt{z}}dz\right)  dy. \label{6bhw}%
\end{align}

Divide into two cases based on $\gamma$. When $\gamma\leq0$, $e^{-\gamma
\sqrt{z}}\geq1$, then
\begin{align}
\widetilde v(0)  &  \geq\frac{2}{\sigma^{2}} \int_{0}^{c} y^{\frac{\alpha
-3}{2}} e^{\gamma\sqrt{y}} \left(  \int_{0}^{y} z^{-\frac{\alpha+1}{2}}
dz\right)  dy\nonumber\\
&  =\frac{2}{\sigma^{2}} \int_{0}^{c} y^{\frac{\alpha-3}{2}} e^{\gamma\sqrt
{y}} \left(  \frac{2}{1-\alpha} y^{\frac{1-\alpha}{2}} \right)  dy\nonumber\\
&  =\frac{4}{\sigma^{2}(1-\alpha)} \int_{0}^{c} y^{-1} e^{\gamma\sqrt{y}}
dy.\nonumber
\end{align}
Apply a change of variable $z=\sqrt{y}$, then
\begin{align}
\widetilde v(0)  &  \geq\frac{4}{\sigma^{2}(1-\alpha)} \int_{0}^{c} y^{-1}
e^{\gamma\sqrt{y}} dy=\frac{8}{\sigma^{2}(1-\alpha)} \int_{0}^{\sqrt{c}}
z^{-1} e^{\gamma z} dz=\infty.\nonumber
\end{align}
The last equality is from the property of the gamma function. Then $\widetilde
v(0)=\infty$ holds when $\gamma\leq0$.

Since $\gamma\leq0$ is assumed, then $e^{-\gamma\sqrt{z}}\leq e^{-\gamma
\sqrt{y}}$ for $0\leq z \leq y$, and
\begin{align}
\widetilde v_{b}(0)  &  \leq\frac{2}{\sigma^{2}} \int_{0}^{c} y^{\frac
{\alpha-1}{2}} e^{\gamma\sqrt{y}} \left(  \int_{0}^{y} z^{-\frac{\alpha+1}{2}}
e^{-\gamma\sqrt{y}}dz\right)  dy =\frac{4c}{\sigma^{2} (1-\alpha)}%
<\infty.\nonumber
\end{align}
Then $\widetilde v_{b}(0)<\infty$ holds when $\gamma\leq0$.

When $\gamma>0$, $e^{-\gamma\sqrt{z}}> e^{-\gamma\sqrt{y}}$ for $0\leq z \leq
y$, then
\begin{align}
\widetilde v(0)  &  > \frac{2}{\sigma^{2}} \int_{0}^{c} y^{\frac{\alpha-3}{2}}
e^{\gamma\sqrt{y}} \left(  \int_{0}^{y} z^{-\frac{\alpha+1}{2}}e^{-\gamma
\sqrt{y}} dz\right)  dy =\frac{4}{\sigma^{2}(1-\alpha)} \int_{0}^{c}
y^{-1}dy=\infty.\nonumber
\end{align}
Then $\widetilde v(0)=\infty$ holds when $\gamma> 0$.

When $\gamma>0$, $e^{-\gamma\sqrt{z}}<1$ for $0\leq z \leq y$, then
\begin{align}
\widetilde v_{b}(0)  &  =\frac{2}{\sigma^{2}} \int_{0}^{c} y^{\frac{\alpha
-1}{2}} e^{\gamma\sqrt{y}} \left(  \int_{0}^{y} z^{-\frac{\alpha+1}{2}}
e^{-\gamma\sqrt{z}}dz\right)  dy\nonumber\\
&  < \frac{2}{\sigma^{2}} \int_{0}^{c} y^{\frac{\alpha-1}{2}} e^{\gamma
\sqrt{y}} \left(  \int_{0}^{y} z^{-\frac{\alpha+1}{2}} dz\right)
dy\nonumber\\
&  =\frac{2}{\sigma^{2}} \int_{0}^{c} y^{\frac{\alpha-1}{2}} e^{\gamma\sqrt
{y}} \left(  \frac{2}{1-\alpha} y^{\frac{1-\alpha}{2}} \right)  dy\nonumber\\
&  =\frac{4}{\sigma^{2} (1-\alpha)} \int_{0}^{c} e^{\gamma\sqrt{y}}
dy\nonumber\\
&  <\infty.\nonumber
\end{align}
Then $\widetilde v_{b}(0)<\infty$ holds when $\gamma>0$.

To summarize, we have that $\widetilde v(0)=\infty$ and $\widetilde
v_{b}(0)<\infty$ hold when $\alpha<1$.

Consider the case when $\alpha<1$ and $\gamma>0$. From the definition in
\eqref{vhw}
\begin{align}
\widetilde v(\infty)  &  =\frac{2}{\sigma^{2}} \int_{c}^{\infty}
y^{\frac{\alpha-3}{2}} e^{\gamma\sqrt{y}} \left(  \int_{y}^{\infty}
z^{-\frac{\alpha+1}{2}} e^{-\gamma\sqrt{z}}dz\right)  dy, \label{7hw}%
\end{align}
and
\begin{align}
\widetilde v_{b}(\infty)  &  =\frac{2}{\sigma^{2}} \int_{c}^{\infty}
y^{\frac{\alpha-1}{2}} e^{\gamma\sqrt{y}} \left(  \int_{y}^{\infty}
z^{-\frac{\alpha+1}{2}} e^{-\gamma\sqrt{z}}dz\right)  dy. \label{7bhw}%
\end{align}
Since $\gamma>0$ is assumed, then $\lim_{y\rightarrow\infty} y^{-\frac{\alpha
}{2}}e^{-\gamma\sqrt{y}}=0$, and from L'Hôpital's rule
\begin{align}
\lim_{y\rightarrow\infty} \frac{ \int_{y}^{\infty} z^{-\frac{\alpha+1}{2}}
e^{-\gamma\sqrt{z}} dz }{y^{-\frac{\alpha}{2}}e^{-\gamma\sqrt{y}}}  &
=\lim_{y\rightarrow\infty} \frac{ 1 }{\frac{\alpha}{2}y^{-\frac{1}{2}}%
+\frac{\gamma}{2}}=\frac{2}{\gamma}>0. \label{8hw}%
\end{align}
As $y\rightarrow\infty$, $\int_{y}^{\infty} z^{-\frac{\alpha+1}{2}}
e^{-\gamma\sqrt{z}}dz\sim\frac{2}{\gamma} y^{-\frac{\alpha}{2}}e^{-\gamma
\sqrt{y}}$, and there exists $M>0$, such that for $y>M$
\begin{align}
\int_{y}^{\infty} z^{-\frac{\alpha+1}{2}} e^{-\gamma\sqrt{z}}dz  &  < \frac
{4}{\gamma} y^{-\frac{\alpha}{2}}e^{-\gamma\sqrt{y}}. \label{9hw}%
\end{align}
Substitute \eqref{9hw} into \eqref{7hw}
\begin{align}
\widetilde v(\infty)  &  =\frac{2}{\sigma^{2}} \int_{c}^{M} y^{\frac{\alpha
-3}{2}} e^{\gamma\sqrt{y}} \left(  \int_{y}^{\infty} z^{-\frac{\alpha+1}{2}}
e^{-\gamma\sqrt{z}}dz\right)  dy+\frac{2}{\sigma^{2}} \int_{M}^{\infty}
y^{\frac{\alpha-3}{2}} e^{\gamma\sqrt{y}} \left(  \int_{y}^{\infty}
z^{-\frac{\alpha+1}{2}} e^{-\gamma\sqrt{z}}dz\right)  dy\nonumber\\
&  <\frac{2}{\sigma^{2}} \int_{c}^{M} y^{\frac{\alpha-3}{2}} e^{\gamma\sqrt
{y}} \left(  \int_{y}^{\infty} z^{-\frac{\alpha+1}{2}} e^{-\gamma\sqrt{z}%
}dz\right)  dy+\frac{2}{\sigma^{2}} \int_{M}^{\infty} y^{\frac{\alpha-3}{2}}
e^{\gamma\sqrt{y}} \left(  \frac{4}{\gamma} y^{-\frac{\alpha}{2}}%
e^{-\gamma\sqrt{y}}\right)  dy\nonumber\\
&  =\frac{2}{\sigma^{2}} \int_{c}^{M} y^{\frac{\alpha-3}{2}} e^{\gamma\sqrt
{y}} \left(  \int_{y}^{\infty} z^{-\frac{\alpha+1}{2}} e^{-\gamma\sqrt{z}%
}dz\right)  dy +\frac{8}{\gamma\sigma^{2}} \int_{M}^{\infty} y^{-\frac{3}{2}}
dy\nonumber\\
&  =\frac{2}{\sigma^{2}} \int_{c}^{M} y^{\frac{\alpha-3}{2}} e^{\gamma\sqrt
{y}} \left(  \int_{y}^{\infty} z^{-\frac{\alpha+1}{2}} e^{-\gamma\sqrt{z}%
}dz\right)  dy +\frac{16}{\sqrt{M}\gamma\sigma^{2}}\nonumber\\
&  <\infty.\nonumber
\end{align}
Then $\widetilde v(\infty)<\infty$, for $\alpha<1$ and $\gamma>0$.

From \eqref{8hw}, there exists $M^{\prime}>c>0$, such that for $y>M^{\prime}$
\begin{align}
\int_{y}^{\infty} z^{-\frac{\alpha+1}{2}} e^{-\gamma\sqrt{z}}dz  &  > \frac
{1}{\gamma} y^{-\frac{\alpha}{2}}e^{-\gamma\sqrt{y}}. \label{9bhw}%
\end{align}
Substitute \eqref{9bhw} into \eqref{7bhw} to obtain
\begin{align}
\widetilde v_{b}(\infty)  &  \geq\frac{2}{\sigma^{2}} \int_{M^{\prime}%
}^{\infty} y^{\frac{\alpha-1}{2}} e^{\gamma\sqrt{y}} \left(  \int_{y}^{\infty}
z^{-\frac{\alpha+1}{2}} e^{-\gamma\sqrt{z}}dz\right)  dy\nonumber\\
&  >\frac{2}{\sigma^{2}} \int_{M^{\prime}}^{\infty} y^{\frac{\alpha-1}{2}}
e^{\gamma\sqrt{y}} \left(  \frac{1}{\gamma} y^{-\frac{\alpha}{2}}%
e^{-\gamma\sqrt{y}}\right)  dy =\frac{2}{\gamma\sigma^{2}} \int_{M^{\prime}%
}^{\infty} y^{-\frac{1}{2}} dy =\infty.\nonumber
\end{align}
Then $\widetilde v_{b}(\infty)=\infty$, for $\alpha<1$ and $\gamma>0$.\qed

\subsection{Proof of Proposition \ref{uihw} \label{uihwapp}}

\proof
From the proof in Proposition \ref{hwmarthm}, we study separately the
following three cases (I), (II) and (III).

(I) $\mu>\frac{1}{2}\sigma^{2}$. Then we have the following classification:
\begin{align}
\widetilde{s}(\infty)  &
\begin{cases}
<\infty, \quad & \text{if $\gamma\geq0$},\notag\\
=\infty, \quad & \text{if $\gamma<0$},
\end{cases}
\end{align}
and $\widetilde{s}(0)=-\infty$ for $\gamma\in\mathbb{R}$. This result,
combined with the classification in Table \ref{tablehw1}, gives us the first
three rows of Table \ref{tablehw4}. From Table \ref{tablehw4} and Proposition
\ref{ui}, we have that when $\mu>\frac{1}{2}\sigma^{2}$, $(S_{t})_{0\leq
t\leq\infty}$ is not a uniformly integrable $P$-martingale.

(II) $\mu=\frac{1}{2}\sigma^{2}$. Then we have the following classification:
\begin{align}
\widetilde{s}(\infty)  &
\begin{cases}
<\infty, \quad & \text{if $\gamma> 0$},\notag\\
=\infty, \quad & \text{if $\gamma\leq0$},
\end{cases}
\end{align}
and $\widetilde{s}(0)=-\infty$ for $\gamma\in\mathbb{R}$.

This result, combined with the classification in Table \ref{tablehw1}, gives
us the three middle rows of the classification in Table \ref{tablehw4}. From
Table \ref{tablehw4} and Proposition \ref{ui}, we have that when $\mu=\frac
{1}{2}\sigma^{2}$, $(S_{t})_{0\leq t\leq\infty}$ is not a uniformly integrable
$P$-martingale.


\vspace{0.1cm} (III) $\mu<\frac{1}{2}\sigma^{2}$. Then we have the following
classification:
\begin{align}
\widetilde{s}(\infty)  &
\begin{cases}
<\infty, \quad & \text{if $\gamma> 0$},\notag\\
=\infty, \quad & \text{if $\gamma\leq0$},
\end{cases}
\end{align}
and $\widetilde{s}(0)>-\infty$ for $\gamma\in\mathbb{R}$.

This result, combined with the classification in Table \ref{tablehw1}, gives
us the last three rows of the classification in Table \ref{tablehw4}. From
Table \ref{tablehw4} and Proposition \ref{ui}, we have that when $\mu<\frac
{1}{2}\sigma^{2}$, $(S_{t})_{0\leq t\leq\infty}$ is a uniformly integrable
$P$-martingale if and only if $\gamma\leq0$, or equivalently $\rho\leq0$.
\qed



\end{document}